%% file: CoGaSu.tex
\documentclass[12pt,leqno]{amsart}
\usepackage{amsmath}
\usepackage{graphicx,color}
\usepackage[all]{xy}
\newtheorem{theorem}{Theorem}[subsection]
\newtheorem{corollary}[theorem]{Corollary}
\newtheorem{lemma}[theorem]{Lemma}
\newtheorem{proposition}[theorem]{Proposition}

\newtheorem{example}[theorem]{Example}
\newtheorem{definition}[theorem]{Definition}
\theoremstyle{remark}
\newtheorem{remark}[theorem]{Remark}
\numberwithin{equation}{section}

\newcommand{\CC}{\mathcal C}
\newcommand{\B}{\mathbb B}
\newcommand{\C}{\mathbb C}
\newcommand{\R}{\mathbb R}

\begin{document}

\title{Some aspects of analysis on almost complex manifolds with boundary}

\author{Bernard Coupet, Herv{\'e} Gaussier and Alexandre Sukhov}

\address{
\begin{tabular}{lllll}
Bernard Coupet & & Herv{\'e} Gaussier & & Alexandre Sukhov\\
C.M.I. & & I.M.J. & & U.S.T.L. \\
39, rue Joliot-Curie & & 4, place Jussieu & & Cit{\'e} Scientifique \\
13453 Marseille Cedex 13 & & 75252 Paris Cedex & & 59655 Villeneuve
d'Ascq Cedex \\
\rm coupet@cmi.univ-mrs.fr & & \rm gaussier@math.jussieu.fr & & \rm
sukhov@agat.univ-lille1.fr
\end{tabular}
}

%\email{\begin{array}{lll}
%coupet@cmi.univ-mrs.fr & & \\
%gaussier@mah.jussieu.fr & & \\
%sukhov@agat.univ-lille1.fr
%\end{array}}


\date{\number\year-\number\month-\number\day}

\begin{abstract}

We present some results dealing with the local geometry of almost
complex manifolds. We establish mainly the complete hyperbolicity of
strictly pseudoconvex domains, the extension of plurisubharmonic
functions through generic submanifolds and the elliptic regularity of
some diffeomorphisms in almost complex manifolds with boundary.
\end{abstract}

\maketitle

\tableofcontents

\section*{Introduction}
Geometric invariants are natural objects in the problem of
classifying varieties. For instance, the Riemann
curvature is such an invariant in Riemannian geometry and isometries are the
corresponding transformations. In symplectic geometry, the
transformations are called symplectomorphisms but the Darboux Theorem
states the non existence of local symplectic invariants.

In complex geometry,
Cauchy-Riemann invariants are local invariants interferring in the
classification of complex manifolds with boundaries. The generic
corresponding situation is the equivalence problem between
strictly pseudoconvex hypersurfaces in complex manifolds. Such
classifications, initiated and completely treated by {\'E}.Cartan
\cite{car} in complex dimension two, was
intensively studied and two important approaches in higher dimension are
the theory of normal forms due to
SS.Chern-J.Moser \cite{che-mos} and the Fefferman theorem \cite{fe}
connecting complex and Cauchy-Riemann geometries.

The integral representation theory, which can be considered
as the extension of the well-known Cauchy-Green formula, was
elaborated mainly by
Grauert-Lieb \cite{gra-lie} and Henkin \cite{hen} for strictly pseudoconvex domains
in the Euclidean complex space and enabled
to solve the $\bar \partial$-problem.
In contrast with complex manifolds, there is generically no
pseudoholomorphic map between almost complex
manifolds, since the holomorphicity condition is given by a
nonsolvable overdetermined system. In particular, the lack of
holomorphic coordinates prevents from developing a straightforward
integral representation theory.

The local existence of pseudoholomorphic discs in any almost complex
manifold goes back to the work of N.Nijenhuis-W.Woolf
\cite{ni-wo}. With the idea of generalizing the Newlander-Nirenberg
Theorem to integrable almost complex manifolds with low regularity, they
presented both an analytic and a geometric approach to the
$\bar{\partial}_J$-equation, satisfied by pseudoholomorphic discs in
an almost complex manifold $(M,J)$. This equation is locally a
perturbation of the standard elliptic $\bar \partial$-equation; the
local existence of pseudoholomorphic discs relies on the stability of
the solutions to elliptic systems of partial differential equations
under small quasi-linear elliptic deformation. Pseudoholomorphic discs
played a fundamental r{\^o}le in symplectic geometry where some global
symplectic invariants are associated to the moduli space of
holomorphic curves for a compatible almost complex structure. Almost
complex manifolds are viewed as natural manifolds for a deformation
theory and pseudoholomorphic discs carry some information on the
symplectic geometry by means of compactness phenomena.

The main purpose of this survey is to present some bases of
local analysis in almost complex manifolds which consists in studying, essentially
by analytic
methods, local equivalence problems between strictly pseudoconvex
domains. A large part of the survey is devoted to a precise study of
pseudoholomorphic discs and of different objects whose existence is direcly issued from
them. These are the Kobayashi-Royden metric and plurisubharmonic functions. 

Most of our approach follows the concept of structures
deformation. Nonisotropic dilations fitted to the geometry of
strictly pseudoconvex hypersurfaces were introduced in complex
analysis by S.Pinchuk \cite{pi} (related ideas were already used by Kuiper and
Benzecri in
the context of affine and projective geometry). The essential difference in almost
complex
manifolds relies on the non holomorphicity of the dilation maps. Hence
the scaling procedure involves both a deformation of domains and of
the ambiant almost complex structure. This general scheme drives to
the presentation of some exotic non integrable almost complex
structures, appearing as cluster points of dilated structures. These
model almost complex structures imply new phenomena which distinguish
almost complex manifolds from complex ones. A consequent part of our
works relies on the study of these model structures.

\vskip 0,1cm
This survey is for a large part an overview of different results obtained in a
series of
papers \cite{ga-su, CoGaSu, co-ga-su, ga-su2}. This is organized as
follows.

In Section 1, we mainly present the basic properties of almost complex
manifolds. The first two subsections are devoted to generalities. In
subsection 1.3, we explain how to attach pseudoholomorphic discs to
totally real submanifolds and we introduce the model structures. These
were introduced in \cite{ga-su}. We define complex hypersurfaces for
these structures. In Subsection 1.4 we focus on plurisubharmonic
functions. We first establish the Hopf Lemma and, as an application,
we obtain the boundary distance preserving property for
biholomorphisms between strictly pseudoconvex domains (proposition
1.4.8). Then we study the problem of removable singularities for
plurisubharmonic functions. In Subsection 1.5 we present two
constructions of almost complex structures on the tangent and on the
cotangent bundles of an almost complex manifold. These canonical lifts
play a major r{\^o}le in our studies and they will be used throughout
the survey.

Section 2 concerns the study of stationary discs. The notion of
stationary discs was introduced by L.Lempert \cite{le81} in the complex
setting. Our results can be considered as a local analogue of
Lempert's theory in almost complex manifolds. We prove the existence
of stationary discs in the unit ball for small almost complex
deformation of the standard complex structure (Propositions 2.2.1 and
2.2.2).
We show that the stationary discs form a foliation of the unit ball,
singular at the origin (Proposition 2.3.4). Then we define an analogue
of the Riemann map and we give its main properties (Theorem
2.3.10). We end Section 2 with three applications : the boundary study
of biholomorphisms (Corollary 2.3.11), a partial generalization of
Cartan's theorem (Corollary 2.3.12) and the local biholomorphic
equivalence problem between almost complex manifolds (Theorem 2.3.13).

Section 3 is devoted to the study of the Kobayashi metric in almost
complex manifolds. Proving precise estimates of the Kobayashi-Royden
metric in strictly pseudoconvex domains, similar to the estimates
obtained by I.Graham \cite{gr75} in the complex setting, we answer a
question by S.Kobayashi about the existence of basis of complete
hyperbolic neighborhoods at every point in an almost complex
manifold. The general idea is to use a boundary localization principle
(Proposition 3.1.5) and rescaling. 

In Section 4 we study the boundary behaviour of a biholomorphism
between two strictly pseudoconvex domains. We present the scaling
procedure in details and we study the properties of the dilated
objects. For a better understanding, we start with the case of four
real dimensional almost complex manifolds in Subsections 4.1 and
4.2. One of the main results concerns the behaviour of the lift of a
biholomorphism to the cotangent bundle presented in Proposition 4.2.6.
The corresponding results in higher dimension are treated in
Subsection 4.3. The analogue of Proposition 4.2.6 is Proposition
4.3.6. We end this Subsection with a compactness principle (Theorem 4.3.8). This can
be considered as an almost complex analogue to the classical
Wong-Rosay Theorem. 

Section 5 deals with the question of the regularity of
pseudoholomorphic maps attached to totally real submanifolds in an
almost complex manifold. These results are consequences of a geometric
elliptic theory. In Subsection 5.1 we show that a pseudoholomorphic
disc attached, in the sense of the cluster set, to a smooth totally
real submanifold extends smoothly up to the boundary (Theorem 5.1.1). We point out
that similar results have been established in the almost complex case
under stronger assumptions on the initial boundary regularity of the
disc. In Subsection 5.2 we establish the regularity of
a pseudoholomorphic map, defined in a wedge attached to a totally real
submanifold and whose image is contained, in the sense of the cluster
set, into a totally real submanifold of an almost complex manifold
(Proposition 5.2.1). As an application we obtain a partial version of
the Fefferman Theorem (Corollary 5.2.4).
Subsection 5.3 is devoted to the proof of the Fefferman Theorem, one
of the main results of our survey (Theorem 5.3.1). This is a
consequence of the results established in the previous sections.

\section{Basic properties of almost complex structures}

\subsection{Almost complex structures}

Everywhere in this paper $\Delta$ denotes the unit disc in $\C$ and $\B$
the unit ball in $\C^n$. Let $M$ be a smooth $\mathcal C^\infty$ real
manifold of real dimension $2n$.

\begin{definition}

$(i)$ An almost complex structure on $M$ is a smooth $\mathcal
C^\infty$-field $J$ on the tangent bundle $TM$ of $M$, satisfying
$J^2=-I$. 

$(ii)$ If $J$ is an almost complex structure on $M$ then the 2-tuple
$(M,J)$ is called an almost complex manifold.

\end{definition}
 
By an abuse of notation $J_{st}$ is the standard structure on
$\R^{2k}$ for every positive integer $k$. The standard structure
$J_{st}^{(2)}$ of $\mathbb R^2$ has the form
$$
J_{st}^{(2)} =\left(
\begin{array}{cc}
0   & -1 \\
1   & 0
\end{array}
\right).
$$
In $\mathbb R^{2n}$ with
standard coordinates
$(x^1,y^1,\dots,x^n,y^n)$
this is given by the block diagonal matrix~:
$$
J_{st}^{(2n)} = \left(
\begin{array}{cccc}
J_{st}^{(2)}   & & &   \\
 & J_{st}^{(2)} & &   \\
 &  & .  &   \\ 
 &  &   &  J_{st}^{(2)} 
\end{array}
\right).
$$
In what follows we will just write $J_{st}$ since the dimension of the space
will be clear from the context.
The first examples of almost complex manifolds are provided by complex
manifolds. Namely, a complex manifold is a smooth real manifold $M$ of
positive dimension $2n$ with local complex analytic (holomorphic)
charts $f=(f^1,\dots,f^n)$ from $M$ to $\C^n$. The almost complex
structure is locally defined by $J:= df \circ J_{st} \circ
df^{-1}$. We point out that $J$ does not depend on the choice of $f$.
Thus an almost complex manifold $(M,J)$ is a complex manifold if for every
point $p$ in $M$ there exists a neighborhood $U$ of $p$ and a coordinate
diffeomorphism $z : U \rightarrow \mathbb R^{2n}$ such that $dz \circ J \circ
dz^{-1} = J_{st}$ on $z(U)$.
 
\vskip 0,1cm
An important special case of an almost complex manifold is a bounded
 domain $D$ in $\C^n$ equipped with an almost complex structure $J$,
 defined in a neighborhood of $\bar{D}$, and sufficiently close to the
 standard structure $J_{st}$ in the standard $\mathcal C^k$ norm on
 $\bar{D}$. Every almost complex manifold may be represented locally
 in such a form. More precisely, we have the following statement. 

\begin{lemma}
\label{suplem1}
Let $(M,J)$ be an almost complex manifold. Then for every point $p \in
M$, every $\lambda_0 > 0$ and every real $k \geq 0$ there exist a neighborhood $U$
of $p$ and a
coordinate diffeomorphism $z: U \rightarrow \mathbb B$ such that
$z(p) = 0$, $dz(p) \circ J(p) \circ dz^{-1}(0) = J_{st}$  and the
direct image $z^*(J) = dz(p) \circ J(p) \circ dz^{-1}$ satisfies $\vert\vert
z^*(J) - J_{st}
\vert\vert_{\CC^k(\overline {\mathbb B})} \leq \lambda_0$.
\end{lemma}
Here $\CC^k(\overline {\mathbb B})$ denotes the standard norm on $\mathbb B$

\proof There exists a diffeomorphism $z$ from a neighborhood $U'$ of
$p \in M$ onto $\mathbb B$ satisfying $z(p) = 0$ and $dz(p) \circ J(p)
\circ dz^{-1}(0) = J_{st}$. For $\lambda > 0$ consider the dilation
$d_{\lambda}: t \mapsto \lambda^{-1}t$ in $\C^n$ and the composition
$z_{\lambda} = d_{\lambda} \circ z$. Then $\lim_{\lambda \rightarrow
0} \vert\vert z_{\lambda}^{*}(J) - J_{st} \vert\vert_{\CC^2(\overline
{\mathbb B})} = 0$. Setting $U = z^{-1}_{\lambda}(\mathbb B)$ for
$\lambda > 0$ small enough, we obtain the desired statement. \qed

\vskip 0,1cm
In particular, every almost complex structure $J$ sufficiently close 
to the standard
structure $J_{st}$ will be written locally 
$J=J_{st} + \mathcal O(\|z\|)$. 
Finally by a small perturbation 
(or deformation) of the
standard structure $J_{st}$ defined in a neighborhood of $\bar{D}$, where
$D$ is a domain in $\C^n$, we will mean a smooth one parameter family
$(J_\lambda)_\lambda$ of almost complex structures defined in a
neighborhood of $\bar{D}$, the real parameter $\lambda$ belonging to a
neighborhood of the origin, and satisfying~: $\lim_{\lambda
\rightarrow 0}\|J_\lambda - J_{st}\|_{\CC^k(\bar{D})} = 0$.

\subsubsection{$\partial_J$ and $\bar \partial_J$ operators}

Let $(M,J)$ be an almost complex manifold. We denote by $TM$ the real 
tangent bundle of $M$ and by $T_\C M:=\C \otimes TM$ its complexification.
Recall that $T_\C M = T^{(1,0)}M \oplus T^{(0,1)}M$ where
$T^{(1,0)}M:=\{ X \in T_\C M : JX=iX\} = \{\zeta -iJ \zeta, \zeta \in
TM\},$ 
and $T^{(0,1)}M:=\{ X \in T_\C M : JX=-iX\} = \{\zeta +
iJ \zeta, \zeta \in TM\}$.
 Let $T^*M$ denote the cotangent bundle of $M$.
Identifying $\C \otimes T^*M$ with
$T_\C^*M:=Hom(T_\C M,\C)$ we define the set of complex
forms of type $(1,0)$ on $M$ by~:
$
T^*_{(1,0)}M=\{w \in T_\C^* M : w(X) = 0, \forall X \in T^{(0,1)}M\}
$
and the set of complex forms of type $(0,1)$ on $M$ by~:
$
T^*_{(0,1)}M=\{w \in T_\C^* M : w(X) = 0, \forall X \in T^{(1,0)}M\}
$.
Then $T_\C^*M=T^*_{(1,0)}M \oplus T^*_{(0,1)}M$.

This allows to define the operators $\partial_J$ and
$\bar{\partial}_J$ on the space of smooth functions defined on
$M$~: given a complex smooth function $u$ on $M$, we set $\partial_J u =
du_{(1,0)} \in T^*_{(1,0)}M$ and $\bar{\partial}_Ju = du_{(0,1)}
\in T^*_{(0,1)}M$. As usual,
differential forms of any bidegree $(p,q)$ on $(M,J)$ are defined
by means of the exterior product.

We point out that there is a one-to-one correspondence between almost
complex structures and independent complex one forms. More precisely,
to any almost complex structure one can associate a basis
$(w^1,\dots,\omega^n)$of $(1,0)$-forms.  The conjugated forms
$\bar{\omega}^1,\dots,\bar{\omega}^n$ define a basis of $T^*_{(0,1)}M$
and $w^1,\dots,\omega^n,\bar{\omega}^1,\dots,\bar{\omega}^n$ are
$\mathbb C$-linearly independent.  Conversely, if $n$ complex one
forms $\omega^1,\dots,\omega^n$ on $T^*_\mathbb C M$ are such that
$w^1,\dots,\omega^n,\bar{\omega}^1,\dots,\bar{\omega}^n$ are $\mathbb
C$-linearly independent, one can define an almost complex structure on $M$ by
asserting that $(\omega^1,\dots,\omega^n)$ is a basis of $T^*_{(1,0)}M$.
The corresponding almost complex structure $J$ is defined as follows.
Set $\omega^j=\zeta^j + i \eta^j$. The one forms
$\zeta^1,\dots,\zeta^n,\eta^1,\dots,\eta^n$ define a basis of $TM$.
Then define $J$ on $TM$ by $J\zeta^j = \eta^j$ for $j=1,\dots,n$.

\subsubsection{Integrability}

Let $(X_{\bar 1},\dots,X_{\bar n})$ be a basis of $T^*_{(0,1)}M$.
If $M$ is a complex manifold, then one can find local
charts $f=(f^1,\dots,f^n)$ such that 

\begin{equation}\label{vf-eq}
X_{\bar j}f^k = 0 \ {\rm for \ every}\  j,k=1,\dots,n.
\end{equation}

Equation (\ref{vf-eq}) is equivalent to the equation $df^j(X_{\bar j})
 = 0$ for $j,k=1,\dots,n$, meaning that $(df^1,\dots,df^n)$ form a
 basis of $T^*_{(1,0)}M$. As a direct consequence of (\ref{vf-eq}) we
 have 
$[X_{\bar j},X_{\bar k}] = 0
$
for every $j,k=1,\dots,n$. This last condition is equivalent to the
integrability of the two fiber bundles $T^{(1,0)}M$ and $T^{(0,1)}M$.
One can rewrite the integrability of $T^{(1,0)}M$ as
\begin{equation}\label{nij-eq}
N_J(\zeta,\eta)=0 \ {\rm for \ every} \ \zeta,\eta \in TM
\end{equation}
where $N_J$ is the Nijenhuis tensor defined on $TM \times TM$ by~:
$$
N_J(\zeta,\eta) = [J\zeta,J\eta] - J [J\zeta,\eta] - J [\zeta,J\eta] -
  [\zeta,\eta].
$$
 
The content of the Newlander-Nirenberg theorem is the following (see
\cite{new-nir,ha77,we89})~:
\begin{theorem}\label{ne-ni-theo}
An almost complex manifold $(M,J)$ is a complex manifold if and only
 if $T^{(1,0)}M$ is integrable.
\end{theorem}

We point out that the integrability condition can also be interpreted in
terms of forms as follows : $(M,J)$ is a complex manifold if and only
if $d\omega$ has no $(0,2)$ component for every $(1,0)$ form $\omega$.

\vskip 0,2cm Let $J$ be an almost complex structure defined in a
neighborhood of the origin in $\mathbb R^{2n}$ endowed with the usual
standard complex coordinates $z^1,\dots,z^n$. We assume that $J(0) =
J_{st}$. Since $(dz^1,\dots,dz^n)$ form a basis of $(1,0)$ forms at
the origin, one can find a basis of $(1,0)$ forms
$(\omega^1,\dots,\omega^n)$ such that $\omega^j(0) = dz^j$ for
$j=1,\dots,n$.  Hence one can assume that $\omega^j = dz^j +
\sum_{k=1}A^j_j(z,\bar{z})d\bar{z}^k$, where $A^j_k$ are smooth
functions satisfying $A^j_k(0,0) = 0$. One cannot impose the condition
$(\partial A^j_k / \partial \bar{z})(0) = 0$ unless the structure $J$
is integrable at the origin. However, if ... then by the change of
variables ... one can add the normalization condition $(\partial A^j_k
/ \partial z)(0) = 0$. This normalization will be used in the local
description of strictly pseudoconvex hypersurfaces, see
Subsection~1.3.

\vskip 0,2cm
We point out that ny almost complex structure on a real  surface is integrable.

\vskip 0,1cm
Below in this first section we will describe some examples of almost complex
structures which will play essential role thoughhout the paper. These are model
almost complex structures (defined and studied in Subsection~1.3) and
canonical almost complex structures on the tangent and cotangent bundles
(Subsection~1.5). These different structures will be the center of our
study in the forecoming Sections.

\subsection{Pseudoholomorphic discs}

A smooth map $f$ between two almost complex manifolds $(M',J')$ and
$(M,J)$ is holomorphic if its differential satisfies the following
holomorphicity condition~: 

\begin{equation}\label{ho-eq1}
df \circ J = J' \circ df \  {\rm on} \ TM'.
\end{equation}

\begin{lemma}\label{ho-lem1}
The map $f$ is $(J',J)$ holomorphic if and only if
\begin{equation}\label{ho-eq2}
\forall \omega \in T^*_{(1,0)}M', f^*w \in T^*_{(1,0)}M.
\end{equation}
\end{lemma}

Here $f^*w$ is the complex one form defined on $T_\mathbb C M$ by
$f^*\omega = \omega \circ df$.

\proof
System~(\ref{ho-eq1}) implies that $df(X) \in T^{(0,1)}M$ for every
$X \in T^{(0,1)}M'$. In particular, if $\omega \in T^*_{(1,0)}M$ and
$X \in T^{(0,1)}M'$ then $(f^*\omega)(X) = \omega (df(X)) = 0$.
Conversely, assume that condition~(\ref{ho-eq2}) is satisfied. If $X \in
T^{(0,1)}M'$ and $\omega \in T^*_{(1,0)}M$, then $\omega(df(X)) = 0$.
Hence, $df(X) _in T^{(0,1)}M$ and $df(JX) = -i df(X) = J'(df(X))$.
One can prove similarly the equality $df \circ J = J' \circ df$ on
$T^{(1,0)}M'$, implying system~(\ref{ho-eq1}). \qed

\vskip 0,1cm
 Generically, if $dim_{\mathbb R} M'=2k >2$ then system~(\ref{ho-eq1})
 is overdetermined. If $k=1$ it follows from the
 previous subsection that $J'$ is integrable. In particular, one can
 view locally $(M',J')$ as  $(\Delta, J_{st})$. In case
$(M',J') = (\Delta, J_{st})$ the map $f$ is called a $J$-holomorphic
disc. We denote by $\zeta$ the complex variable in $\C$: $\zeta = x +
 iy$. Since $J_{st}(\partial / \partial x) = \partial / \partial y$
 system~(\ref{ho-eq1}) can be written~:
$$
\frac{\partial f}{\partial y} = J(f) \frac{\partial f}{\partial x},
$$
or equivalently
$$
(J+J_{st})\frac{\partial f}{\partial \bar{\zeta}} =
  (J-J_{st})\frac{\partial f}{\partial \zeta}.
$$

  In view of
Lemma~\ref{suplem1}, $J+J_{st}$ is locally invertible. Then the
holomorphicity condition is usually written as
$$
\frac{\partial f}{\partial \bar{\zeta}} + Q_J(f) \frac{\partial
f}{\partial \zeta} = 0,
$$
 where $Q_J$ is an endomorphism of $\mathbb R^{2n}$
 given by $Q_J=(J_{st}+J)^{-1}(J_{st}-J)$.
However it is easy to see that $Q_J$ is an anti $\mathbb C$-linear
endomorphism of $\mathbb C^n$. 
One can locally write a basis $w:=(w^1,\dots,w^n)$
of $(1,0)$ forms on $M$ as $w^j=dz^j +
\sum_{k=1}^nA^j_k(z,\bar{z})d\bar{z}$ where $A^j_k$ is a smooth
function satisfying the normalization conditions $A^j_k(0,0)=0$.
According to condition~(\ref{ho-eq2}) the disc $f$ being
$J$-holomorphic if $f^*(w^j)$ is a
$(1,0)$ form for $j=1,\dots,n$ (see~\cite{ch1}), meaning that
$f^*(w^j)(\partial / \partial \bar{\zeta}) = 0$, $f$ satisfies the
following equation on $\Delta$~:
\begin{equation}\label{ho-eq3}
\frac{\partial f}{\partial \bar{\zeta}} + A(f)
\overline{\frac{\partial f}{\partial \zeta}} = 0,
\end{equation}
where $A=(A_{j,k})_{1 \leq j,k \leq n}$. 
We will use equation~(\ref{ho-eq3}) to characterize the
$J$-holomorphicity in the survey.

The Nijenhuis-Woolf theorem gives the existence of pseudoholomorphic
discs and their smooth dependence on initial data.
The following general form is due to S.Ivashkovich and J.P.Rosay.
\begin{proposition}\label{I-R}
Let $k \in \mathbb N$, $k \geq 1$, and $0 < \alpha < 1$. Let $(M,J)$ be an
almost complex manifold with $J$ of class
$\mathcal C^{k,\alpha}$ and let $p \in M$. Then for every sufficiently
small $V=(v_1,\dots,v_k) \in \mathcal J^k_pM$,
there exists a $\mathcal C^{k+1,\alpha}$ $J$-holomorphic map $u_{p,V}$
from $\Delta$ into $M$ such that the $k^{th}$ jet of $u_{p,V}$ at the origin
is equal to $(p,V)$.
Moreover, $u_{p,V}$ can be chosen with $\mathcal C^1$ dependence (in
$\mathcal C^{k,\alpha}$) on the parameters $(p,V)$ in $\mathcal J^kM$.
\end{proposition}
Here $\mathcal J^kM$ denotes the space of jets of order $k$ of maps from
the unit disc $\Delta$ to $M$.

\proof We follow the proof given in~\cite{iv-ro04}.
If we fix a chart $B$ containing $p$, we can assume that
$p= 0 \in \mathbb R^{2n}$ and $J(0) = J_{st}$. Hence there exists a
neighborhood $U$ of 0 such that equation~(\ref{ho-eq3}) is defined for every
disc $u$, defined on $\Delta$, with values in $U$.

Consider the Cauchy-Green operator $T_{CG}$ for maps $g$
continuous on $\bar{\Delta}$, with values
in a complex vector space~:
$$
\forall \zeta \in \bar{\Delta}, \ T_{CG}(g)(\zeta) = 
\frac{1}{2\pi i}\int\int_{\Delta}\frac{g(\tau)}{\zeta-\tau}d\tau \bar{d\tau}.
$$
The operator $T_{CG}$ satisfies the following properties~:

(a) $g \in \mathcal C^{k,\alpha}(\bar{\Delta}) \Rightarrow T_{CG}g \in
\mathcal C^{k+1,\alpha}(\bar{\Delta})$, $ \forall k \in \mathbb N,
\ 0<\alpha<1$,

(b) $\frac{\partial}{\partial \bar{\zeta}}(T_{CG}(g)) = g$ on $\bar{\Delta}$.

By considering $A_t(u) := A(tu)$ and $u_{t,\tau}(\zeta) := t^{-1}u(\tau \zeta)$
we can assume that $\|A\|_{\mathcal C^1}$ is small enough.
%By this we mean that
%for a fixed disc $u$, the $J$-holomorphicity can be described (locally) by 
%equation~(\ref{eq1}), globally satisfied on $\Delta$ by the map $u_{t,\tau}$.
Moreover, in view of equation~(\ref{ho-eq3}), the $J$-holomorphicity of a disc
$u$ is equivalent to the (usual) holomorphicity of the disc
$h=(Id-T_{CG}(A(u)\overline{\frac{\partial}{\partial \zeta}}))u$.
Consider now the $\mathcal C^k$ mapping :
$$
\begin{array}{ccccc}
\Phi & : & (-1,1) \times \mathcal C^k(\Delta,B) & \rightarrow &
\mathcal C^k(\Delta,\mathbb C^n)\\
     &   &  (t,u)                               & \mapsto     &
(Id-T_{CG}(A(tu)\overline{\frac{\partial}{\partial \zeta}}))u.
\end{array}
$$
Since $\Phi(0,u) =u$, it follows from the Implicit Function Theorem that
there exists $0 < t_0 < 1$ such that for $|t| < t_0$ the map $\Phi(t,.)$
is a $\mathcal C^k$ diffeomorphism from a neighborhood of the origin
in $\mathcal C^k(\Delta,B)$ onto a  neighborhood of the origin $V$
in $\mathcal C^k(\Delta,\mathbb C^n)$.

For $w=(w_1,\dots,w_k) \in (\mathbb C^n)^k$ small enough, the holomorphic
map $h_{q,w}$ defined on $\bar{\Delta}$ by 
$$
h_{q,w}(\zeta) = \sum_{l=1}^k \frac{1}{l!}\zeta^l w_l
$$
belongs to $V$. If $u_{t,w}:=\Phi(t,.)^{-1}(h_{q,w})$ then $tu_{t,w}$ is
$J$-holomorphic. Moreover, since $u_{0,w} = h_w$, the $k-th$ jet of
$u_{0,w}$ at the origin is equal to $w$. hence, for sufficiently small
positive $t$ the map $w \mapsto (\frac{\partial u_{t,w}}{\partial Re(\zeta)}(0),
\dots,\frac{\partial^k u_{t,w}}{\partial (Re(\zeta)^k)}(0))$
is a diffeomorphism between neighborhoods of the origin in $(\mathbb C^n)^k$.
\qed

\vskip 0,2cm
\begin{remark}\label{nij-woo}
The statement of Proposition \ref{I-R} means that there exists a on-to-one
correpondence between
sufficiently small $J$-holomorphic discs and standard holomorphic discs.
\end{remark}

The following Proposition establishes the stability of arbitrary
$J$-holomorphic discs under perturbation of the center and of the
derivative at the center.  This result due to Ivashkovich-Rosay
in~\cite{iv-ro04} is fundamental for the upper
semi-continuity of the Kobayashi Royden pseudo-norm.

\begin{proposition}\label{kob-usc-prop}
Let $(M,J)$ be an almost complex manifold with $J$ of class $\mathcal
C^{1,\alpha}$ $(\alpha > 0)$. Let $u$ be a $J$-holomorphic map from a
neighborhood of $\bar{\Delta}$ into $M$. There exists a neighborhood
$V$ of $(u(0),\frac{\partial u}{\partial Re(\zeta)}(0))$ in $TM$ such that for
every $(q,X) \in V$, there exists a $J$-holomorphic map $v : \Delta
\rightarrow M$ with $v(0) =q$, $\frac{\partial v}{\partial Rer(\zeta)}(0) = x$.
\end{proposition}

\proof
We follow the exposition of \cite{iv-ro04}. Assume that $u$ is defined
on $\Delta_r$ for some $r > 1$. Since the map $\tilde{u}$ is $J_{st}
\times J$- holomorphic from $\Delta_r$ into $\mathbb R^2 \times M$, it
is sufficient to prove Proposition~\ref{kob-usc-prop} for imbedded
discs. Moreover, there exist $(n-1)$ smooth vector fields
$Y_1,\dots,Y_{n-1}$, defined in a neighborhood of $u(\Delta_r)$ such
that for every $\zeta \in \Delta_r$, the vectors $\frac{\partial
u}{\partial x}(\zeta),Y_1(\zeta),\dots,Y_{n-1}(\zeta)$ are
$J(u(z))$-linearly independent. Consider now the $\mathcal
C^{2,\alpha}$ change of variables $\Phi$
$$
(z_1,\dots,z_n) \mapsto u(z_1) + \sum_{j=1}^{n-1}z_{j+1}Y_j(u(z_1))
$$
defined for $|z_1| < r$, $|z_j|$ small if $j \geq 2$. The 
structure $\Phi^*(J)$ is a $\mathcal C^{1,\alpha}$ almost complex structure
that coincides with the standard structure on $\mathbb C \times \{0\}
\subset \mathbb C^n$. So Proposition~\ref{kob-usc-prop} reduces to the
following Lemma~:

\begin{lemma}\label{kob-usc-lem}
Let $J$ be a $\mathcal C^{1,\alpha}$ almost complex structure on
$\mathbb R^{2n}$ that coincides with the standard complex structure on
$\mathbb C \times \{0\}$. Let $U$ be a neighborhood of $\bar{\Delta} \times
\{0\}$. For any $(q,t) \in \mathbb C^n \times \mathbb C^n$ close enough to
$(0,0)$, there exists a $J$-holomorphic map $v : \Delta \rightarrow U$
such that $v(0) = q$ and $\frac{\partial v}{\partial Re(\zeta)}(0) =
(1,0,\dots,0) + t$.
\end{lemma}

\noindent{\bf Proof of Lemma~\ref{kob-usc-lem}}. Consider a neighborhood
of $\bar{\Delta} \times \{0\}$ on which the $J$-holomorphicity condition can
be written :
$$
\frac{\partial u}{\partial \bar{\zeta}} + A_J(u)
\overline{\frac{\partial u}{\partial \zeta}} =0
$$
with $A_J(\zeta,0,\dots,0) = 0$.

Set $\mathcal E_0=\{f : \bar{\Delta} \rightarrow \mathbb C^n /
f \in \mathcal C^{1,\alpha}, f(0) = 0, \nabla f(0) = 0\}$,
$\mathcal F_0 = \{g : \bar{\Delta} \rightarrow \mathbb C^n /
g \in \mathcal C^{\alpha}, g(0) = 0\}$ and $F(z) = (z,0,\dots,0)$.

Define the map

$$
\begin{array}{llccl}
\Psi & : & \mathcal E_0 & \rightarrow & \mathcal F_0\\
     &   &    f         & \mapsto     &  
\displaystyle{\frac{\partial(F + f)}{\partial \bar{\zeta}} + A_J(F + f)
\overline{\frac{\partial(F + f)}{\partial \zeta}}}.
\end{array}
$$

Since $\frac{\partial F}{\partial \bar{\zeta}} +0$ and
$A_J(F) = 0$ we have~:
$$\Psi(f) = \frac{\partial f}{\partial \bar{\zeta}} + B_J(f)
\overline{\frac{\partial F}{\partial \zeta}} + o(|f|),
$$
where $B_J$ is a $(2n \times 2n)$ matrix with $\mathcal C^\alpha$
entries in $z$ and $\mathbb R$-linear in $f$. We want to show that
the derivative $D\Psi_0$ of the map $\Psi$ at $f=0$ is onto.

In complex notations we can write~:

$$
D \Psi_0(f)(\zeta) = \frac{\partial f}{\partial \bar{\zeta}}(\zeta) + B_1(\zeta)
f(\zeta)
 + B_2(\zeta) \overline{f(\zeta)}
$$
where $B_1$ and $B_2$ are $(n \times n)$ complex matrices with
$\mathcal C^\alpha$ coefficients.

The surjectivity of $D\Psi_0$ follows from the following classical
result (see \cite{iv-ro04} for a direct proof)

\begin{lemma}
If $B_1$, $B_2$ are $(n \times n)$ complex matrices with $\mathcal C^\alpha$
 coefficients on $\bar{\Delta}$, for every $g \in \mathcal F_0$, there
exists $f \in \mathcal E_0$ such that
$$
\frac{\partial f}{\partial \bar{\zeta}}(\zeta) + B_1(\zeta) f(\zeta)
 + B_2(\zeta) \overline{f(\zeta)} = g(\zeta),
$$
 for every $\zeta \in \Delta$.
\end{lemma}

\qed

\subsection{Real submanifolds}
This Subsection deals with the study of totally real and of strictly
pseudoconvex submanifolds in an almost complex manifold.
In the first part, we attach Bishop's discs to a totally real
submanifold. In the second part, we describe locally strictly
pseudoconvex hypersurfaces in real dimension four as deformations
of strictly pseudoconvex hypersurfaces for the standard structure.
The associated scaling procedure is our main tool for the study of strictly
pseudoconvex domains in almost complex manifolds.

\vskip 0,2cm
Let $\Gamma$ be a real smooth submanifold in $M$ and let $p \in
\Gamma$. We denote by $H^J(\Gamma)$ the $J$-holomorphic tangent bundle
$T\Gamma \cap JT\Gamma$. 

\begin{definition}
The real submanifold $\Gamma$ is called totally real if $H^J(\Gamma)=\{0\}$
and is called $J$-complex if $H^J(\Gamma)=T\Gamma$.
\end{definition}

We note that if $\Gamma$ is a real
hypersurface in $M$ defined by $\Gamma=\{r=0\}$ and $p \in \Gamma$
then by definition $H_p^J(\Gamma) = \{v \in T_pM : dr(p)(v) =
dr(p)(J(p)v) = 0\}$.

\vskip 0,2cm
As usual, if $\theta$ is a one form on $M$ then $J^*\theta$ is the form acting
on a vector field $X$ by $(J^*\theta)X = \theta(JX)$.

We recall the notion of the Levi form of a hypersurface~:

\begin{definition}\label{DEF}
Let $\Gamma=\{r=0\}$ be a smooth real hypersurface in $M$ 
($r$ is any smooth defining function of $\Gamma$) and let $p \in \Gamma$. 

$(i)$ The {\sl Levi form} of $\Gamma$ at $p$ is the map defined on
$H^J_p(\Gamma)$  by ${\mathcal L}_{\Gamma}^J(X_p) = J^\star dr[X,JX]_p$,
where the vector field $X$ is any section of the $J$-holomorphic tangent
bundle  $H^J \Gamma$ such that $X(p) = X_p$.

$(ii)$ A real smooth hypersurface $\Gamma=\{r=0\}$ in $M$ is 
{\sl strictly $J$-pseudoconvex} at $p$ if ${\mathcal L}_{\Gamma}^J(X_p)> 0$
for any nonzero $X_p \in H^J_p(\Gamma)$. The hypersurface $\Gamma$ is called
{\sl strictly $J$-pseudoconvex} if it is strictly $J$-pseudoconvex at every point.
\end{definition}
  
\begin{remark} 
$(i)$ the ``strict $J$-pseudoconvexity'' condition does not depend on
the choice of a smooth defining function of $\Gamma$. Indeed if $\rho$
is an other smooth defining function for $\Gamma$ in a neighborhood of
$p \in \Gamma$ then there exists a positive smooth function $\lambda$
defined in a neighborhood of $p$ such that $\rho=\lambda r$.  In
particular $(J^\star dr)(p) = \lambda(p)(J^\star d\rho)(p)$.

$(ii)$ since the map $(r,J) \mapsto J^\star dr$ is smooth the ``strict
$J$-pseudoconvexity'' is stable under small perturbations of both the
hypersurface and the almost complex structure.
\end{remark}

Let $X \in TM$. It follows from the identity
$d(J^\star dr)(X,JX)=X(<J^\star dr,JX>) - JX(<J^\star dr,X>) - 
(J^\star dr)[X,JX]$
that 
$
(J^\star dr)[X,JX] = -d(J^\star dr)(X,JX)
$ 
for every $X \in H^J\Gamma$, since $<dr,JX>=<dr,JX> = 0$ in that case. 
Hence we set

\begin{definition}\label{psh-def} Let $p \in M$. If $r$ is a $\CC^2$ function on $M$
then the
Levi form of $r$ at $p$ is defined on $T_pM$ by ${\mathcal L_r}^J(p,v):=
-d(J^\star dr)_p(X,JX)$ where $X$ is any section of $TM$ 
such that $X(p) = v$.
\end{definition}

\subsubsection{Almost complex perturbation of discs}

In this subsection we attach Bishop's discs to a totally real
submanifold in an almost complex manifold.
The following statement is an almost complex
analogue of the well-known Pinchuk's construction \cite{Pi74} of a
family of holomorphic discs attached to a totally real manifold.

\begin{lemma}\label{lem-discs}
For any $\delta > 0$ there exists a family of
$J$-holomorphic discs $h(\tau,t) = h_t(\tau)$ smoothly depending on the
parameter $t \in  \R^{2n}$ such that $h_t(\partial \Delta^+) \subset E$,
$h_t(\Delta) \subset W(\Omega,E)$, $W_{\delta}(\Omega,E) \subset \cup_t
h_t(\Delta)$ and $C_1(1 - \vert \tau \vert) \leq dist (h_t(\tau),E)
\leq C_2 ( 1- \vert \tau \vert)$ for any $t$ and any $\tau \in \Delta^+$,
with constants $C_j > 0$ {\it independent} of $t$.
\end{lemma}
For $\alpha > 1$, noninteger, we denote by $\mathcal C^\alpha(\bar \Delta)$
the Banach space of functions of class $\mathcal C^\alpha$ on $\bar{\Delta}$
and by $\mathcal A^\alpha$ the Banach subspace of
$\mathcal C^\alpha(\bar \Delta)$ of functions holomorphic on $\Delta$.

First we consider the situation where $E=\{r:=(r_1,\dots,r_n)=0\}$
is a smooth totally real submanifold in $\C^n$.
Let $J_{\lambda}$ be an almost complex deformation of the standard
structure $J_{st}$ that is a one-parameter family of almost complex
structures so that $J_0 = J_{st}$. 
We recall that for $\lambda$ small enough the
$(J_{st},J_{\lambda})$-holomorphicity condition for a map $f:\Delta
\rightarrow \C^n$ may be written in the form

\begin{eqnarray}\label{equa0}
\bar\partial_{J_{\lambda}} f = \bar\partial f +
q(\lambda,f)\overline{\partial f} = 0
\end{eqnarray}
where $q$ is a smooth matrix satisfying $q(0,\cdot) \equiv 0$, uniquely
determined by $J_\lambda$ (\cite{si}).

A disc $f \in  (\mathcal C^\alpha(\bar{\Delta}))^n$ is attached
to $E$ and is $J_\lambda$-holomorphic if and only if it
satisfies the following nonlinear boundary Riemann-Hilbert type  problem~:
$$
\left\{
\begin{array}{lll}
r(f(\zeta)) = 0,& & \zeta \in \partial \Delta\\
\bar{\partial}_{J_\lambda}f(\zeta) = 0,& & \zeta \in \Delta.
\end{array}
\right.
$$
Let $f^0 \in \mathcal (\mathcal A^\alpha)^n$ be a disc attached to
$E$ and let $\mathcal U$ be a neighborhood of $(f^0,0)$ in the 
space $(\mathcal C^\alpha(\bar{\Delta}))^n \times \R$.
Given $(f,\lambda)$ in $\mathcal U$ define the maps 
$v_{f}: \zeta \in \partial 
\Delta \mapsto r(f(\zeta))$
and

\begin{eqnarray*}
& &u : \mathcal{U}  \rightarrow  (\mathcal C^\alpha(\partial \Delta))^n
\times \mathcal C^{\alpha-1}(\Delta)\\
  & & (f,\lambda)  \mapsto  (v_{f}, 
\bar{\partial}_{J_\lambda}f).
\end{eqnarray*}

Denote by $X$ the Banach space $(\mathcal C^\alpha(\bar \Delta))^n$.
Since $r$ is of class $\CC^{\infty}$, 
the map
$u$ is smooth and the tangent map $D_Xu(f^0,0)$ (we consider 
the derivative
with respect to the space $X$) is a linear map from $X$ to 
$(\CC^\alpha(\partial \Delta))^n \times \CC^{\alpha-1}(\Delta)$,
defined for every $h \in X$ by 
$$
D_Xu(f^0,0)(h) = 
\left(
\begin{array}{l}
2 Re [G h] \\
\bar\partial_{J_0} h
\end{array}
\right),$$
where for $\zeta \in \partial \Delta$
$$
G(\zeta) = \left(
\begin{array}{lll}
\frac{\partial r_1}{\partial z^1}(f^0(\zeta))  &\cdots&\frac{\partial
r_1}{\partial z^n}(f^0(\zeta))\\
\cdots&\cdots&\cdots\\
\frac{\partial r_n}{\partial z^1}(f^0(\zeta))& \cdots&\frac{\partial r_n}
{\partial z^n}(f^0(\zeta))
\end{array}
\right)$$
(see \cite{gl94}).
\begin{lemma}\label{tthh}
Assume that for some $\alpha > 1$
the linear map from $(\mathcal A^{\alpha})^n$ to
$(\mathcal C^{\alpha-1}(\Delta))^n$
given by $h \mapsto 2 Re [G h]$
is surjective and has a $d$-dimensional kernel.
Then there exist $\delta_0, \lambda_0 >0$ such that for every
$0 \leq \lambda \leq \lambda_0$,
the set of $J_\lambda$-holomorphic discs $f$ attached to $E$
and such that $\| f -f^0 \|_{\alpha} \leq \delta_0$ forms
a smooth $d$-dimensional 
submanifold 
$\mathcal A_{\lambda}$ in the Banach space 
$(C^\alpha(\bar{\Delta}))^n$.
\end{lemma} 

\noindent{\it Proof of Lemma~\ref{tthh}.}
According to the implicit function Theorem, 
the proof of Lemma~\ref{tthh} reduces to the proof of the 
surjectivity of $D_Xu$. 
It follows by classical
one-variable results on the resolution of the
$\bar\partial$-problem in the unit disc that the linear map from
$X$ to $\CC^{\alpha-1}(\Delta)$ given by 
$h \mapsto \bar \partial h$
is surjective. More precisely, given $g \in \CC^{\alpha -1}(\Delta)$
consider the Cauchy transform

$$T_{CG}(g) : \tau \in \partial \Delta \mapsto   
\frac{1}{2\pi i}
\int\int_{\Delta} \frac{g(\zeta)}{\zeta - \tau}d\zeta \wedge d\bar{\zeta}.$$

For every function $g \in 
C^{\alpha-1}(\Delta)$ the solutions $h \in X$ of the equation
$\bar\partial h = g$ have the form $h = h_0 + T_{\Delta}(g)$
where $h_0$  is an arbitrary function in $({\mathcal A}^{\alpha})^n$. 
Consider the equation 

\begin{equation}\label{EQU}
D_Xu(f^0,0)(h) = \left(
\begin{array}{l}
g_1 \\
g_2
\end{array}
\right)
\end{equation}
where $(g_1,g_2)$ is a vector-valued function with components 
$g_1 \in \CC^{\alpha-1}(\partial \Delta)$ and
$g_2 \in \CC^{\alpha-1}(\Delta)$.
Solving the $\bar\partial$-equation for the second component, we reduce 
equation~(\ref{EQU}) to    
$$
2 Re [G(\zeta) h_0(\zeta)] = g_1 - 2 Re [G(\zeta) T_{\Delta}(g_2)(\zeta)]
$$
with respect to $h_0 \in (\mathcal A^{\alpha})^n$. 
The surjectivity of the map $ h_0 \mapsto 2 Re [G h_0]$ gives the result.
\qed

\vskip 0,1cm
\noindent{\it Proof of Lemma~\ref{lem-discs}.} We proceed in three steps.

{\it Step 1. Filling the polydisc.} Consider the $n$-dimensional real 
torus $\mathbb T^n = \partial \Delta \times
...\times \partial \Delta$ in $\C^n$ and the
linear disc $f^0(\zeta) = (\zeta,...,\zeta)$, $\zeta \in \Delta$
attached to $\mathbb T^n$. 
In that case, a disc $h^0$ is in the kernel of
$h \mapsto 2 Re [G h]$ if and only if every component $h^0_k$ of $h^0$
satisfies on $\partial \Delta$ the condition $h^0_k +
\zeta^2\overline{h^0_k} = 0$. Considering the Fourier expansion
of $h_k$ on $\partial \Delta$ (recall that $h_k$ is holomorphic on
$\Delta$) and identifying the coefficients, we obtain that the map  
$h \mapsto 2 Re [G h]$ from $({\mathcal A}^{\alpha})^n$ to
$(C^{\alpha - 1}(\Delta))^n$ is surjective and has a $3n$-dimensional
kernel.
By Lemma~\ref{tthh} if $J_\lambda$ is an almost complex
structure close enough to $J_{st}$ in a neighborhood of the closure
of the polydisc $\Delta^n$, there is a $3n$-parameters family of
$J_\lambda$-holomorphic discs attached to $\mathbb T^n$. These
$J_{\lambda}$-holomorphic discs fill the intersection of 
a sufficiently small neighborhood of the point $(1,...,1)$
with $\Delta^n$.

{\it Step 2. Isotropic dilations.} Consider a smooth totally real
submanifold $E$ in an almost complex manifold $(M,J)$. Fixing local
coordinates, we may assume that $E$ is a submanifold in a neighborhood
of the origin in $\C^n$, $J = J_{st} + 0(\vert z \vert)$ and $E$ is
defined by the equations $y = \phi(x)$, where $\nabla \phi(0) =
0$. For every $\varepsilon > 0$, consider the isotropic dilations
$\Lambda_{\varepsilon}: z \mapsto z' = \varepsilon^{-1}z$. Then
$J_{\varepsilon}:= \Lambda_{\varepsilon}(J) \rightarrow J_{st}$ as
$\varepsilon \rightarrow 0$. In the $z'$-coordinates $E$ is defined by
the equations $y' = \psi(x',\varepsilon):=
\varepsilon^{-1}\phi(\varepsilon x')$ and $\psi \rightarrow 0$ as
$\varepsilon \rightarrow 0$. Consider the local diffeomorphism
$\Phi_{\varepsilon}: z' = x' +iy' \mapsto z''= x'
+i(y'-\psi(x',\varepsilon))$. Then in new coordinates (we omit the
primes) $E$ coincides with a neighborhood of the origin in $\R^n = \{
y = 0\}$ and $\hat J_{\varepsilon}: =
(\Phi_{\varepsilon})_*(J_{\varepsilon}) \rightarrow J_{st}$ as
$\varepsilon \rightarrow 0$. Assume for instance that $E=(]-1,1[
\times \{0\})^n$. For $j=1,\dots,n$, let $\Gamma_j$ be a smooth simple
curve in the real plane $\{(x_j,y_j) : x_j \leq 0\}$, containing
$]-1/2,1/2[$ and bounding a domain $G_j$. If $\psi_j$ is a
$J_{st}$-biholomorphism from $G_j$ to $\Delta$ then the map
$\psi:=(\psi_1,\dots,\psi_n)$ is a $J_{st}$-biholomorphism from $G_1
\times \cdots \times G_n$ to $\Delta^n$. Hence we may assume that $E$
is a neighborhood of the point $(1,...,1)$ on the torus $\mathbb T^n$
and the almost complex structure $J_{\varepsilon}$ is a small
deformation of the standard structure in a neighborhood of
$\bar{\Delta}^n$.  By Step 1, we may fill a neighborhood of the point
$(1,...,1)$ in the polydisc $\Delta^n$ by
$J_{\varepsilon}$-holomorphic discs (for $\varepsilon$ small enough)
which are small perturbations of the disc $\zeta \mapsto
(\zeta,...,\zeta)$.  Returning to the initial coordinates, we obtain a
family of $J$-holomorphic discs attached to $E$ along a fixed arc
(say, the upper semicircle $\partial \Delta^+$) and filling the
intersection of a neighborhood of the origin with the wedge $\{y -
\phi(x) < 0\}$.

{\it Step3.} Let now $W(\Omega,E) = \{ r_j < 0 , j=1,...,n\}$ be a
wedge with edge $E$; we assume that $0 \in E$ and $J(0) =
J_{st}$. We may assume that $E=\{y = \phi(x)\}$, $\nabla \phi(0) = 0$,
since the linear part of every $r_j$ at the origin is equal to $y_j$. So
shrinking $\Omega$ if necessary, we obtain that for any $\delta > 0$
the wedge $W_{\delta}(\Omega, E) = \{z \in \Omega: r_j(z) -
\delta\sum_{k \neq j} r_k(z) < 0 , j=1,...,n \}$ is contained in the
wedge $\{z \in \Omega: y - \phi(x) < 0 \}$. By Step 2 there is a family of
$J$-holomorphic discs attached to $E$ along the upper semicircle and
filling the wedge $W_{\delta}(\Omega,E)$. These discs are  smooth up to the
boundary and smoothly depend on the parameters. \qed

%We recall the notion of conormal bundle of a real submanifold in $\C^n$ 
%(\cite{tu01}).
%Let $\pi: T^*(\C^n) \rightarrow
%\C^n$ be the natural projection. Then $T^*_{(1,0)}(\C^n)$ can be 
%canonically identified with 
%the cotangent bundle of $\C^n$.  
%In the canonical complex coordinates $(z,t)$ on $T^*_{(1,0)}(\C^n)$
%an element of the fiber at $z$ is a holomorphic form $w=\sum_j t_j dz^j$.

%Let $N$ be a real smooth generic submanifold in $\C^n$. 
%The conormal bundle $\Sigma(N)$ of $N$ is
%a real subbundle of $T^*_{(1,0)}(\C^n)|_N$ whose fiber at $z\in N$
%is defined by $\Sigma_{z}(N) = \{ \phi \in T^*_{(1,0)}(\C^n): Re
%\,\phi \vert T_{(1,0)}(N) = 0 \}$.

%Let $\rho_1,\dots,\rho_d$ be local defining functions of $N$.
%Then the forms $\partial \rho_1,\dots, \partial \rho_d$ form
%a basis in $\Sigma_{z}(N)$ and every section $\phi$ of the bundle
%$\Sigma(N)$ has the form $\phi = \sum_{j=1}^d c_j \partial \rho_j$, 
%$c_1,\dots,c_d \in
%\R$. We will use the following
%(see\cite{tu01}):
%\begin{lemma}\label{lemma}
%Let $\Gamma$ be a real $\CC^2$ hypersurface in $\C^n$. 
%If $\Gamma$ is strictly $J$-pseudoconvex then the conormal bundle
%$\Sigma(\Gamma)$ is a totally real
%submanifold of dimension $2n$ in $T^*_{(1,0)}(\C^n)$.
%\end{lemma} 

\subsubsection{Local description of strictly pseudoconvex domains.} If
$\Gamma$ is a germ of a real hypersurface in $\C^n$ strictly
pseudoconvex with respect to $J_{st}$, then $\Gamma$ remains strictly
pseudoconvex for any almost complex structure $J$ sufficiently close
to $J_{st}$  in the $\CC^2$-norm. 
Conversely a strictly
pseudoconvex hypersurface in an almost complex manifold of real dimension
four can be represented, in suitable local coordinates, as a
strictly $J_{st}$-pseudoconvex hypersurface equipped with a small deformation
of the standard structure. Indeed, according to \cite{si} Corollary~3.1.2,
there exist a neighborhood $U$ of $q$ in $M$ and complex coordinates
$z=(z^1,z^2) : U \rightarrow  B \subset \C^2$, $z(q) =
0$ such that $z_*(J)(0) = J_{st}$ and moreover, a map $f: \Delta
\rightarrow  B$ is $J':= z_*(J)$-holomorphic if it satisfies the
equations 

\begin{eqnarray}
\label{Jhol}
\frac{\partial f^j}{\partial \bar \zeta} =
A_j(f^1,f^2)\overline{\left ( \frac{\partial f^j}
{\partial \zeta}\right ) }, j=1,2
\end{eqnarray} 
where $A_j(z) =  O(\vert
z \vert)$, $j=1,2$.

To obtain such coordinates, one can consider two transversal
foliations of the ball $\mathbb B$ by $J'$-holomorphic curves
(see~\cite{ni-wo})and then take these curves into the lines $z^j = const$
by a local diffeomorphism (see Figure 1).

The direct image of the almost complex structure
$J$ under such a diffeomorphism has a diagonal matrix $ J'(z^1,z^2) =
(a_{jk}(z))_{jk}$ with $a_{12}=a_{21}=0$ and $a_{jj}=i+\alpha_{jj}$
where $\alpha_{jj}(z)=\mathcal O(|z|)$ for $j=1,2$.
We point out that the lines $z^j = const$ are
$J$-holomorphic after a suitable parametrization (which, in general,
is not linear). 

\bigskip
\begin{center}
\input{figure4.pstex_t}
\end{center}
\bigskip

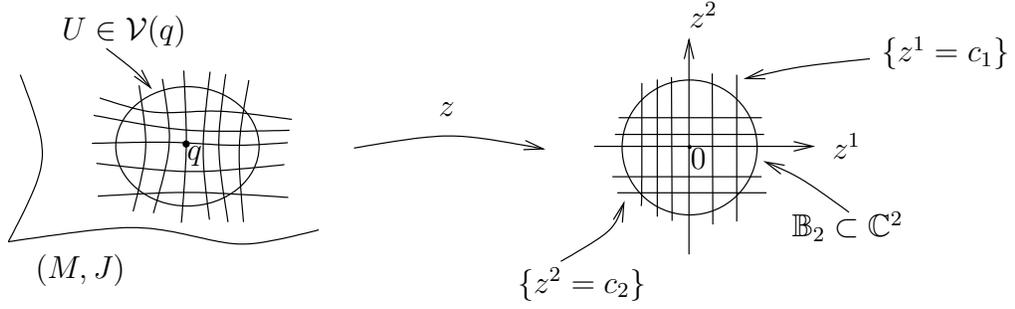
\centerline{Figure 1}
\bigskip

In what follows we omit the prime and denote this structure again by
$J$. We may assume that the complex tangent space $T_0(\partial D)
\cap J(0) T_0(\partial D) = T_0(\partial D) \cap i T_0(\partial D)$ is
given by $\{ z^2 = 0 \}$.
In particular, we have the following expansion for the defining
function $\rho$ of $D$ on $U$~:
$\rho(z,\bar{z}) = 2 Re(z^2) + 2Re K(z) + H(z) + \mathcal O(\vert z
\vert^3)$, where
$K(z)  = \sum k_{\nu\mu} z^{\nu}{z}^{\mu}$, $k_{\nu\mu} =
k_{\mu\nu}$ and 
$H(z) = \sum h_{\nu\mu} z^{\nu}\bar z^{\mu}$, $h_{\nu\mu} =
\bar h_{\mu\nu}$.

\begin{lemma}
\label{PP}
The domain $D$  is strictly $J_{st}$-pseudoconvex near the origin.
\end{lemma}

\noindent{\it Proof of Lemma~\ref{PP}.} Consider a complex vector
 $v=(v_1,0)$ tangent to $\partial D$ at the origin.  Let $f:\Delta
 \rightarrow \C^2$ be a $J$-holomorphic disc centered at the
 origin and tangent to $v$: $f(\zeta) = v\zeta + \mathcal O(\vert
 \zeta \vert^2)$.  Since $A_2 = \mathcal O(\vert z \vert)$, it follows
 from the $J$-holomorphicity equation (\ref{Jhol}) that
 $(f^2)_{\zeta\bar\zeta}(0) = 0$. This implies that $(\rho \circ
 f)_{\zeta\bar\zeta}(0) = H(v).$ Thus, the Levi form with respect to
 $J$ coincides with the Levi form with respect to $J_{st}$ on the
 complex tangent space of $\partial D$ at the origin. This proves
Lemma~\ref{PP}. \qed

\vskip 0,1cm
Consider the non-isotropic dilations $\Lambda_{\delta}: (z^1,z^2) \mapsto
(\delta^{-1/2}z^1,\delta^{-1}z^2) = (w^1,w^2)$ with $\delta > 0$. 
If $J$ has the above
diagonal form in the coordinates $(z^1,z^2)$ in $\C^2$, then
its direct image  $J_{\delta}= (\Lambda_{\delta})_*(J)$ has the form
$J_{\delta}(w^1,w^2) =(a_{jk}(\delta^{1/2}w^1,\delta w^2))_{jk}$
and so $J_{\delta}$ tends to $J_{st}$ in the $\mathcal C^2$ norm as $\delta
\rightarrow 0$. On the other hand, $\partial D$ is, in the $w$ coordinates,
the zero set of the function 
$\rho_{\delta}= \delta^{-1}(\rho \circ \Lambda_{\delta}^{-1})$.
As $\delta \rightarrow 0$, the function $\rho_{\delta}$ tends to 
the function $2 Re w^2 + 2 Re K(w^1,0) + H(w^1,0)$ which defines a
strictly $J_{st}$-pseudoconvex domain by Lemma~\ref{PP} and proves the claim.

This also proves that if $\rho$
is a local defining function of a strictly $J$-pseudoconvex domain, then
$\tilde{\rho}:=\rho + C \rho^2$
is a strictly $J$-plurisubharmonic function, quite similarly to the standard
case.

In conclusion we point out that extending $\tilde \rho$ by a suitable
negative constant, we obtain that if $D$ is a strictly
$J$-pseudoconvex domain in an almost complex
manifold, then there exists a neighborhood $U$ of $\bar{D}$ and a
function $\rho$, $J$-plurisubharmonic on $U$ and strictly
$J$-plurisubharmonic in a neighborhood of $\partial D$, such that
$D=\{ \rho <0\}$.

\subsubsection{Model almost complex structures}

The scaling process in complex manifolds deals with deformations of
domains under holomorphic transformations called dilations. The usual
nonisotropic dilations in complex manifolds, associated with strictly
pseudoconvex domains, provide the unit ball (after biholomorphism) as the
limit domain. In almost complex manifolds dilations are generically no more
holomorphic with respect to the ambiant structure. The scaling process
consists in deforming both the structure and the domain.
This provides, as limits, a quadratic domain and a linear deformation
of the standard structure in $\R^{2n}$, called {\it model structure}.
We study some invariants of such structures.
Let $(x^1,y^1,\dots,x^n,y^n)=(z^1,\dots,z^n)=('z,z^n)$ denote the canonical
coordinates of $\R^{2n}$.
\begin{definition}\label{def-model}
{\it Let $J$ be an almost complex structure on $\C^n$. We call
$J$ a {\rm model structure} if $J(z) = J_{st} + L(z)$ where $L$ is given by
a linear matrix $L=(L_{j,k})_{1 \leq j,k \leq 2n}$ such that $L_{j,k} = 0$ for
$1 \leq j \leq 2n-2, \ 1 \leq k \leq 2n$, $L_{j,k}=0$ for $j,k=2n-1,2n$ and
$L_{j,k} = \sum_{l=1}^{n-1} (a_l^{j,k} z^l + \bar{a}_l^{j,k}\bar{z}^l)$,
$a_l^{j,k} \in \C$, for $j=2n-1,2n$ and $k=1,\cdots,2n-2$.}
\end{definition}
The complexification $J_\C$ of a model structure $J$ can be written
as a $(2n \times 2n)$ complex matrix
\begin{equation}\label{complex}
J_\C=\left(
\begin{array}{ccccccc}
i      & 0      & 0      & 0      & \cdots & 0      & 0 \\
0      & -i     & 0      & 0      & \cdots & 0      & 0 \\
0      & 0      & i      & 0      & \cdots & 0      & 0 \\
0      & 0      & 0      & -i     & \cdots & 0      & 0 \\
\cdots & \cdots & \cdots & \cdots & \cdots & \cdots & \cdots \\
0      & \tilde{L}_{2n-1,2} & 0 & \tilde{L}_{2n-1,4} & \cdots & i & 0\\
\tilde{L}_{2n,1} & 0 & \tilde{L}_{2n,3} & 0 & \cdots &  0 & -i
\end{array}
\right),
\end{equation}
where $\tilde{L}_{2n-1,2k}(z,\bar{z}) = \sum_{l=1,\ l \neq k}^{n-1}
(\alpha_l^{k} z^l + \beta_l^{k} \bar{z}^l)$
with $\alpha_l^k,\ \beta_l^k \in \C$.
Moreover, $\tilde{L}_{2n,2k-1} = \overline{\tilde{L}_{2n-1,2k}}$.

\vskip 0,1cm
With a model structure we associate model domains.

\begin{definition}
Let $J$ be a model structure on $\C^n$ and $D=\{z \in \C^n : Re z^n +
P_2('z,'\bar{z})<0\}$, where $P_2$ is homogeneous second degree real
polynomial on $\C^{n-1}$.  The pair $(D,J)$ is called a {\it model
domain} if $D$ is strictly $J$-pseudoconvex in a neighborhood of the
origin.
\end{definition}

The aim of this Section is to define the complex hypersurfaces for model
structures in $\R^{2n}$.

Let $J$ be a model structure on $\R^{2n}$ and let $N$ be a germ of a
$J$-complex hypersurface in $\R^{2n}$.
\begin{proposition}\label{prop-hyp}
\hfill

$(i)$ The model structure $J$ is integrable if and only if
$ \tilde{L}_{2n-1,j}$ satisfies the compatibility conditions
$$
\frac{\partial \tilde{L}_{2n-1,2k}}{\partial \bar{z}^j} =
\frac{\partial \tilde{L}_{2n-1,2j}}{\partial \bar{z}^k}
$$
for every $1 \leq j,k \leq n-1$. 

In that case there exists a global diffeomorphism of $\R^{2n}$
which is $(J,J_{st})$ holomorphic. In that case the germs of any $J$-complex
hypersurface are given by one of the two following forms~:

\hskip 0,5cm $(a)$ $N= A \times \C$ where $A$ is a germ of a $J_{st}$-complex
hypersurface in $\C^{n-1}$,

\hskip 0,5cm $(b)$ $N =\{('z,z^n) \in \C^n : z^n = \frac{i}{4}
\sum_{j=1}^{n-1}\bar{z}^j\tilde{L}_{2n-1,2j}('z,'\bar{z})
+\frac{i}{4}
\sum_{j=1}^{n-1}\bar{z}^j \tilde{L}_{2n-1,2j}('z,0))
+ \tilde{\varphi}('z)\}$
where $\tilde{\varphi}$ is a holomorphic function locally defined in
$\C^{n-1}$.

\vskip 0,1cm
$(ii)$ If $J$ is not integrable then $N= A \times \C$ where $A$ is a germ
of a $J_{st}$-complex hypersurface in $\C^{n-1}$.
\end{proposition}
\vskip 0,2cm
\noindent{\it Proof of Proposition~\ref{prop-hyp}}. Let $N$ be a germ
of a $J$-complex hypersurface in $\R^{2n}$.
If $\pi:\R^{2n} \rightarrow \R^{2n-2}$ is the projection on the $(2n-2)$ first
variables, it follows from Definition~\ref{def-model}, or similarly
from condition~(\ref{complex}) that
$\pi(T_zN)$ is a $J_{st}$-complex hypersurface in $\C^{n-1}$.

It follows that either $dim_\C\pi(N) = n-1$ or $dim_\C\pi(N) = n-2$. 
\vskip 0,1cm
\noindent{\it Case one : $dim_\C\pi(N) = n-1$}.
We prove the following Lemma~:

\begin{lemma}\label{lem-hyp}
There is a local holomorphic function $\tilde{\varphi}$ in $\C^{n-1}$
such that
$N =\{('z,z^n) : z^n = \frac{i}{4}
\sum_{j=1}^{n-1}\bar{z}^j\tilde{L}_{2n-1,2j}('z,'\bar{z})
+\frac{i}{4}
\sum_{j=1}^{n-1}\bar{z}^j \tilde{L}_{2n-1,2j}('z,0))
+ \tilde{\varphi}('z)\}.$
\end{lemma}

{\it Proof of Lemma~\ref{lem-hyp}}. A germ $N$ can be represented as a
graph $N=\{z^n = \varphi('z,'\bar{z})\}$ where $\varphi$ is a smooth
local complex function. Hence
$T_zN=\{v_n = \sum_{j=1}^{n-1}(\frac{\partial \varphi}{\partial z^j}('z)v_j 
+ \frac{\partial \varphi}{\partial {\bar z}^j}('z) \bar{v}_j)\}$.
A vector $v=(x^1,y^1,\dots,x^n,y^n)$ belongs to $T_zN$ if and only if
the complex components $v^1:=x^1+iy^1,\dots,v^n:=x^n + i y^n$ satisfy
\begin{equation}\label{tan}
iv_n = i\sum_{j=1}^{n-1}(\frac{\partial \varphi}{\partial z^j}('z)v_j +
\frac{\partial \varphi}{\partial \bar{z}^j}('z)\bar{v}_j).
\end{equation}
Similarly, the vector $J_zv$ belongs to $T_zN$ if and only if
\begin{equation}\label{Jtan}
\sum_{j=1}^{n-1}\tilde{L}_{2n,2j-1}('z) \bar{v}_j + i v_n =
i(\sum_{j=1}^{n-1} \frac{\partial \varphi}{\partial z^j}('z) v_j
-\sum_{j=1}^{n-1}\frac{\partial \varphi}{\partial
\bar{z}^j}('z)\bar{v}_j).
\end{equation}
It follows from (\ref{tan}) and (\ref{Jtan}) that $N$ is $J$-complex if and
only if
$$
\sum_{j=1}^{n-1}(\tilde{L}_{2n,2j-1}('z)
\bar{v}_j + 2i \frac{\partial \varphi}{\partial \bar{z}^j}(z)\bar{v}_j) = 0
$$
for every $'v \in \C^{n-1}$, or equivalently if and only if
$$
\tilde{L}_{2n,2j-1} = -2i \frac{\partial \varphi}{\partial \bar{z}^j}
$$
for every $j=1,\cdots,n-1$.
This last condition is equivalent to the compatibility
conditions
\begin{equation}\label{compat}
\frac{\partial \tilde{L}_{2n,2j-1}}{\partial \bar{z}^k} =
\frac{\partial \tilde{L}_{2n,2k-1}}
{\partial \bar{z}^j}\ {\rm for}\ j,k = 1,\cdots,n-1.
\end{equation}
In that case
there exists a local holomorphic function $\tilde{\varphi}$ in $\C^{n-1}$
such that
$$
\varphi('z,'\bar{z})=\frac{i}{2}
\sum_{j=1}^{n-1}\bar{z}^j(\sum_{k \neq j}\alpha_k^j z^k)
-\frac{i}{2}
\sum_{j=1}^{n-2}\bar{z}^j(\sum_{k > j}\beta_k^j \bar{z}^k)
+ \tilde{\varphi}('z),
$$
meaning that such $J$-complex hypersurfaces are parametrized by holomorphic
functions in the variables $'z$.
Moreover we can rewrite $\varphi$ as
$$
\varphi('z,'\bar{z})=\frac{i}{4}
\sum_{j=1}^{n-1}\bar{z}^j\tilde{L}_{2n-1,2j}('z,'\bar{z})
+\frac{i}{4}
\sum_{j=1}^{n-1}\bar{z}^j \tilde{L}_{2n-1,2j}('z,0))
+ \tilde{\varphi}('z).
$$ \qed

We also have the following
\begin{lemma}\label{lem-form}
The $(1,0)$ forms of $J$ have the form
$\omega = \sum_{k=1}^n c_kdz^k -\frac{i}{2}c_n\sum_{k=1}^{n-1}
\tilde{L}_{2n-1,2k} d\bar{z}^k$ with complex numbers $c_1,\dots,c_n$.
\end{lemma}

{\it Proof of Lemma~\ref{lem-form}}.
Let $X=\sum_{k=1}^n(x^k \frac{\partial}{\partial z^k} +
y^k\frac{\partial}{\partial \bar{z}^k})$ be a $(0,1)$ vector field.
In view of (\ref{complex}), we have~:
$$
J_\C(X) = -iX \Leftrightarrow \left\{
\begin{array}{lll}
x^k & = & 0, \ \ {\rm for \ k =1,\dots,n-1}\\
 & & \\
x^n & = & \frac{i}{2} \sum_{k=1}^{n-1}y_k \tilde{L}_{2n-1,2k}.
\end{array}
\right.
$$
Hence the $(0,1)$ vector fields are given by
$$
X=\sum_{k=1}^ny^k\frac{\partial}{\partial \bar{z}^k} +
\frac{i}{2}\frac{\partial}{\partial z^n} \sum_{k=1}^{n-1}y^k\tilde{L}_{2n-1,2k}.
$$
A $(1,0)$ form $\omega=\sum_{k=1}^n(c_kdz^k + d_nd\bar{z}^k)$ 
satisfying $\omega(X)=0$ for every $(0,1)$ vector field $X$ it satisfies
$d_n = 0$ and $d_k + (i/2)c_n \tilde{L}_{2n-1,2k}=0$ for every
$k=1,\dots,n-1$.
This gives the desired form for the $(1,0)$ forms on $\C^n$. \qed
\vskip 0,1cm
Consider now the global diffeomorphism of $\C^n$ defined by
$$
F('z,z^n) = ('z,z^n - \frac{i}{4}
\sum_{j=1}^{n-1}\bar{z}^j\tilde{L}_{2n-1,2j}('z,'\bar{z})
-\frac{i}{4}
\sum_{j=1}^{n-1}\bar{z}^j \tilde{L}_{2n-1,2j}('z,0)).
$$
The map $F$ is $(J,J_{st})$ holomorphic if and only if $F^*(dz^k)$
is a $(1,0)$ form with respect to $J$, for every $k=1,\dots,n$.

Then $F^*(dz^k) = dz^k$ for $k=1,\dots,n-1$ and
$$
\begin{array}{lll} 
F^*(dz^n) & = & \displaystyle
dz^n + \sum_{k=1}^{n-1} \frac{\partial F_n}{\partial z^k} dz^k
+ \sum_{k=1}^{n-1} \frac{\partial F_n}{\partial \bar{z}^k} d\bar{z}^k\\
& & \\
 & = & \displaystyle
dz^n + \sum_{k=1}^{n-1} \frac{\partial F_n}{\partial z^k} dz^k\\
& & \ \ \ \ \ - \displaystyle \frac{i}{4} \sum_{k=1}^{n-1}
(\tilde{L}_{2n-1,2k}('z,'\bar{z})+
\sum_{j \neq k}\bar{z}^j
\frac{\partial \tilde{L}_{2n-1,2j}}{\partial \bar{z}^k}('z,'\bar{z})
+ \tilde{L}_{2n-1,2k}('z,'0) d\bar{z}^k.
\end{array}$$
By the compatibility condition~(\ref{compat}) we have
$$
\begin{array}{lll}
F^*(dz^n) & = & \displaystyle
dz^n + \sum_{k=1}^{n-1} \frac{\partial F_n}{\partial z^k} dz^k\\
& & \\
 & & -\frac{i}{4} \sum_{k=1}^{n-1}
(\tilde{L}_{2n-1,2k}('z,'\bar{z}) +
\tilde{L}_{2n-1,2k}('0,'\bar{z}) +
\tilde{L}_{2n-1,2k}('z,'0)) d\bar{z}^k\\
& & \\
 & = & \displaystyle dz^n 
- \frac{i}{2} \sum_{k=1}^{n-1}
\tilde{L}_{2n-1,2k}('z,'\bar{z}) d\bar{z}^k
+ \sum_{k=1}^{n-1} \frac{\partial F_n}{\partial z^k} dz^k.
\end{array}
$$
These equalities mean that $F$ is a local $(J,J_{st})$-biholomorphism of
$\C^n$, and so that $J$ is integrable.

\vskip 0,1cm
\noindent{\it Case two : $dim_\C\pi(N) = n-2$}. In that case we may write
$N=\pi(N) \times \C$, meaning that $J$-complex hypersurfaces
are parametrized by $J_{st}$-complex hypersurfaces of $\C^{n-1}$.

\vskip 0,1cm
We can conclude now the proof of Proposition~\ref{prop-hyp}. We proved in
Case one that if there exists a $J$-complex hypersurface in $\C^n$ such
that $dim \pi(N) = n-1$ (this is equivalent to the compatibility conditions
(\ref{compat})) then $J$ is integrable.
Conversely, it is immediate that if $J$
is integrable then there exists a $J$-complex hypersurface whose form
is given by Lemma~\ref{lem-hyp} and hence that the compatibility conditions
(\ref{compat}) are satisfied. This gives part $(i)$ of
Proposition~\ref{prop-hyp}.

To prove part $(ii)$, we note that if $J$ is not integrable then in view
of part $(i)$ the form of any $J$-complex hypersurface is given by Case two.
 \qed

\subsection{Plurisubharmonic functions}

In this Subsection, we present essentially two results. The first establishes
the Hopf Lemma. As a consequence, we obtain the boundary equivalence property
for biholomorphisms between relatively compact, strictly pseudoconvex domains.
The second result deals with the removability of singularities
for plurisubharmonic functions.

\vskip 0,2cm
We first recall the following 
definition~:
\begin{definition}\label{d6}
An upper semicontinuous function $u$ on $(M,J)$ is called 
{\sl $J$-plurisubharmonic} on $M$ if the composition $u \circ f$ 
is subharmonic on $\Delta$ for every $J$-holomorphic disc $f:\Delta
\rightarrow M$.
\end{definition} 
If $M$ is a domain in $\C^n$ and $J=J_{st}$ then a 
$J_{st}$-plurisubharmonic function is a plurisubharmonic function 
in the usual sense. 

\vskip 0,1cm
\begin{proposition}\label{di-su}
Let $r$ be a real function of class $C^2$ in a neighborhood of a point $p \in M$.
\begin{itemize}
\item[(i)] If $F: (M,J) \longrightarrow (M', J')$ is a $(J,
  J')$-holomorphic map,  and $\varphi$ is a real function of class
  $C^2$ in a neighborhood of $F(p)$, then for any $v \in T_p(M)$ we have
$L^J_{\varphi \circ F}(p;v) = L^{J'}_\varphi(F(p),dF(p)(v))$.
\item[(ii)] If $z:\mathbb D \longrightarrow M$ is a $J$-holomorphic disc satisfying
  $z(0) = p$, and $dz(0)(e_1) = v$ (here $e_1$ denote the vector
  $\frac{\partial}{\partial Re(\zeta)}$ in
  $\R^2$), then $L^J_r(p;v) = \Delta (r \circ z) (0)$.
\end{itemize}
\end{proposition}
 The property (i) expresses the invariance of the Levi form with
respect to biholomorphic maps. The property (ii) is often useful in order to compute
the Levi form if a vector
$t$ is given. 
\proof (i) Since the map $F$ is $(J, J')$-holomorphic, we have ${J'}^* dr(dF(X))
= dr(J' dF(X)) = dr ( dF(J X)) = d(r \circ F)(JX)$ that is
$F^*({J'}^* dr) = J^* d(r \circ F)$. By the invariance of the exterior
derivative we obtain that $F^*(d{ J'}^* dr) = d J^* d (r \circ F)$. Again
using the holomorphy of $F$, we get $d{ J'}^* dr(dF(X),J'dF(X)) = 
F^*(d{ J'}^* dr)(X,JX) = d J^* d(r \circ F)(X,JX)$ which
implies (i).

(ii) Since $z$ is a
$(J_{st},J)$-holomorphic map, (i) implies that $L_r^J(p,v) = L_{r \circ
  z}^{J_{st}}(0,e_1) = \Delta(r \circ z)(0)$. This proves proposition. \qed

It follows from Proposition~\ref{di-su} that a $\CC^2$ real valued
function $u$ on $M$ is
$J$-plurisubharmonic on $M$ if and only if $\mathcal L^J(u)(p)(v) \geq 0$
for every $p \in M$, $v \in T_pM$.
\vskip 0,1cm
This leads to the definition~:
\begin{definition}
A $\CC^2$ real valued function $u$ on $M$ is {\sl strictly 
$J$-plurisubharmonic} on $M$ if  $\mathcal L^J(u)(p)(v)$
is positive for every $p \in M$, $v \in T_pM \backslash \{0\}$.
\end{definition}

We have the following example of a
$J$-plurisubharmonic function on an almost complex manifold $(M,J)$~:
\begin{example}\label{example}
For every point $p\in (M,J)$ there exists a neighborhood $U$ of $p$
 and a diffeomorphism $z:U \rightarrow \mathbb
B$ centered at $p$ (ie $z(p) =0$) such that the function $|z|^2$ is
$J$-plurisubharmonic on $U$. 
\end{example}

\vskip 0,1cm     
We also have the following 

\begin{lemma}
A function  $u$ of class $\CC^2$ in a neighborhood
of a point $p$ of $(M,J)$  is strictly $J$-plurisubharmonic
if and only there exists a neighborhood $U$ of $p$  with local
complex coordinates $z:U \rightarrow \mathbb B$ centered at $p$, such
that the function $u - c|z|^2$ is $J$-plurisubharmonic on $U$ for some
constant $c > 0$.
\end{lemma}

The function $\log \vert z \vert$ is
$J_{st}$-plurisubharmonic on $\C^n$ and plays an important role in the
pluripotential theory as the Green function for the complex
Monge-Amp{\`e}re operator on the unit ball. In particular, this function
is crucially used in Sibony's method  in order to localize and
estimate the Kobayashi-Royden metric on a complex manifold. Unfortunately,
after an arbirarily small general almost complex deformation of the
standard structure this function is {\it not} plurisubharmonic with
respect to the new structure (in any neighborhood of the origin), see
for instance \cite{de99}. So we will need the following statement
communicated to the authors by E.Chirka~:
\begin{lemma}
Let $p$ be a point in an almost complex manifold $(M,J)$. There exist
a neighborhood $U$ of $p$ in $M$, a diffeomorphism $z : U \rightarrow
\B$ centered at $p$ and positive constants $\lambda_0,\ A$, such that
the function $\log|z| + A|z|$ is
$J'$-plurisubharmonic on $U$ for every almost complex structure $J'$
satisfying $\|J'-J\|_{\CC^2(\bar{U})} \leq \lambda_0$.
\end{lemma}

\noindent{\it Proof.} Consider the function $u=|z|$ on $\mathbb B$.
Since $\mathcal L^{J_{st}}(u \circ z^{-1})(w)(v) \geq \|v\|^2/4|w|$
for every $w \in \B\backslash \{0\}$ and every $v \in \C^n$, it
follows by a direct expansion of ${\mathcal L}^{J'}(u)$ that there
exist a neighborhood $U$ of $p$, $U \subset\subset U_0$, and a
positive constant $\lambda_0$ such that ${\mathcal L}^{J'}(u)(q)(v)
\geq \|v\|^2/5|z(q)|$ for every $q \in U \backslash\{p\}$, every $v\in
T_qM$ and every almost complex structure $J'$ satisfying
$\|J'-J\|_{\CC^2(\bar{U})} \leq \lambda_0$. Moreover, computing the
Laplacian of $\log|f|$ where $f$ is any $J$-holomorphic disc we
obtain, decreasing $\lambda_0$ if necessary, that there exists a
positive constant $B$ such that $\mathcal L^{J'}(\log|z|)(q)(v) \geq
-B\|v\|^2/|z(q)|$ for every $q \in U\backslash\{p\}$, every $v \in
T_qM$ and every almost complex structure $J'$ satisfying
$\|J'-J\|_{\CC^2(\bar{U})} \leq \lambda_0$. We may choose $A = 2B$ to
get the result. \qed

\vskip 0,2cm
We point out that such constructions of plurisubharmonic functions were
generalized recently by J.P.Rosay. He constructed a plurisubharmonic function
in $\mathbb C^2$, whose polar set is a two real dimensional submanifold in
$\mathbb C^2$.

\subsubsection{Hopf lemma and the  boundary distance preserving property}

In what follows we need an analogue to the Hopf lemma 
for almost complex manifolds. This can be proved quite similarly 
to the standard one. 

\begin{lemma}\label{hopf}
(Hopf lemma) Let $G$ be a relatively compact domain with a $\mathcal C^2$
boundary on an almost complex manifold $(M,J)$. Then for any negative
$J$-psh function $u$ on $D$ there exists a constant $C > 0$ such that
$\vert u(p) \vert \geq C dist(p,\partial G)$ for any $p \in G$ ($dist$
is taken with respect to a Riemannian metric on $M$).
\end{lemma}
\noindent{\it Proof of Lemma~\ref{hopf}}.
{\it Step 1.} We have the following precise version on the unit
disc: let $u$ be a subharmonic function on $\Delta$, $K$ be a fixed
compact on $\Delta$. Suppose that $u < 0$ on $\Delta$ and $u \vert K
\leq -L$ where $ L > 0$ is constant. Then there exists $C(K,L) > 0$
(independent of $u$) such that $\vert u(p) \vert \geq C
dist(p,\partial \Delta)$ (see \cite{Ra}). 

{\it Step 2.} Let $G$ be a domain in $\C$ with $\mathcal C^2$-boundary.
Then there exists an $r > 0$ (depending on the curvature of the boundary)
such that for any boundary point $q \in \partial G$ the ball $B_{q,r}$
of radius $r$ centered on
the interior normal to $\partial G$ at $q$, such that  $q \in
\partial B_{q,r}$, is
contained in $G$. Applying Step 1 to the restriction of $u$ on every
such a ball (when $q$ runs over $\partial G$) we obtain the Hopf lemma
for a domain with $\CC^2$ boundary: 
let $u$ be a subharmonic function on $G$, $K$ be a fixed
compact on $G$. Suppose that $u < 0$ on $G$ and $u \vert K
\leq -L$ where $ L > 0$ is constant. Denote by $k$ the curvature of
$\partial G$. Then there exists $C(K,L,k) > 0$
(independent of $u$) such that $\vert u(p) \vert \geq C
dist(p,\partial \Delta)$.

{\it Step 3.} Now we can prove the Hopf lemma for almost complex
manifolds. Fix a normal field $v$ on $\partial G$ and consider the family
of $J$-holomorphic discs $d_v$ satisfying $d'_0(\partial_x) =
v(d(0))$. The image of such a disc is a real surfaces intesecting
$\partial G$ transversally, so its pullback gives a $\mathcal C^2$-curve in
$\Delta$. Denote by $G_v$ the component of  $\Delta$ defined by the
condition $d_v(G_v) \subset G$. Then every $G_v$ is a domain with
$\mathcal C^2$-boundary in $\C$ and the curvatures of boundaries depend
continuously on $v$. We conclude by applying Step 2 to the composition
$u \circ d_v$ on $G_v$. \qed

\vskip 0,2cm
As an application, we obtain the boundary distance preserving property for
biholomorphisms between strictly pseudoconvex domains. 

\begin{proposition}\label{equiv}
Let $D $ and $D'$ be two  smoothly bounded strictly pseudoconvex
domains in almost complex manifolds $(M,J)$ and
$(M',J')$ respectively and let $f:D \rightarrow D'$ be a
$(J,J')$-biholomorphism. Then 
there exists a constant $C > 0$ such that

$$
(1/C) dist(f(z),\partial D') \leq dist(z,\partial D) \leq C dist
(f(z),\partial D').
$$
\end{proposition}
\noindent{\it Proof of Proposition~\ref{equiv}}.
According to the previous section, we may assume that  $D = \{ p:
\rho(p) < 0 \}$ where $\rho$ is a  $J$-plurisubharmonic
function on $D$, strictly $J$-plurisubharmonic in a neighborhood of the
boundary; similarly $D'$ can be defined by means of a function
$\rho'$. Now it suffices
to apply the Hopf lemma to the functions $\rho' \circ f$ and $\rho
\circ f^{-1}$. \qed

\subsubsection{Removable singularities for plurisubharmonic functions}
In this subsection we prove the following

\begin{theorem}
\label{rem-theo}
Let $(M,J)$ be an almost complex manifold of real dimension four and let
$E \subset M$ be a generic submanifold of $M$ of real
codimension two. Then for any continuous plurisubharmonic function $u$ on
$M \backslash E$ the function $u^*$ defined by $u^*(x) = u(x)$ for $x \in
M \backslash E$ and $u^*(x) = \lim \sup_{y \in M \backslash E, y
\longrightarrow x} u(x)$ for $x \in E$, is $J$-plurisubharmonic on $M$.
\end{theorem}

As usual, by a generic manifold we mean a real
submanifold $E$ of $(M,J)$ such that the tangent space of $E$ at
every point $p \in E$ spans the tangent space $T_p(M)$ (considered as a
complex space with the structure $J(p)$). In real dimension four, $E$ is
generic if and only if it is totally real (for any vector $v \in T_p(E)
\backslash \{ 0 \}$ the vector $J(p)v$ is not in $T_p(E)$).

Theorem~\ref{rem-theo} deserves some comments.

Firstly, the definition of $u^*$ in Theorem~\ref{rem-theo} is correct since $E$
has the empty interior. The function $u^*$ is the unique possible plurisubharmonic
extension of $u$ on $M$.

Secondly, Theorem~\ref{rem-theo} is obtained by N.Karpova \cite{Ka} in the case
where $M$ is a complex manifold. Related results in the
integrable case are due to B.Shiffman \cite{Sh}, P.Pflug \cite{Pf},
U.Cegrell \cite{Ce}.

Thirdly, our proof admits easy generalizations to higher
dimensions. These are given at the end of this Subsection.

%{\bf Remark 5.} The assumption of continuity of $u$ is technical; we
%believe that it may be removed in the statement of Theorem~\ref{maintheo}.

\vskip 0,1cm
Our method is similar to the approach of \cite{Ka} and
includes two main steps. First, we show that $u$ is upper
bounded. The proof uses a filling of $M \backslash E$ by $J$-holomorphic
discs. This shows that $u^*$ is defined correctly. In order to
prove its plurisubharmonicity, we can not use the characterization of
upper semi-continuous $J$-plurisubharmonic functions in term of positivity of their
Levi current since this
result is not yet established for almost complex structures. So we use another
approach based on the construction of the plurisubharmonic envelope of an
upper semicontinuous function by means of $J$-holomorphic curves; this
construction is quite elementary and in the case of the standard structure
this is due to Edgar \cite{Ed} and Bu-Schachermayer
\cite{BuSch}; from our point of view, it is of independent interest.

\vskip 0,2cm
\noindent{\bf Step 1. Filling by J-holomorphic discs and upper bound
for $u^*$}

\vskip 0,1cm
The statement of Theorem~\ref{rem-theo} is local. So everywhere below we
may assume $M$ is the unit ball $\B_2 \subset \mathbb C^2$, $J = J_{st} +
O(\vert z \vert)$ where $z = x + iy$ are standard coordinates in $\mathbb
C^2$ and $J_{st}$ is the standard complex structure on $\mathbb C^2$. We
also may assume that $E$ coincides with $\R^2 = \{z \in \mathbb C^2 : y =
0 \}$.

%\subsubsection{Wedges and support hypersurfaces}
\vskip 0,1cm
Consider the
``wedge-type'' domains $M_{11} = \{ y_1 > 0, y_2 > 0 \}$, $M_{12} = \{ y_1
> 0, y_2 < 0 \}$, $M_{21} = \{ y_1 < 0, y_2 > 0 \}$, $M_{22} = \{ y_1 < 0,
y_2 < 0 \}$. Then $\B_2 \backslash \R^ 2 = \cup M_{ij}$ so it is enough to
show that $u$ is upper bounded on every
$M_{ij}$. We prove it for instance on $M_{11}$.

Consider the ``support'' hypersurface $\Gamma = \{ \rho = 0 \}$ where
$\rho = y_1 + y_2 + y_1^2 +y_2^2$ so that $M_{11} \subset
\Gamma^+ = \{ \rho > 0 \}$ and $\overline M_{11} \cap \Gamma = \R^2$. We
may assume that the norm $\| J - J_{st} \|_{C^2}$ is small enough
so that the function $\rho$ is strictly $J$-plurisubharmonic on $\B_2$.

For $t \geq 0$ consider the translated hypersurface $\Gamma_t = \{ z:
\rho(z) - t > 0 \}$ which is strictly $J$-pseudoconvex for small $t$.
 According to the result of J.F.Barraud-E.Mazzilli \cite{ba-ma} and
 S.Ivashkovich-J.P.Rosay \cite{iv-ro04} there exists a family  $(f_{p,t} :
 \Delta \rightarrow \B_2$ of
$J$-holomorphic discs with the following properties~(see Figure 2) :
\begin{itemize}
\item[(a)] the discs depend smoothly on parameters $p \in \Gamma_t$, $t \geq
0$ 
\item[(b)] there exists an $r > 0$ such that for every $p$, $t$ the
surface $f_{p,t}(r\Delta)$ is contained in $\{ p \} \cup \{ z:\rho(z) - t > 0
\}$ and $f_{p,t}(0) = p$.
\end{itemize}

\bigskip
\begin{center}
\input{figure1.pstex_t}
\end{center}
\bigskip

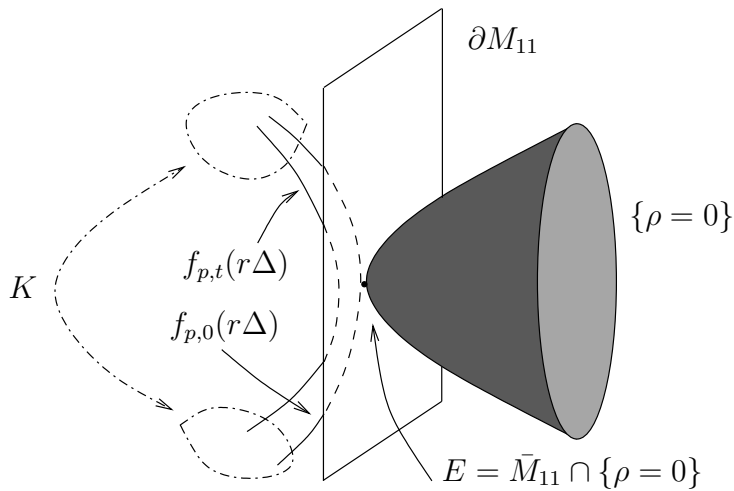
\centerline{Figure 2}

\bigskip
This family of discs fills the ``wedge''
$M_{11}$ and every disc is contained in $\B_2 \backslash \R^2$. It follows
by construction that their boundaries $f(r\partial \Delta)$ form a compact 
set
$K = \cup_{p,t}f(r\partial \Delta)$ in $\B_2 \backslash
\R^2$. Since $u$ is bounded above on $K$, say, $u \vert K \leq C$,
applying the maximum principle for every subharmonic function $\rho \circ
f_{pt}$ we obtain that $u$ is bounded from above by $C$ on $M_{11}$. Repeating
this construction for other wedges we obtain the following

\begin{proposition}
\label{pro1}
The function $u$ is upper bounded on $\B_2 \backslash \R^2$.
\end{proposition}

The standard argument allows now to reduce the proof of
Theorem~\ref{rem-theo} to the case where $u$ is bounded. Indeed,
suppose that Theorem~\ref{rem-theo} is proved for bounded
functions. For any positive integer $n$, consider the continuous
$J$-plurisubharmonic function $u_n = \sup (u, - n)$. So $u_n^*$ is
$J$-plurisubharmonic on $\B_2$. Since the sequence $(u_n)$ is
decreasing, it converges to a $J$-plurisubharmonic function $\hat u$
on $\B_2$ and $u = \hat u$ on $\B_2 \backslash \R^2$. We prove that $u^* =
\hat u$. Clearly, by the uppersemicontinuity of $\hat{u}$ we have $u^*
\leq \hat u$. Fix now a point $x_0 \in \R^2$ and a vector $v \in
T_{x_0}(\B_2)$ generating a complex line $L$ such that $L \cap
T_{x_0}\R^2 = \{ 0 \}$. According to Nijenhuis-Woolf~\cite{ni-wo}, there
exists a $J$-holomorphic disc $f$ such that $f(0) = x_0$ and $L$
is tangent to the surface $f(\Delta)$ at $x_0$. So the bounded function $u
\circ f$ is subharmonic on the punctured disc $\Delta^*$ and so extends
uniquely as a subharmonic function on $\Delta$ setting $(u\circ f)(0) = \lim
\sup_{\zeta \longrightarrow 0 } u(\zeta)$. Since this extension is unique,
we obtain that $\hat u(x_0) = (\hat u \circ f)(0) = \lim \sup_{\zeta
\longrightarrow 0 } (u \circ f)(\zeta) \leq \lim \sup_{z \longrightarrow
x_0} u(z) = u^*(x_0)$.  Applying this argument to every point in $\mathbb R^2$
we obtain that $\hat u = u^*$ on
$\B_2$.
Thus, {\it everywhere below we assume that $u$ is bounded on $\B_2
\backslash \R^2$.}

\vskip 0,2cm
\noindent{\bf Step two. $J$-plurisubharmonic envelopes of upper semicontinuous
functions}

\vskip 0,1cm
Denote by $\mu_r$ the normalized Lebesgue measure of the disc $r\Delta$ (we
simply write $\mu$ if $r = 1$).

\begin{proposition}
\label{pro2.1}
Let $v$ be an upper semicontinuous function on
$\B_2$. Consider the sequence $(v_n)$ defined as follows: $v_0 = v$ and
for $n \geq1$, for $z \in \B^2$,
$$
v_n(z) = \inf \int_{\Delta} (v_{n-1} \circ f)
(\zeta)d\mu,
$$
where $\inf$ is taken over all $J$-holomorphic discs $f : \Delta \rightarrow
\mathbb B^2$ such that $f(0) = z$ (without loss of generality we always
assume that every disc $f$ is continuous on $\overline \Delta$ and $f(\overline
\Delta) \subset \B_2$). Then the sequence $(v_n)$
decreases pointwise to the largest $J$-plurisubharmonic function $\hat v$ on $\B_2$
bounded from above by $v$.
\end{proposition}

\noindent{\it Proof of Proposition~\ref{pro2.1}.} {\it Step 1.} The sequence $(v_n)$
decreases.
Indeed, for every $z$, the constant disc $f_z(\zeta) \equiv z$ is
$J$-holomorphic so
$$v_n(z) = \inf_f \int_{\Delta} (v_{n-1} \circ f)d\mu \leq
\int_{\Delta} v_{n-1} \circ f_z d\mu = v_{n-1}(z).
$$

\begin{lemma}\label{usc}
The function $\hat v$ is upper semicontinuous.
\end{lemma}

\noindent{\it Proof of Lemma~\ref{usc}.}
We proceed by induction on $n$. For $n = 0$ the statement is
correct. Suppose that $v_{n-1}$ is upper semicontinuous. Let $(z_k) \subset \B_2$ be a
sequence of points converging to $z_0 \in \B_2$.

The following claim is a direct consequence of Proposition \ref{I-R}.

{\bf Claim.} {\it Let $f:\Delta \rightarrow \B_2$ be a $J$-holomorphic disc
  continuous on $\overline \Delta$ such that $f(0) = z_0$ and $f(\overline \Delta)
\subset \B_2$. Then there
exists a sequence of $J$-holom{\`u}orphic discs $f_k : \Delta \rightarrow \B_2$,
continuous on $\overline \Delta$ such that $f_k(\overline \Delta) \subset \B_2$,
$f_k(0) = z_k$
for every $k$ and $f_k \longrightarrow f$ uniformly on $\Delta$.}

\vskip 0,1cm
Consider the compact set $K$ defined as the closure of the union $\cup_k
f_k(\overline \Delta)$. Since $v_{n-1}$ is upper semicontinuous it is bounded from
above by a
constant $C$ on $K$ and
$$
(v_{n-1}\circ f)(\zeta) \geq \lim \sup_{k \longrightarrow \infty}
(v_{n-1} \circ f_k)(\zeta), \zeta \in \Delta.
$$

So by the Fatou lemma
$$
\int_{\Delta} (v_{n-1} \circ f d\mu \geq \lim \sup_{k
\longrightarrow \infty} \int_{\Delta} (v_{n-1} \circ f_k) d\mu \geq \lim
\sup_{k \longrightarrow \infty} v_n(z_k).
$$
This implies that
$$
v_n(z_0) = \inf_f \int_{\Delta} v_{n-1} \circ f d\mu \geq \lim \sup_{k
\longrightarrow \infty} v_n(z_k)
$$
and proves Lemma~\ref{usc}.

Therefore, the function $\hat v$ is also upper semicontinuous as a decreasing
limit of upper semicontinuous functions.

{\it Step 2.} We prove by induction that for any $J$-plurisubharmonic function
$\phi$ satisfying $\phi
\leq v$ we have $ \phi \leq v_n$ for any $n$.  This is true for $n = 0$.
Suppose that $\phi \leq v_{n-1}$. Fix an arbitrary point $z_0 \in \B_2$.
For every $J$-holomorphic disc $f$ continuous on $\overline \Delta$ and satisfying
$f(0) =z_0$,
$f(\overline \Delta) \subset \B_2$ we have~:
$$
\phi(z_0) \leq
\int_{\Delta} \phi \circ f d\mu \leq \int_{\Delta} v_{n-1} \circ f d\mu.
$$
  So
$v_n(z_0) \leq \phi(z_0)$.

{\it Step 3.} We show that the restriction of $\hat v$ on a
$J$-holomorphic disc is subharmonic. Given $z_0$ and $f$ as above, we have
$$
\hat v(z_0) = \lim_{n \longrightarrow \infty} v_n(z_0) \leq \lim_{n
\longrightarrow \infty} \int_{\Delta} v_{n-1} \circ f d\mu = \int_{\Delta} v \circ f
d\mu
$$
by the Beppo Levi theorem. This proves Proposition~\ref{pro2.1}. \qed

\vskip 0,1cm
Let now $u$ be as in the previous subsection (that is a bounded continuous
$J$-plurisubharmonic function on $\B_2 \backslash \R^2$) and let $u^*$ be its upper
semicontinuous extension defined in the statement of
Theorem~\ref{rem-theo}. Setting $v = u^*$, we apply Proposition
\ref{pro2.1} and obtain the largest $J$-plurisubharmonic function $\hat v$ on $\B_2$
bounded from above by $v$. In order to conclude the proof of the
Theorem~\ref{rem-theo}, it is sufficient to show that $\hat v =
u^*$ on $\B_2$.

The substantial part of the proof is contained in the following

\begin{proposition}
\label{prop2.2}
For any $z \in \B_2 \backslash \R^2$ we have $\hat v(z) = u(z)$.
\end{proposition}
\noindent{\bf Proof.}
Since $\hat v \leq u$ ($u=u^*$ on $\B_2 \backslash \R^2$),
we just need to prove the inverse
inequality. Using the induction, we show that for every $n$, $v_n \geq
u$. For $n=0$ this is clear. Suppose that $v_{n-1} \geq u$ on $\B_2
\backslash \R^2$ (and so $v_{n-1} \geq u^*$ on $\B_2$). Fix a point $z_0
\in \B_2 \backslash \R^2$ and a $J$-holomorphic disc $f$ satisfying $f(0)
= z_0$.

\vskip 0,1cm
{\bf Claim.} {\it The interior of the set $f^{-1}(\R^2) \subset \Delta$  is
empty and $\mu(f^{-1}(\R^2)) = 0$.}

\vskip 0,1cm
Indeed, the set $S \subset \Delta$ of critical points of $f$ is
discrete (see~\cite{mc-sa}). If $p \in f^{-1}(\R^2)$ is not a critical point, the
tangent
space to $f(\Delta)$ at $f(p)$ is not contained in $T_{f(p)}\R^2 = \R^2$  since
$\R^2$ is totally real. So in a neighborhood of $p$ the pullback
$f^{-1}(\R^2)$ is contained in a smooth real curve. This implies the
desired assertion. \qed

\begin{lemma}\label{ext}
The function $u \circ f$ defined and subharmonic on $\Delta \backslash
f^{-1}(\R^2)$ extends as a subharmonic function on $\Delta$.
\end{lemma}

\noindent{Proof.} Let $p \in f^{-1}(\R^2)$ be a non-critical point of $f$.
Consider a small open disc $U = p + r_0\Delta$ centered at $p$. It follows by
~\cite{ni-wo} that

\vskip0,1cm
{\bf Claim.} {\it There exists a sequence of $(J_{st},J)$-holomorphic maps
$f_k : U \rightarrow \mathbb B^2$ uniformly converging to $f$ and such that $U \cap
f^{-1}_k(\R^2) = \{ p \}$ for every $k$.}

Since $u \circ f_k$ is a bounded subharmonic function on $U \backslash \{
p \}$, it extends as a subharmonic function on $U$.  By continuity of $u$
(this is the only place of the proof where we use this assumption!), the
sequence $(u \circ f_k)$ converges to $u \circ f$ on $U \backslash \{ p
\}$. Therefore, by the Lebesgue dominated convergence theorem, $(u \circ f_k)$ tends to
$u \circ f$ on $U$ in the distribution sense and the
generalized laplacian $\Delta(u \circ f)$ is a positive
distribution. Therefore, $u \circ f$ admits a subharmonic extension
$\widetilde{u \circ f}$ on $U$ (given by $\widetilde{u \circ f}(q) = \lim
\sup_{\zeta \longrightarrow q} u(\zeta)$ for $q \in
f^{-1}(\R^2)$). Therefore, since the set of critical points of $f$ is
discrete, $u \circ f$ extends subharmonically on $\Delta$, proving
Lemma~\ref{ext}. \qed

\vskip 0,1cm
Using Lemma~\ref{ext} and the induction step we have~:

$$
\int_{\Delta} v_{n-1} \circ f d\mu \geq \int_{\Delta} u^* \circ f d\mu =
\int_{\Delta \backslash f^{-1}(\R^2)} u \circ f d\mu \geq u(z_0).
$$
Hence, $v_n(z_0) \geq u(z_0)$ for every $z_0 \in \B_2 \backslash
\R^2$. This proves Proposition~\ref{prop2.2}. \qed

\vskip 0,1cm
\noindent{\it Proof of Theorem~\ref{rem-theo}}.
Since $\hat v$ is upper semicontinuous on $\B_2$ and $\hat{v} = u$ on $\B_2 \backslash
\R^2$ by Proposition~\ref{prop2.2}, we have $\hat v  \geq u^*$ on $\B_2$.
However, $\hat{v} \leq u^*$ by construction (see
Proposition~\ref{pro2.1}). Hence $\hat u = u^*$ on $\B_2$. This proves
Theorem~\ref{rem-theo}. \qed

%\vskip 0,1cm
%\noindent{\bf Related results.}

\vskip 0,2cm
Our method can be easily carried out to almost complex manifolds of any
dimension.

\begin{theorem}
Let $(M,J)$ be an almost complex manifold  and
let $E \subset M$ be a genericl submanifold of $M$ of real
codimension $2$. Then for any continuous plurisubharmonic function $u$ on
$M \backslash E$ the function $u^*$ defined by $u^*(x) = u(x)$ for $x \in
M \backslash E$ and $u^*(x) = \lim \sup_{y \in M \backslash E, y
\longrightarrow x} u(x)$ for $x \in E$, is $J$-plurisubharmonic on $M$.
\end{theorem}

A slight modification of these technics can be used in order to obtain the
following almost complex analogues of well-known results. Denote by
${\mathcal H}_d$ the Hausdorff measure of dimension $d$.

\begin{theorem}
Let $E$ be a subset of an almost complex manifold $(M,J)$ of real
dimension $2n$.
\begin{itemize}
\item[(a)] If ${\mathcal H}_{2n-2}(E) \leq \infty$, then $E$ is removal
for every bounded continuous $J$-plurisubharmonic function.
\item[(b)] If ${\mathcal H}_{2n-3}(E) = 0$, then $E$ is removal for every
continuous $J$-plurisubharmonic function.
\end{itemize}
\end{theorem}

\bigskip
\subsection{Canonical almost complex bundles}

We present two constructions of almost complex structures on the tangent
and on the cotangent bundles of an almost complex manifold.

\subsubsection{Lift to the tangent bundle}
We endow the tangent bundle $TM$ with the complete lift $J^c$ of $J$
(see~\cite{ya-is} for its definition). We recall that $J^c$ is an
almost complex structure on $TM$. Moreover, if $\nabla$ is any $J$
complex connection on $M$ (ie $\nabla J=0$) and $\bar{\nabla}$ is the
connection defined on $M$ by $\bar\nabla_XY = \nabla_YX +[X,Y]$ then
$J^c$ is the horizontal lift of $J$ with respect to
$\bar\nabla$. Another definition of $J^c$ is given in \cite{le-sz}
where this is characterized by a deformation property.
The equality between the two definitions given in \cite{ya-is} and in
\cite{le-sz} is obtained by their (equal) expression in the local
canonical coordinates on $TM$~:
$$
J^c=\left(
\begin{array}{ccc}
J_i^h & & (0)\\
  & &    \\
t^a\partial_a J_i^h & & J_i^h
\end{array}
\right).
$$
(Here $t^a$ are fibers coordinates).

\subsubsection{Canonical lift to the cotangent bundle}
We recall the definition of the canonical lift of an almost
complex structure $J$ on $M$ to the cotangent bundle $T^*M$, following~
\cite{ya-is}. Set $m=
2n$. We use the following notations. Suffixes A,B,C,D take the values
$1$ to $2m$, suffixes $a,b,c,\dots$,$h,i,j,\dots$ take the values $1$ to
$m$ and $\bar j = j+ m$, $\dots$ The summation notation for
repeated indices is used. If the notation $(\varepsilon_{AB})$,
$(\varepsilon^{AB})$, $(F_{B}^{{}A})$ is used for matrices, the suffix
on the left indicates the column and the suffix on the right indicates
the row. We denote local coordinates on $M$ by $(x^1,\dots,x^n)$ and by
$(p_1,\dots,p_n)$ the fiber coordinates.

Recall that the cotangent space $T^*(M)$ of $M$ possesses the {\it
canonical contact form} $ \theta$ given in local coordinates by
$\theta = p_idx^i.$
The cotangent lift $\varphi^*$ of any diffeomorphism $\varphi$ of $M$
is contact with respect to $\theta$, that is $\theta$ does not depend on
the choice of local coordinates on $T^*(M)$. 

The exterior derivative $d\theta$ of $\theta$ defines {\it the
canonical
symplectic structure} of $T^*(M)$:
$d\theta = dp_i \wedge dx^i$
which is also independent of local coordinates in view of the
invariance of the exterior derivative. Setting $d\theta =
(1/2)\varepsilon_{CB}dx^C \wedge dx^B$ (where $dx^{\bar j} =
dp_j$), we have

$$
(\varepsilon_{CB}) = \left(
\begin{array}{cc}
0 & I_n \\
-I_n  & 0
\end{array}
\right).
$$

Denote by $(\varepsilon^{BA})$ the inverse matrix and write
$\varepsilon^{-1}$ for the tensor field of type (2,0) whose component
are $(\varepsilon^{BA})$. By construction, this definition does not
depend on the choice of local coordinates.

Let now $E$ be a tensor field of type (1,1) on $M$. If $E$ has
components $E_i^{\, h}$ and $E_i^{*h}$ relative to local coordinates $x$
and $x^*$ respectively, then
$p_a^*E_{i}^{*\,a} = p_aE_j^{\, b}\frac{\partial x^j}{\partial x^{*i}}.$
If we interpret a change of coordinates as a diffeomorphism $x^* =
x^*(x) = \varphi(x)$ we denote by $E^*$ the direct image of the tensor $E$
under the action of $\varphi$. In the case where $E$ is an almost
complex structure (that is $E^2 = -Id$), then $\varphi$ is a biholomorphism
between $(M,E)$ and $(M,E^*)$. Any (1,1) tensor field $E$ on $M$
canonically defines a contact form on $E^*M$ via $
\sigma = p_aE_b^{\, a}dx^b.$
Since $(\varphi^*)^*( p_a^*E_b^{{*}\,a}dx^{*b}) = \sigma,$
$\sigma$ does not depend on a choice of local coordinates (here
$\varphi^*$ is the cotangent lift of $\varphi$). Then this canonically
defines the symplectic form

$$
d\sigma = p_a\frac{\partial E_b^{\, a}}{\partial x^c}dx^c \wedge dx^b
+ E_b^{\, a}dp_a \wedge dx^b.
$$
The cotangent lift $\varphi^*$ of a diffeomorphism $\varphi$ is a
symplectomorphism for $d\sigma$. We may write
$d\sigma = (1/2)\tau_{CB}dx^C \wedge dx^B$
where $x^{\bar i} = p_i$; so we have

$$
\tau_{ji} = p_a \left ( \frac{\partial E_i^{\, a}}{\partial x^j} -
\frac{\partial E_j^{\, a}}{\partial x^i} \right ), \tau_{\bar{j} i} =
E_i^{\, j}, \tau_{j\bar{i}} = -E_j^{\, i}, \tau_{\bar{j}\bar{i}} = 0.
$$

We write $\widehat{E}$ for the tensor field of type (1,1) on $T^*(M)$ whose
components $\widehat{E}_B^{{}A}$  are given by 
$\widehat{E}_B^{{}A} = \tau_{BC}\varepsilon^{CA}.$
Thus $
\widehat{E}_i^{\, h} = E_i^{\, h}, \ \widehat{E}_{\bar{i}}^{\, h} = 0$
and
$\widehat{E}_i^{\, \bar{h}} = 
p_a \left( \frac{\partial E_i^{\, a}}{\partial x^j} -
\frac{\partial E_j^{\, a}}{\partial x^i} \right ), 
\widehat{E}_{\bar{i}}^{\,\bar{h}} = E_h^{\, i}.$
In the matrix form we have

$$
\widehat{E} = \left(
\begin{array}{cll}
E_i^{\, h} & & 0 \\
 p_a \left ( \frac{\partial E_i^{\, a}}{\partial x^j} -
\frac{\partial E_j^{\, a}}{\partial x^i} \right ) & & E_h^{\, i}
\end{array}
\right).
$$

By construction, the complete lift $\widehat{E}$ has the following
{\it invariance property}~: if $\varphi$ is a local diffeomorphism of $M$
transforming $E$ to $E'$, then the direct image of $\widehat{E}$ under the
cotangent lift $\psi: = \varphi^*$ is $\widehat{E'}$.
In general, $\widehat{E}$ is not an almost complex
structure, even if $E$ is.  Moreover, one can show \cite{ya-is} that
$\widehat{J}$ is a complex structure if and only if $J$ is integrable.
One may however construct an almost complex structure on $T^*(M)$ as follows.

Let $S$ be a tensor field of type (1,s) on $M$. We may consider the
tensor field $\gamma S$ of type $(1,s-1)$ on $T^*M$, defined in local
canonical coordinates on $T^*M$ by the expression

$$
\gamma S = p_aS_{i_s...i_2i_1}^{\,\,\,a}dx^{i_s} \otimes
\cdots \otimes dx^{i_2}\otimes \frac{\partial}{\partial p_{i_1}}.
$$

In particular, if $T$ is a tensor field of type (1,2) on $M$, then $\gamma T$
has components

$$
\gamma T = \left(
\begin{array}{cll}
0 & & 0\\
 p_a T_{ji}^{\,\,a} & & 0 
\end{array}
\right)
$$
in the local canonical coordinates on $T^*M$.

Let $F$ be a (1,1) tensor field on $M$. Its Nijenhuis tensor $N$
is the tensor field of type (1,2) on $M$ acting on two vector fields $X$ and
$Y$ by

$$
N(X,Y) = [FX,FY] - F[FX,Y] - F[X,FY] + F^2[X,Y].
$$

By $NF$ we denote the tensor field acting by $(NF)(X,Y) =
N(X,FY)$. The following proposition is proved in \cite{ya-is} (p.256).

\begin{proposition}\label{prop-lift}
Let $J$ be an almost complex structure on $M$. Then
\begin{equation}\label{EEE}
\tilde J := \widehat J + (1/2)\gamma(NJ)
\end{equation}
is an almost complex structure on the cotangent bundle $T^*(M)$.
\end{proposition}

We stress out that the definition of the tensor $\tilde{J}$ is independent of
the choice of coordinates on $T^*M$. Therefore if $\phi$ is a biholomorphism
between two almost complex manifolds $(M,J)$ and $(M',J')$, then its
cotangent lift $\tilde \phi:= (\phi,{}^td\phi^{-1})$is a biholomorphism
between $(T^*(M),\tilde J)$ and
$(T^*(M'), \tilde J')$. Indeed one can view $\phi$ as a change of coordinates
on $M$, $J'$ representing $J$ in the new coordinates. The cotangent lift
$\phi^*$ defines a change of coordinates on $T^*M$ and $\tilde{J}'$
represents $\tilde{J}$ in the new coordinates.
So the assertion $(i)$ of Proposition~\ref{prop-lift} holds. Property
$(ii)$ of Proposition~\ref{prop-lift} is immediate in view of the definition
of $\tilde{J}$ given by~(\ref{EEE}).

\subsubsection{Conormal bundle of a submanifold}
The conormal bundle of a strictly pseudoconvex hypersurface in $(M,J)$
provides an important example of a totally real submanifold in the
cotangent bundle $T^*M$, endowed with the canonical almost complex
structure $\tilde{J}$ defined in the last Subsection.
%Moreover, if $f$ is a biholomorphism
%between $(M,J)$ and $(M',J')$ then the cotangent map
%$\tilde f:= (f,{}^tdf^{-1})$ is a biholomorphism between $(T^*M,\tilde J)$
%and $(T^*M',\tilde J')$.
If $\Gamma$ is a real submanifold in $M$, the conormal bundle
$\Sigma_J(\Gamma)$ of $\Gamma$ is the real subbundle of $T^*_{(1,0)}M$
defined by
$$
\Sigma_J(\Gamma) = \{ \phi \in T^*_{(1,0)}M: Re \,\phi
\vert_{T\Gamma} = 0\}.
$$
One can identify the
conormal bundle $\Sigma_J(\Gamma)$ of $\Gamma$ with any of the
following subbundles of $T^*M$~: $N_1(\Gamma)=\{\varphi \in T^*M :
\varphi_{|T\Gamma}=0\}$ and $N_2(\Gamma)=\{\varphi \in T^*M :
\varphi_{|JT\Gamma}=0\}$.

\begin{proposition}\label{prop-tot-real}
Let $\Gamma$ be a $\CC^2$ real hypersurface in $(M,J)$.  If $\Gamma$
is strictly $J$-pseudoconvex, then the bundles $N_1(\Gamma)$ and
$N_2(\Gamma)$ (except the zero section) are totally real submanifolds
of dimension $2n$ in $T^*M$ equipped with $\tilde{J}$.
\end{proposition} 
%There is in fact an equivalence between the nondegeneracy of the Levi form
%of $\Gamma$ and the fact that the conormal bundle of $\Gamma$ is totally
%real. However we focus on the implication suitable for our purpose.
Proposition~\ref{prop-tot-real} is due to A.Tumanov~\cite{tu01} in the
integrable case. The question whether a similar result was true in the
almost complex case was asked by the second author to A.Spiro who
gave a positive answer~\cite{sp}. For completeness we give an
alternative proof of this fact. 

\vskip 0,1cm
\noindent{\it Proof of Proposition~\ref{prop-tot-real}}. Let $x_0 \in \Gamma$.
We consider local coordinates $(x,p)$ for the real cotangent bundle
$T^*M$ of $M$ in a neighborhood of $x_0$.
The fiber of $N_2(\Gamma)$ is given by $c(x)J^*d\rho(x)$, where $c$ is a real
nonvanishing function. In what follows we denote $-J^*d\rho$
by $d^c_J\rho$. For every $\varphi \in N_2(\Gamma)$
we have $\varphi|_{J(T\Gamma)} \equiv 0$. It is equivalenty to prove that
$N_1(\Gamma)$ is totally real in $(T^*M,\tilde J)$ or that $N_2(\Gamma)$ is
totally real in $(T^*M,\tilde J)$. We recall
that if $\Theta=p_idx^i$ in local coordinates then $d\Theta$ defines the
canonical symplectic form on $T^*M$.
If $V,W \in T(N_2(\Gamma))$ then
$d\Theta(V,W)=0$. Indeed the projection $pr_1(V)$ of $V$ (resp. $W$) on $M$ is
in $J(T\Gamma)$ and the projection of $V$ (resp. $W$) on the fiber
annihilates $J(T\Gamma)$ by definition. It follows that $N_2(\Gamma)$ is a
Lagrangian submanifold of $T^*M$ for this symplectic form. 

Let $V$ be a vector field in $T(N_2(\Gamma)) \cap \tilde JT(N_2(\Gamma))$.
We wish to prove that $V=0$.
According to what preceeds we have $d\Theta(V,W) = d\Theta(JV,W)=0$ for every
$W \in T(N_2(\Gamma))$. We restrict to $W$ such that $pr_1(W) \in T\Gamma
\cap J(T\Gamma)$. Since $\Theta$ is defined over $x_0 \in \Gamma$ by
$\Theta = cd^c_J \rho$, then
$d\Theta = dc \wedge d^c_J \rho + cdd^c_J\rho$.
Since $d^c_J \rho(pr_1(V)) = d^c_J \rho(Jpr_1(V)) = d^c_J \rho(pr_1(V)) =
d^c_J \rho(Jpr_1(V)) =0$ it follows that
$dd^c_J \rho (pr_1(V),pr_1(\tilde{J}W)) =0$. However, by the definition of
$\tilde{J}$, we know that $pr_1(\tilde{J}W) = J pr_1(W)$.
Hence, choosing $W = V$, we obtain that $dd^c_J \rho (pr_1(V),J pr_1(V)) = 0$.
Since $\Gamma$ is strictly $J$-pseudoconvex, it follows that $pr_1(V) = 0$.
In particular, $V$ is given in local coordinates by $V=(0,pr_2(V))$.
It follows now from the form of $\tilde J$ that $JV=(0,J pr_2(V))$
(we consider $pr_2(V)$ as a vector in $\mathbb R^{2n}$ and $J$ defined
on $\R^{2n}$). Since $N_2(\Gamma)$ is a real bundle of rank one, then
$pr_2(V)$ is equal to zero. \qed

\section{Riemann mapping theorem for almost complex manifolds with boundary}

The aim of this Section is to prove the existence of stationary discs
in the ball for small 
almost complex deformations of the standard structure.
Stationary discs are natural global biholomorphic invariants of complex
manifolds with boundary. L.Lempert proved in \cite{le81} that for a
strictly convex domain stationary discs coincide with extremal discs
for the Kobayashi metric and studied their basic properties. 
Using these discs he introduced a multi dimensional analogue of the
Riemann map. We define here a local
analogue of the 
Riemann map and establish its main properties. These results were
obtained in \cite{CoGaSu}.

\subsection{Existence of discs attached to a real submanifold 
of an almost complex manifold}

\subsubsection{ Partial indices and the Riemann-Hilbert problem}

In this section we introduce basic tools of the linear Riemann-Hilbert 
problem.

Let $V \subset \C^N$ be an open set. We denote
by $\CC^k(V)$ the Banach space of (real or complex valued) functions
of class $\CC^k$ on $V$ with the standard norm
$$
\parallel r \parallel_k =
\sum_{\vert \nu \vert \leq k}
\sup \{ \vert D^{\nu} r(w) \vert : w \in V \}.
$$
For a positive real 
number $\alpha <1$ and a Banach space $X$, we denote by
$\CC^{\alpha}(\partial \Delta,X)$ the Banach space of all functions
$f: \partial \Delta \rightarrow X$ such that

$$
\parallel f \parallel_{\alpha} :=
\sup_{\zeta \in \partial \Delta} \parallel f(\zeta) \parallel +
\sup_{\theta,\eta \in \partial \Delta, \theta \neq \eta}
\frac{\parallel f(\theta) - f(\eta)\parallel}{\vert \theta - \eta
\vert^{\alpha}} < \infty.
$$

\vskip 0,1cm \noindent If $\alpha = m +
\beta$ with an integer $m \geq 0$ and $\beta \in ]0,1[$, then we consider
the Banach space
$$
\CC^{\alpha}(V):=\{r \in \CC^m(V,\R): D^{\nu} r \in C^{\beta}(V), 
\nu, \vert \nu \vert \leq m\}
$$
and we set $\parallel r
\parallel_{\alpha} = 
\sum_{\vert \nu \vert \leq m} \parallel D^{\nu}r \parallel_{\beta}$.

Then a map $f$ is in $\CC^{\alpha}(V,\C^k)$ if and only if its components 
belong to $\CC^{\alpha}(V)$ and we say that $f$ is of class $\CC^{\alpha}$. 

\vskip 0,1cm
Consider the following situation:

\vskip 0,1cm
\noindent $\bullet$ $B$ is an open ball centered at the origin in $\C^N$ and 
$r^1,\dots,r^N$ are smooth $\CC^\infty$ functions defined in a neighborhood 
of  $\partial \Delta \times B$ in $\C^N \times \C$

\noindent $\bullet$ $f$ is a map of class $\CC^{\alpha}$ 
from $\partial \Delta$ to $B$, where $\alpha>1$ is a noninteger real
number

\noindent $\bullet$ for every $\zeta \in \partial \Delta$

\begin{itemize}
\item[(i)] $E(\zeta) = \{ z \in B: r^j(z,\zeta) = 0, 1 \leq j \leq N \}$ 
is a maximal totally real submanifold in $\C^N$,
\item[(ii)] $f(\zeta) \in E(\zeta)$,
\item[(iii)] $\partial_z r^1(z,\zeta) \wedge \cdots \wedge \partial_z
r^N(z,\zeta) \neq 0$ on $B \times \partial \Delta$.
\end{itemize}

Such a family ${E} = \{ E(\zeta) \}$
of manifolds with a fixed disc $f$ is called a
{\it totally real fibration} over the unit circle.
A disc attached to a fixed totally real manifold ($E$ is independent of 
$\zeta$) is a special case of a totally real fibration. 

\vskip 0,1cm
Assume that the defining function $r:=(r^1,\dots,r^N)$ of $E$ 
depends smoothly on a small real parameter $\varepsilon$, namely $r =
r(z,\zeta,\varepsilon)$, 
and that the fibration $E_0:=E(\zeta,0)$ corresponding
to $\varepsilon = 0$ coincides with the above fibration
$E$. Then for every sufficiently small $\varepsilon$ and for every
$\zeta \in \partial \Delta$ the manifold $E_\varepsilon:=E(\zeta,\varepsilon) 
:= \{ z \in B:
r(z,\zeta,\varepsilon) = 0\}$ is totally real. 
We call $E_{\varepsilon}$ a {\it smooth totally real deformation }  
of the totally
real fibration $E$. By a holomorphic disc $\tilde f$ attached to
$E_{\varepsilon}$  we mean a holomorphic map $\tilde f:\Delta
\rightarrow B$, continuous on $\bar\Delta$, satisfying
$r(f(\zeta),\zeta,\varepsilon) = 0$ on $\partial \Delta$. 

For every positive real noninteger $\alpha$ we denote by
$(\mathcal A^\alpha)^N$ the space of maps defined on
$\bar{\Delta}$, $J_{st}$-holomorphic on $\Delta$, and belonging to
$(\mathcal C^\alpha(\bar{\Delta}))^N$.

\subsubsection{Almost complex perturbation of discs}
Consider a smooth deformation $(J_\lambda)$ of $J_{st}$. We recall that for
$\lambda$ small enough the
$J_{\lambda}$-holomorphicity condition for a map $f:\Delta
\rightarrow \C^N$ may be written in the form

\begin{equation}
\bar\partial_{J_{\lambda}} f = \bar\partial f +
q(\lambda,f)\overline{\partial f} = 0
\end{equation}
where $q$ is a smooth $(n \times n)$ complex matrix satisfying $q(0,\cdot) \equiv 0$.

Let $E_\varepsilon=\{r_j(z,\zeta,\varepsilon) = 0, 1 \leq j \leq
N\}$ be a smooth totally real deformation of a totally real fibration $E$. 
A disc $f \in \mathcal (\CC^\alpha(\bar{\Delta}))^N$ is attached
to $E_\varepsilon$ and is $J_\lambda$-holomorphic if and only if it
satisfies the following nonlinear boundary Riemann-Hilbert type  problem~:
$$
\left\{
\begin{array}{lll}
r(f(\zeta),\zeta,\varepsilon) &=& 0, \ \ \ \zeta \in \partial \Delta\\
 & & \\
\bar{\partial}_{J_\lambda}f(\zeta) &=& 0, \ \ \ \zeta \in \Delta.
\end{array}
\right.
$$
Let $f^0 \in \mathcal (\mathcal A^\alpha)^N$ be a disc attached to
$E$ and let $\mathcal U$ be a neighborhood of $(f^0,0,0)$ in the 
space $(\mathcal C^\alpha(\bar{\Delta}))^N \times \R \times \R$.
Given $(f,\varepsilon,\lambda)$ in $U$ define the maps 
$v_{f,\varepsilon,\lambda}: \zeta \in \partial 
\Delta \mapsto r(f(\zeta), \zeta, \varepsilon)$
and
$$
\begin{array}{llcll}
u &:& \mathcal U & \rightarrow & (\mathcal C^\alpha(\partial \Delta))^N \times
 \mathcal C^{\alpha-1}(\Delta)\\
  & & (f,\varepsilon,\lambda) & \mapsto & (v_{f,\varepsilon,\lambda}, 
\bar{\partial}_{J_\lambda}f).
\end{array}
$$

Denote by $X$ the Banach space $(\mathcal C^\alpha(\bar \Delta))^N$.
Since $r$ is of class $\CC^{\infty}$, 
the map
$u$ is smooth and the tangent map $D_Xu(f^0,0,0)$ (we consider 
the derivative
with respect to the space $X$) is a linear map from $X$ to 
$(\mathcal C^\alpha(\partial \Delta))^N \times \mathcal C^{\alpha-1}(\Delta)$,
defined for every $h \in X$ by 
$$\begin{array}{llllll}
D_Xu(f^0,0,0)(h) = \left(
\begin{matrix}
2 Re [G h] \\
\bar\partial_{J_0} h
\end{matrix}
\right),
\end{array}$$
where for $\zeta \in \partial \Delta$
$$\begin{array}{lllllll}
G(\zeta) = \left(
\begin{matrix}
\frac{\partial r_1}{\partial z^1}(f^0(\zeta),\zeta,0)  &\cdots&\frac{\partial
r_1}{\partial z^N}(f^0(\zeta),\zeta,0)\\
\cdots&\cdots&\cdots\\
\frac{\partial r_N}{\partial z^1}(f^0(\zeta),\zeta,0)& \cdots&\frac{\partial r_N}
{\partial z^N}(f^0(\zeta),\zeta,0)
\end{matrix}
\right)
\end{array}$$
(see \cite{gl94}).
\begin{proposition}\label{Tthh}
Let $f^0:\bar \Delta \rightarrow \C^N$ be a $J_{st}$-holomorphic
disc attached to a totally real fibration $E$ in $\C^N$.  Let
$E_\varepsilon$ be a smooth totally real deformation 
of $E$ and
$J_\lambda$ be a smooth almost complex deformation of $J_0$ in a
neighborhood of $f(\bar{\Delta})$.
Assume that for some $\alpha > 1$
the linear map from $(\mathcal A^{\alpha})^N$ to $(\mathcal 
C^{\alpha-1}(\Delta))^N$ given by 
$h \mapsto 2 Re [G h]$
is surjective and has a $k$ dimensional kernel.
Then there exist $\delta_0,
\varepsilon_0, \lambda_0 >0$ such that for every $0 \leq \varepsilon
\leq \varepsilon_0$ and for every $0 \leq \lambda \leq \lambda_0$,
the set of $J_\lambda$-holomorphic discs $f$ attached to $E_\varepsilon$
and such that $\parallel f -f^0 \parallel_{\alpha} \leq \delta_0$ forms
a smooth $k$-dimensional 
submanifold 
$\mathcal A_{\varepsilon,\lambda}$ in the Banach space 
$(\mathcal C^\alpha(\bar{\Delta}))^N$.
\end{proposition} 

\noindent{\bf Proof.} According to the implicit function Theorem, 
the proof of Proposition~\ref{Tthh} reduces to the proof of the 
surjectivity of $D_Xu$. 
It follows by classical
one-variable results on the resolution of the
$\bar\partial$-problem in the unit disc that the linear map from
$X$ to $\mathcal C^{\alpha-1}(\Delta)$ given by 
$h \mapsto \bar \partial h$
is surjective. More precisely, given $g \in \mathcal C^{\alpha-1}(\Delta)$ 
consider
the Cauchy transform 

$$T_{CG}(g) : \tau \in \Delta \mapsto   
\frac{1}{2\pi i}\int\int_{\Delta} \frac{g(\zeta)}{\zeta - \tau}d\zeta d\bar{\zeta}.$$

For every function $g \in \mathcal
C^{\alpha-1}(\Delta)$ the solutions $h \in X$ of the equation
$\bar\partial h = g$ have the form $h = h_0 + T_{CG}(g)$
where $h_0$  is an arbitrary function in $({\mathcal A}^{\alpha})^N$. 
Consider the equation

\begin{equation}\label{equa1}
D_Xu(f^0,0,0)(h) = \left(
\begin{matrix}
g_1 \\
g_2
\end{matrix}
\right),
\end{equation}
where $(g_1,g_2)$ is a vector-valued function with components 
$g_1 \in \mathcal C^{\alpha-1}(\partial \Delta)$ and $g_2 \in \mathcal 
C^{\alpha-1}(\Delta)$. Solving the
$\bar\partial$-equation for the second component, we reduce 
(\ref{equa1}) to 
$$
2 Re [G(\zeta) h_0(\zeta)] = g_1 - 2 Re [G(\zeta) T_{CG}(g_2)(\zeta)]
$$
with respect to $h_0 \in (\mathcal A^{\alpha})^N$. 
The surjectivity of the map $ h_0 \mapsto 2 Re [G h_0]$ gives the result. \qed

\subsubsection{Riemann-Hilbert problem on the cotangent bundle of an almost 
complex manifold} 
Let $(J_{\lambda})_\lambda$ be an almost complex deformation 
of the standard structure $J_{st}$ on $\B_n$, satisfying $J_\lambda(0)=
J_{st}$.Consider the canonical lift $\tilde{J}_\lambda$ on the cotangent
bundle, defined by ~\ref{EEE}. In the $(x,y)$ coordinates we identify this with
the $(4n \times
4n)$-matrix

$$
\tilde{J}_{\lambda} = \left(
\begin{matrix}
J_{\lambda}(x) & 0\\
\sum y_k A^k_{\lambda}(x)&  ^tJ_{\lambda}(x)
\end{matrix}
\right),
$$
where $A_{\lambda}^k(x) = A^k(\lambda,x)$, $A_{\lambda}^k(x)$ are smooth $(2n \times
2n)$-matrix functions. 
In what follows we always assume that: 

\begin{equation}
\label{norm1}
A^k_0(x) \equiv 0, \ \ \ {\rm for \ every}\ k.
\end{equation}

The trivial bundle $\mathbb B \times \R^{2n}$
over the unit ball is a local coordinate
representation of the cotangent bundle of an almost complex
manifold. We denote by $x = (x^1,\dots,x^{2n})\in \mathbb B_n$
and $y = (y_1,\dots,y_{2n}) \in \R^{2n}$ the coordinates on the base and
fibers respectively. We identify the base space $(\R^{2n},x)$
with $(\C^n,z)$. 
Since ${}^tJ_{st}$ is orthogonally equivalent to $J_{st}$ we may identify
$(\mathbb R^{2n},{}^tJ_{st})$ with $(\mathbb C^n,J_{st})$. After this identification
the ${{}^tJ_{st}}$-holomorphicity is expressed by the
$\bar{\partial}$-equation in the usual coordinates in $\C^n$.

Consider a smooth map
$\hat{f}=(f,g): \Delta \rightarrow \mathbb B \times \R^{2n}$ which is 
{\it $(J_{st},\tilde J_{\lambda})$-holomorphic} : 

$$
\tilde J_{\lambda}(f,g) \circ  d\hat{f} = d\hat{f} \circ J_{st}
$$
on $\Delta$.

For $\lambda$ small enough this can be rewritten as the
following Beltrami type quasilinear elliptic equation:

$$
(\mathcal E)\left\{
\begin{array}{cll}
\bar \partial f + q_1(\lambda,f))\overline{\partial f} & = & 0\\
 & & \\
\bar \partial g + q_2(\lambda,f))\overline{\partial g} + q_3(\lambda,f)g +
q_4(\lambda,f) \overline{g}
 & = & 0,
\end{array}
\right.
$$
where the first equation coincides with the $(J_{st},J_{\lambda})$-holomorphicity
condition for $f$ that is $\bar\partial_{J_\lambda} f 
= \bar \partial f + q_1(\lambda,f))\overline{\partial f}$.
The coefficient $q_1$ is uniquely determined by $J_{\lambda}$ and, in view 
of (\ref{norm1}),
the coefficient $q_k$ satisfies, for $k=2,3,4$:

\begin{equation}
\label{nprm3}
q_k(0,\cdot) \equiv 0, \ q_k(\cdot,0) \equiv 0.
\end{equation}

We point out that in $(\mathcal E)$ 
the equations for the fiber component $g$ are
obtained as a small perturbation of the standard $\bar\partial$-operator. 
An important feature of this system is that the second equation is
{\it linear} with respect to the fiber component $g$. 

We consider the operator

$$\bar\partial_{\tilde J_{\lambda}} :     \left(
\begin{matrix}
f\\
g
\end{matrix}
\right) \mapsto  \left(
\begin{matrix}
\bar \partial f + q_1(\lambda,f))\overline{\partial f}\\
\bar \partial g + q_2(\lambda,f))\overline{\partial g} + q_3(\lambda,f) g+
q_4(\lambda,f)\overline{g}
\end{matrix}
\right).$$

\vskip 0,1cm 
Let $r_j(z,t,\lambda)$ , $j = 1,\dots,4n$ be $\CC^{\infty}$-smooth real 
functions
on $\mathbb B \times \mathbb B \times [0,\lambda_0]$ and let 
$r:=(r_1,\dots,r_{4n})$. 
Consider the following nonlinear boundary
Riemann-Hilbert type problem for the operator
$\bar\partial_{\tilde J_{\lambda}}$:
$$
({\mathcal
BP}_{\lambda})
\left\{
\begin{array}{cll}
r(f(\zeta),\zeta^{-1}g(\zeta),\lambda) &=& 0 \ \ {\rm on} \ \partial \Delta\\
 \\
\bar \partial_{\tilde J_\lambda} (f,g) &=& 0 \ \ {\rm on} \ \Delta
\end{array}
\right.
$$
on the space $\CC^{\alpha}(\bar\Delta,\B_n\times \B_n)$.

The boundary problem $({\mathcal
BP}_{\lambda})$ has the following geometric meaning.
Consider the disc $(\hat f,\hat g) := (f, \zeta^{-1}g)$ on 
$\Delta\backslash \{0\}$ and the set
$E_{\lambda} := \{ (z,t): r(z,t,\lambda) = 0 \}$. The
boundary condition in $({\mathcal BP}_{\lambda})$ means that 
$$
(\hat f, \hat g)(\partial \Delta) \subset E_{\lambda}.
$$
This boundary problem has the following {\it invariance
property}. Let $(f,g)$ be a solution of
$({\mathcal BP}_{\lambda})$ and let $\phi$ be a automorphism of 
$\Delta$. 
Then $(f \circ \phi,c g \circ \phi)$ also satisfies the
$\bar\partial_{\tilde J_{\lambda}}$ equation for every complex constant $c$.
In particular, if $\theta \in[0,2\pi]$ is fixed, then the disc
$(f(e^{i\theta}\zeta),e^{-i\theta}g(e^{i\theta}\zeta))$ satisfies the
$\bar\partial_{\tilde J_{\lambda}}$-equation on $\Delta \backslash \{0\}$ and
the boundary of the disc
$(f(e^{i\theta}\zeta),e^{-i\theta}\zeta^{-1}g(e^{i\theta}\zeta))$ is
attached to $E_{\lambda}$.  This implies the following 

\begin{lemma}
\label{circle}
If $(f,g)$ is a solution of $({\mathcal BP}_{\lambda})$,
then 
$\zeta \mapsto (f(e^{i\theta}\zeta),e^{-i\theta}g(e^{i\theta}\zeta))$ is also a
solution of $({\mathcal BP}_{\lambda})$.
\end{lemma}

Suppose that this problem has a solution $(f^0,g^0)$ for $\lambda = 0$
(in view of the above assumptions this solution is holomorphic on
$\Delta$ with respect to the standard structure on $\C^n \times
\C^n$). Using the implicit function theorem we study, 
for sufficiently small $\lambda$, the solutions of
$({\mathcal BP}_{\lambda})$ close to $(f^0,g^0)$. 
As above consider the map $u$
defined in a neighborhood of $(f^0,g^0,0)$ in 
$(\mathcal C^\alpha(\bar{\Delta}))^{4n} \times \mathbb R$ by:  

$$u: (f,g,\lambda) \mapsto \left(
\begin{matrix}
\zeta \in \partial \Delta \mapsto r(f(\zeta),\zeta^{-1}g(\zeta),\lambda)\\
\bar \partial f + q_1(\lambda,f)\overline{\partial f}\\
\bar \partial g + q_2(\lambda,f) \overline{\partial g} + q_3(\lambda,f)  g +
q_4(\lambda,f) \overline{g}
\end{matrix}
\right).
$$

\vskip 0,1cm
If $X:=(\mathcal C^ \alpha(\bar{\Delta}))^{4n}$ 
then its tangent map at $(f^0,g^0,0)$ has the form 

$$\begin{array}{llllll}
D_Xu(f^0,g^0,0): h=(h_1,h_2) \mapsto  \left(
\begin{matrix}
\zeta \in \partial \Delta \mapsto 
2 Re [G(f^0(\zeta),\zeta^{-1}g^0(\zeta),0)h]  \\
\bar \partial h_1 \\
\bar \partial h_2  
\end{matrix}
\right)
\end{array}$$
where for $\zeta \in \partial \Delta$ one has
$$\begin{array}{lllllll}
G(\zeta) = \left(
\begin{matrix}
\frac{\partial r-1}{\partial w_1}(f^0(\zeta),\zeta^{-1}g^0(\zeta),0) 
&\cdots&\frac{\partial r_1}{\partial w_N}(f^0(\zeta),\zeta^{-1}g^0(\zeta),0)\\
\cdots&\cdots&\cdots\\
\frac{\partial r_N}{\partial w_1}(f^0(\zeta),\zeta^{-1}g^0(\zeta),0)&
\cdots&\frac{\partial r_N}
{\partial w_N}(f^0(\zeta),\zeta^{-1}g^0(\zeta),0)
\end{matrix}
\right)
\end{array}$$
with $N = 4n$ and $w = (z,t)$.

If the tangent map $D_Xu(f^0,g^0,0): (\mathcal A^{\alpha})^N
\longrightarrow (\mathcal 
C^{\alpha-1}(\Delta))^N$ is surjective  and has a
finite-dimensional kernel,  we may apply the implicit function theorem
as in Section~2.2 (see Proposition~\ref{Tthh}) 
and conclude to the existence of a
finite-dimensional variety of nearby discs. 
In particular, consider the fibration $E$ over the disc $(f^0,g^0)$
with fibers
$E(\zeta) = \{ (z,t) : r^j(z,\zeta^{-1}t,0) = 0 \}$. Suppose that this
fibration is totally real. Then we have: 

\begin{proposition}\label{tthhh}
Suppose that the fibration $E$ is totally real. 
If the tangent map $D_Xu(f^0,g^0,0): (\mathcal A^{\alpha})^{4n}
\longrightarrow (\mathcal 
C^{\alpha-1}(\Delta))^{4n}$ is surjective 
 and has a finite-dimensional kernel,
then for every sufficiently small $\lambda$
the solutions of the boundary problem $({\mathcal BP}_{\lambda})$ form 
a smooth finite dimensional submanifold in the space $(\CC^{\alpha}(\Delta))^{4n}$.
\end{proposition}

In the next Section we present a sufficient condition for the 
surjectivity of the map $D_Xu(f^0,g^0,0)$. This is due to 
J.Globevnik~\cite{gl94,gl96} for the integrable case and
relies on the partial indices of the totally real fibration along $(f^0,g^0)$.

\subsection{Generation of stationary discs}

Let $D$ be a smoothly bounded domain in $\C^n$ with the
boundary $\Gamma$. According to \cite{le81} 
a continuous map $f:\bar\Delta
\rightarrow \bar D$, holomorphic on $\Delta$,
is called a {\it stationary} disc for $D$ (or
for $\Gamma$) if there
exists a holomorphic map
$\hat{f}: \Delta \backslash \{ 0 \} \rightarrow
T^*_{(1,0)}(\C^n)$, $\hat{f} \neq 0$, continuous on $\bar
\Delta \backslash \{ 0 \}$ and such that
\begin{itemize}
\item[(i)] $\pi \circ \hat{f} = f$
\item[(ii)] $\zeta \mapsto \zeta \hat{f}(\zeta)$ is in ${\mathcal O}(\Delta)$
\item[(iii)] $\hat{f}(\zeta) \in \Sigma_{f(\zeta)}(\Gamma)$ for every $\zeta$
in $\partial \Delta$.
\end{itemize}

We call $\hat{f}$ a {\it lift} of $f$ to the conormal bundle of $\Gamma$
(this is a meromorphic map from $\Delta$  into
$T^{*}_{(1,0)}(\C^n)$ whose values on the unit circle lie on 
$\Sigma(\Gamma)$).

We point out that originally Lempert gave this definition in a different form,
using the natural coordinates on the cotangent bundle of $\C^n$. The
present more geometric version in terms of the conormal
bundle is due to Tumanov~\cite{tu01}. This form is particularly useful
for our goals
since it can be transferred to the almost complex case.
\vskip 0,1cm
Let $f$ be a stationary disc for $\Gamma$.
It follows from Proposition~\ref{prop-tot-real} that if $\Gamma$ is a Levi
nondegenerate
hypersurface, the conormal bundle $\Sigma(\Gamma)$ is a totally real
fibration along $f^*$. Conditions $(i)$, $(ii)$, $(iii)$ may be viewed as
a nonlinear boundary problem considered in Section~2.
If the associated tangent map is surjective, Proposition~\ref{tthhh} 
gives a description of all stationary
discs $\tilde{f}$ close to $f$, for a small deformation of $\Gamma$. 
When dealing with the standard complex structure on $\C^n$, the
bundle $T^*_{(1,0)}(\C^n)$  is a holomorphic vector bundle which can be 
identified, after projectivization of the fibers, with the
holomorphic bundle of complex hyperplanes that is with $\C^n \times
\mathbb P^{n-1}$. The conormal bundle $\Sigma(\Gamma)$ of a real
hypersurface $\Gamma$ may be naturally
identified, after this projectivization, with the bundle 
of holomorphic tangent spaces $H(\Gamma)$
over $\Gamma$. According to S.Webster
\cite{we78} this is a totally real submanifold in $\C^n \times \mathbb
P^{n-1}$. When dealing with the standard structure, we may therefore
work with projectivizations of lifts of stationary discs attached to
the holomorphic tangent bundle $H(\Gamma)$. The technical avantage
is that after such a projectivization lifts of stationary discs become
holomorphic, since the lifts have at most one pole of order 1 at the origin.
This idea was first used
by L.Lempert and then applied by several authors \cite{ba-le98,
ce95, sp-tr02}. 

When we consider almost complex 
deformations of the standard structure (and not just deformations of
$\Gamma$)
the situation is more complicated. If the cotangent bundle $T^*(\R^{2n})$ is
equipped with $\tilde J$, there is no natural
possibility to transfer this structure to the space obtained by the
projectivization of the fibers. Consequently we do not work with
projectivization of the cotangent bundle but we will deal with meromorphic
lifts of stationary discs. Representing such lifts $(\hat{f},\hat{g})$
in the form 
$(\hat{f},\hat{g})=(f,\zeta^{-1}g)$,
we will consider $\tilde J_\lambda$-holomorphic discs close to the 
$J_{st}$-holomorphic disc $(f^0,g^0)$. The disc $(f,g)$
satisfies a nonlinear boundary problem of 
Riemann-Hilbet type $(\mathcal BP_\lambda)$.
When an almost complex structure on the
cotangent bundle is fixed, we may view it as an elliptic prolongation
of an initial almost complex structure on the base and apply the
implicit function theorem as in previous section. This avoids
difficulties coming from the projectivization of almost complex fibre spaces. 

\subsubsection{Maslov index and Globevnik's condition}
We denote by $GL(N,\C)$
the group of invertible $(N \times N)$ complex matrices
and by $GL(N,\R)$ the group of all such matrices with real entries.
Let $0 < \alpha < 1$ and let
$B:\partial \Delta \rightarrow GL(N,\C)$ be of class $\mathcal C^\alpha$. 
According to \cite{ve} (see also \cite{cl-go81}) $B$ admits the factorization
$B(\tau) = F^{+}(\tau)\Lambda(\tau)F^{-}(\tau), \tau \in \partial \Delta$,
where:

$\bullet$ $\Lambda$ is a diagonal matrix of the form
$\Lambda(\tau) = diag(\tau^{k_1},\dots,\tau^{k_N})$, 

$\bullet$ $F^{+}: \bar{\Delta} \rightarrow GL(N,\C)$ is
of class $\mathcal C^\alpha$ on $\bar{\Delta}$ and holomorphic in $\Delta$,

$\bullet$ $F^{-}: [\C \cup \{ \infty \}] \backslash \Delta
\rightarrow GL(N,\C)$
is of class $\mathcal C^\alpha$ on $[\C \cup \{ \infty \}] \backslash \Delta$
and holomorphic on
$[\C \cup \{ \infty \}] \backslash \bar \Delta$.

\vskip 0,1cm The integers $k_1 \geq \cdots \geq k_n$ are called the partial 
indices of $B$. 

\vskip 0,1cm
Let $E$ be a totally real fibration over the unit circle. For every 
$\zeta \in \partial \Delta$ consider the ``normal''
vectors $\nu_j(\zeta) =
(r^j_{\bar{z}^1}(f(\zeta),\zeta),\dots,
r^j_{\bar{z}^N}(f(\zeta),\zeta))$, $j=1,\dots,N$.
We denote by $K(\zeta) \in GL(N,\C)$ the matrix with rows
$\nu_1(\zeta),\dots,\nu_N(\zeta)$ and we set
 $B(\zeta) := -\overline{K(\zeta)}^{-1}K(\zeta)$,
$\zeta \in \partial \Delta$. The partial indices of the map
$B: \partial \Delta \rightarrow GL(N,\C)$
are called {\it the partial indices} of the fibration ${E}$
along the disc $f$ and their sum is
called the {\it total index} or the {\it Maslov index of ${ E}$ along $f$}. 
The following result is due to
J. Globevnik~\cite{gl96}:

\vskip 0,1cm
\noindent{\bf Theorem} : {\it Suppose that all the partial indices of
the totally real fibration ${E}$ along $f$ are $\geq -1$ 
and denote by $k$ the Maslov index of $E$ along 
$f$. Then the linear map from $(\mathcal A^{\alpha})^N$ to $(\mathcal 
C^{\alpha-1}(\Delta))^N$ given by 
$h \mapsto 2 Re [G h]$
is surjective and has a $(N+k)$ dimensional kernel.}

\vskip 0,1cm
Proposition~\ref{Tthh} may be restated in terms of partial indices
as follows~:

\begin{proposition}\label{tthh1}
Let $f^0:\bar \Delta \rightarrow \C^N$ be a $J_{st}$-holomorphic
disc attached to a totally real fibration ${E}$ in $\C^N$. Suppose
that all the partial indices of ${E}$ along $f^0$ are $\geq
-1$. Denote by $k$ the Maslov index of ${E}$ along $f^0$. Let also
$E_\varepsilon$ be a smooth totally real deformation of ${E}$ and
$J_\lambda$ be a smooth almost complex deformation of $J_{st}$ in a
neighborhood of $f(\bar{\Delta})$. Then there exists $\delta_0,
\varepsilon_0, \lambda_0 >0$ such that for every $0 \leq \varepsilon
\leq \varepsilon_0$ and for every $0 \leq \lambda \leq \lambda_0$  
the set of $J_\lambda$-holomorphic discs attached to $E_\varepsilon$
and such that $\parallel f -f^0 \parallel_{\alpha} \leq \delta_0$ forms
a smooth $(N+k)$-dimensional submanifold 
$\mathcal A_{\varepsilon,\lambda}$ in the Banach space 
$(\mathcal C^\alpha(\bar{\Delta}))^N$.
\end{proposition}

Globevnik's result was applied to the study of stationary discs in
some classes of domains in $\C^n$ by M.Cerne~\cite{ce95} and
A.Spiro-S.Trapani~\cite{sp-tr02}.  Since they worked with the
projectivization of the conormal bundle, we explicitely compute, for
reader's convenience and completeness of exposition, partial indices
of {\it meromorphic} lifts of stationary discs for the unit sphere in
$\C^n$.

\vskip 0,1cm
Consider the unit sphere $\Gamma:=\{z \in \C^n : 
z^1\bar{z}^1 + \cdots +z^{n}\bar{z}^{n} -1 = 0\}$ in $\C^n$. 
The conormal bundle $\Sigma(\Gamma)$ is given in the $(z,t)$ coordinates 
 by the equations 
$$
(S)\left\{
\begin{array}{cll}
z^1\bar{z}^1 + \cdots +z^{n}\bar{z}^{n} - 1 = 0, & &\\
t_1 = c \bar{z}^1,\dots,t_{n} = c \bar{z}^{n}, &c \in \mathbb R.&
\end{array}
\right.
$$
According to \cite{le81}, every stationary disc for $\Gamma$
 is extremal for the Kobayashi metric. Therefore,
such a stationary disc $f^0$ centered at the origin is linear by the
Schwarz lemma. So, up to a unitary transformation, we have
 $f^0(\zeta) = (\zeta,0,\dots,0)$ with lift 
$(\widehat{f^0},\widehat{g^0})(\zeta) =
(\zeta,0,\dots,0,\zeta^{-1},0,\dots,0)=(f^0,\zeta^{-1}g^0)$ to the 
conormal bundle.
Representing nearby meromorphic discs in the form
$(z,\zeta^{-1}w)$ and eliminating the parameter $c$ in system 
$(S)$ we obtain that holomorphic discs
$(z,w)$ close to $(f^0,g^0)$ satisfy for $\zeta \in \partial \Delta$:

$$
\begin{array}{lll}
r^1(z,w) & = & z^1\bar{z}^1 + \cdots +z^{n}\bar{z}^{n} - 1 = 0,\\
r^2(z,w) & = & i z^1 w_1\zeta^{-1} - i\bar{z}^1\bar{w}_1\zeta = 0,\\
r^3(z,w) & = & \bar{z}^1w_2\zeta^{-1} - \bar{z}^2w_1\zeta^{-1} + 
z^1\bar{w}_2 \zeta-
z^2\bar{w}_1\zeta = 0,\\
r^4(z,w) & = & i\bar{z}^1w_2\zeta^{-1} - i\bar{z}^2w_1\zeta^{-1} -
 iz^1\bar{w}_2\zeta +
iz^2\bar{w}_1\zeta = 0,\\
r^5(z,w) & = & \bar{z}^1w_3\zeta^{-1} - \bar{z}^3w_1\zeta^{-1} +
 z^1\bar{w}_3 \zeta-
z^3\bar{w}_1\zeta = 0,\\
r^6(z,w) & = & i\bar{z}^1w_3\zeta^{-1} - i\bar{z}^3w_1\zeta^{-1} - 
iz^1\bar{w}_3\zeta +
iz^3\bar{w}_1\zeta = 0,\\
& &\cdots\\
r^{2n-1}(z,w) & = & \bar{z}^1w_n\zeta^{-1} - \bar{z}^nw_1\zeta^{-1} + 
z^1\bar{w}_n \zeta-
z^n\bar{w}_1\zeta = 0,\\
r^{2n}(z,w) & = & i\bar{z}^1w_n\zeta^{-1} - i\bar{z}^nw_1\zeta^{-1} - 
iz^1\bar{w}_n\zeta +
iz^n\bar{w}_1\zeta = 0.
\end{array}
$$

Hence the $(2n \times 2n)$-matrix $K(\zeta)$ has the following expression:
$$
\left( 
\begin{matrix}
%{ccccccccccc}
\zeta&0&0& \cdots &0&0&0&0& \cdots & 0\\
-i\zeta&0&0& \cdots &0&-i&0&0&\cdots &0\\
0&-\zeta^{-1}&0& \cdots &0&0&\zeta^2&0&\cdots &0\\
0&-i\zeta^{-1}&0& \cdots &0&0&-i\zeta^2&0&\cdots &0\\
0&0&-\zeta^{-1}& \cdots &0&0&0&\zeta^2&\cdots &0\\
0&0&-i\zeta^{-1}& \cdots &0&0&0&-i\zeta^2&\cdots &0\\
\cdots & \cdots & \cdots & \cdots & 
\cdots & \cdots & \cdots & \cdots & \cdots &\cdots\\
0&0&0& \cdots &-\zeta^{-1}&0&0&0& \cdots &\zeta^2\\
0&0&0& \cdots &-i\zeta^{-1}&0&0&0& \cdots &-i\zeta^2
\end{matrix}
\right)
$$
and a direct computation shows that $-B =  \bar{K}^{-1} K$ has
the form 
$$
\left( \begin{matrix}C_1&C_2\\
C_3&C_4
\end{matrix}
\right),
$$
where the $(n \times n)$ matrices $C_1,\dots,C_4$ are given by

$$
C_1 = \left( \begin{matrix}\zeta^2&0&.&0\\
0&0&.&0\\
0&0&.&0\\
.&.&.&.\\
0&0&.&0
\end{matrix}
\right), \ \ C_2 = \left( \begin{matrix}0&0&.&0\\
0&-\zeta&.&0\\
0&0&-\zeta&0\\
.&.&.&.\\
0&0&0&-\zeta
\end{matrix}
\right),
$$
$$
C_3 = \left( \begin{matrix}-2\zeta&0&.&0\\
0&-\zeta&.&0\\
0&0&.&0\\
.&.&.&.\\
0&0&.&-\zeta
\end{matrix}
\right), \ \ C_4 = \left( \begin{matrix}-1&0&.&0\\
0&0&.&0\\
0&0&.&0\\
.&.&.&.\\
0&0&.&0
\end{matrix}
\right).
$$
\vskip 0,1cm

We point out that the matrix 

$$\left( \begin{matrix}\zeta^2&0\\
\zeta&1
\end{matrix}
\right)$$
admits the following factorization:

$$\left( \begin{matrix}1&\zeta\\
0&1
\end{matrix} \right ) \times \left( \begin{matrix}-\zeta&0\\
0&\zeta
\end{matrix} 
\right) \times 
\left( \begin{matrix}0&1\\
1&\zeta^{-1}
\end{matrix} \right).$$

Permutating the lines (that is multiplying $B$ by some nondegenerate matrics 
with constant coefficients) and using the above factorization of $(2 \times
2)$ matrices, we obtain the following 

\begin{proposition}\label{PPRRR}
All the partial indices of the conormal bundle of the unit sphere along a 
meromorphic lift of a stationary disc are equal to one 
and the Maslov index is equal to $2n$.
\end{proposition}

Proposition~\ref{PPRRR} enables to apply Proposition~\ref{Tthh}
to construct the family of stationary discs attached to the unit
sphere after a small deformation of the complex structure. Indeed
denote by $r_j(z,w,\zeta,\lambda)$ $\CC^{\infty}$-smooth functions 
coinciding for
$\lambda = 0$ with the above functions $r_1,\dots,r_{2n}$.  

\vskip 0,1cm
In the end of this Subsection we
make the two following assumptions:
\begin{itemize}
\item[(i)] $r_1(z,w,\zeta,\lambda) =  z^1\bar{z}^1 +
\cdots+z^{n}\bar{z}^{n} - 1$, meaning that the sphere is not deformed
\item[(ii)] $r_j(z,tw,\zeta,\lambda)
= tr^j(z,w,\zeta,\lambda)$ for every $j \geq 2, \ t \in \R$.
\end{itemize}

Geometrically this means that given $\lambda$, the set $\{ (z,w):
r_j(z,w,\lambda) = 0 \}$ is a real vector bundle with one-dimensional
fibers over the unit sphere.

Consider an almost
complex deformation $J_{\lambda}$ of the standard structure on $\mathbb B_n$
and its canonical lift $\tilde J_\lambda$ to the cotangent bundle $\mathbb B_n \times 
\mathbb R^{2n}$.
Consider now the corresponding boundary problem:
$$
({\mathcal BP}_{\lambda})
\left\{
\begin{array}{lll}
& &r(f,g,\zeta,\lambda) = 0, \zeta \in \partial \Delta,\\
%& & \\
& &\bar \partial f + q_1(\lambda,f)\overline{\partial f} = 0,\\
%& & \\
& &\bar \partial  g + q_2(\lambda,f)\overline{\partial g} +
q_3(\lambda,f) g + q_4(\lambda,f)\overline{g}= 0.
\end{array}
\right.
$$

Combining Proposition~\ref{PPRRR} with the previous results 
we obtain the following
\begin{proposition}\label{THEO}
For every sufficiently small positive $\lambda$, the set of solutions of
 $({\mathcal BP}_{\lambda})$, close enough to the disc
$(\widehat{f^0},\widehat{g^0})$, forms a smooth $4n$-dimensional
submanifold $V_{\lambda}$ in the space
$\CC^{\alpha}(\bar \Delta)$ (for every noninteger $\alpha > 1$).
\end{proposition}
Moreover, in view of the assumption (ii) and of
the linearity of $({\mathcal E})$ with respect to the fiber component $g$,
we also have the following

\begin{corollary}\label{CORO}
The projections of discs from
$V_{\lambda}$
to the base $(\R^{2n},J_{\lambda})$ form a $(4n-1)$-dimensional subvariety in
$\CC^{\alpha}(\bar\Delta)$.
\end{corollary}

Geometrically the solutions $(f,g)$ of the boundary problem
$({\mathcal BP}_{\lambda})$ are
such that the discs $(f,\zeta^{-1}g)$ are attached to the conormal
bundle of the unit sphere with respect to the standard structure. In
particular, if $\lambda = 0$ then every such disc satisfying $f(0) =
0$ is linear.

\subsection{Canonical Foliation and  the ``Riemann map'' associated with an almost
complex structure}
In this Section we study the geometry of stationary discs in the unit ball 
after
a small almost complex perturbation of the standard structure. 
The idea is simple since these discs are small deformations of the 
complex lines passing through the origin in the unit ball.  
\subsubsection{Foliation associated with an elliptic prolongation} 
Fix a vector $v^0$ with $\vert\vert v^0 \vert\vert = 1$ and consider the
corresponding stationary disc $f^0:\zeta \mapsto \zeta v^0$. Denote by
$(\widehat{f^0},\widehat{g^0})$ its lift to the conormal bundle of
the unit sphere. 
Consider a smooth deformation $J_{\lambda}$ of the standard structure
on the unit ball ${\mathbb B_n}$ in $\C^n$.
For sufficiently small $\lambda_0$ consider the lift $\tilde J_{\lambda}$ on
$\mathbb B_n \times \mathbb R^{2n}$, 
where $\lambda \leq
\lambda_0$. Then the solutions of the associated boundary
problem $({\mathcal BP}_{\lambda})$  form a $4n$-parameter family
of $\tilde J_\lambda$-holomorphic maps from $\Delta$ to
$\C^{n} \times \C^{n}$. Given such a solution $(f^\lambda,g^\lambda)$, 
consider the disc
$(\widehat{f^\lambda},\widehat{g^\lambda}) := 
(f^\lambda, \zeta^{-1}g^\lambda)$. 
In the case where $\lambda = 0$ 
this is just the lift of a stationary disc for the unit
sphere to its conormal bundle. The set of solutions of the problem
$(\mathcal BP_\lambda)$, 
close to $(\widehat{f^0},\widehat{g^0})$, forms a
smooth submanifold of real dimension $4n$ 
in $(\mathcal C^\alpha(\bar{\Delta}))^{4n}$ according
to Proposition~\ref{THEO}. Hence there is a neighborhood $V_0$ of $v^0$ in
$\mathbb R^{2n}$ and a smooth real hypersurface $I_{v^0}^\lambda$ in $V_0$ 
such that for every $\lambda \leq \lambda_0$ and for every
$v \in I_{v^0}^\lambda$ there is one and only one solution 
$(f^\lambda_v,g^\lambda_v)$ of $(\mathcal BP_\lambda)$, up to multiplication
of the fiber component $g^\lambda_v$ by a real constant,
such that 
$f^\lambda_v(0) = 0$ and $df^\lambda_v(0)(\partial / \partial Re(\zeta)) = v$.

We may therefore consider the map
$$
F_0^{\lambda}: (v,\zeta) \in I_{v^0}^\lambda \times \bar{\Delta} \mapsto
(f^\lambda_v,g^\lambda_v)(\zeta).
$$

This is a smooth map with respect to $\lambda$ 
close to the origin in $\mathbb R$.

Denote by $\pi$ the canonical
projection $\pi: \mathbb B_n \times \R^{2n} \rightarrow \mathbb B_n$
and consider the composition $\widehat{F_0^{\lambda}} = \pi \circ
F_0^{\lambda}$. 
This is a smooth map defined for $0 \leq \lambda <
\lambda_0$
and such that

\begin{itemize}
\item[(i)] $\widehat{F_0^{0}}(v,\zeta) 
= v\zeta$, for every $\zeta \in \bar\Delta$ and for every $v \in I_{v^0}^\lambda$.
\item[(ii)] For every $\lambda \leq \lambda_0$, 
$\widehat{F_0^\lambda}(v,0) = 0$.
\item[(iii)] For every fixed $\lambda \leq \lambda_0$ and every 
$v \in I_{v^0}^\lambda$ the map 
$\widehat{F_0^{\lambda}}(v,\cdot)$ is
a $J_{\lambda}$-holomorphic disc attached to the unit sphere.
\item[(iv)] For every fixed $\lambda$, different values of
$v \in I_{v^0}^\lambda$ define different discs.
\end{itemize}

\begin{definition} We call the family 
$(\widehat{F_0^\lambda}(v,\cdot))_{v \in I_{v^0}^\lambda}$ {\it canonical }
discs associated with the boundary problem $({\mathcal
BP}_{\lambda})$.
\end{definition}

 We stress out that by 
a canonical disc {\it we always mean a disc centered at the
origin}. The preceding condition $(iv)$ may be restated as follows:

\begin{lemma}\label{LEMMA1} 
For $\lambda < \lambda_0$ every canonical disc is uniquely 
determined by its tangent vector at the origin.
\end{lemma}

In the next Subsection we glue the sets $I_{v}^\lambda$, depending on vectors
$v \in \mathbb S^{2n-1}$, to define the global indicatrix of $F^\lambda$.

\subsubsection{Indicatrix} For $\lambda < \lambda_0$ consider canonical discs
in $(\mathbb B_n,J_{\lambda})$ centered at the origin and admitting lifts
close to $(\widehat{f^0},\widehat{g^0})$. 

As above we denote by $I_{v^0}^{\lambda}$ 
the set of tangent vectors at the origin of
canonical discs whose lift is close to $(\widehat{f^0},\widehat{g^0})$. 
Since these vectors depend smoothly on
parameters $v$ close to $v^0$ and $\lambda \leq \lambda_0$, 
$I_{v^0}^{\lambda}$ is a smooth deformation of a piece of
the unit sphere $\mathbb S^{2n-1}$. So this is a smooth real
hypersurface in $\C^n$ in a neighborhood of $v^0$.  
Repeating this construction for every vector $v \in \mathbb S^{2n-1}$ 
we may find a finite covering of $\mathbb S^{2n-1}$ by open
connected sets $U_j$ such that for every $j$ the nearby stationary
discs with tangent vectors at the origin close to $v$ are given by
$\widehat{F^{\lambda}_j}$. Since every nearby stationary disc is uniquely
determined by its tangent vector at the origin, we may glue the maps
$\widehat{F^{\lambda}_j}$ to the map $\widehat{F^{\lambda}}$ 
defined for every $v
\in \mathbb S^{2n-1}$ and every $\zeta \in \bar{\Delta}$. The tangent vectors
of the constructed family of stationary discs form a smooth real
hypersurface $I^{\lambda}$ which is a small deformation of the unit
sphere. This hypersurface is an analog of the indicatrix for the
Kobayashi metric (more precisely, its boundary).

We point out that the local indicatrix $I_{v^0}^\lambda$ 
for some fixed $v^0 \in \mathbb S^{2n-1}$ is also useful.

\subsubsection{Circled property and Riemann map}
If $\lambda$  is  small enough, the hypersurface 
$I^{\lambda}$ is strictly pseudoconvex with respect to the standard
structure. Another important property of the ``indicatrix'' is its
invariance with respect to the linear action of the unit circle.

Let $\lambda \leq \lambda_0$, $v \in I^\lambda$ and 
$f_v^\lambda:=\widehat{F^\lambda}(v,\cdot)$. For $\theta \in \mathbb R$ 
we denote by
$f_{v,\theta}^\lambda$ the $J_\lambda$-holomorphic disc in $\mathbb B_n$
defined by $f^\lambda_{v,\theta} : \zeta \in \Delta \mapsto
f_v^\lambda(e^{i\theta}\zeta)$. We have~:

\begin{lemma}\label{LEMMA2}
For every $0 \leq \lambda < \lambda_0$, every $v \in I^\lambda$
and every $\theta \in \R$ we have~: $f_{v,\theta}^\lambda \equiv
f_{e^{i\theta}v}^\lambda$.
\end{lemma}

\proof Since $f_v^\lambda$ is a canonical disc, the disc
$f_{v,\theta}^\lambda$ has a  lift   close to the lift
of the disc $\zeta \mapsto e^{i\theta}v\zeta$. Then according to Lemma
\ref{circle} 
$f_{v,\theta}^\lambda$ is a canonical disc close to the
linear disc $\zeta \mapsto e^{i\theta}v\zeta$. Since the first jet
of $f_{v,\theta}^\lambda$ coincides with the first jet of
$f_{e^{i\theta}v}^\lambda$, these two nearby stationary discs coincide
according to Lemma~\ref{LEMMA1}. \qed

This statement implies that for any $w \in I^{\lambda}$ the vector
$e^{i\theta}w$ is in $I^{\lambda}$ as well. 

It follows from the
above arguments that there exists a natural parametrization of the set
of canonical discs by their tangent vectors at the origin,
that is by the points of $I^{\lambda}$. The map
$$
\begin{array}{llcll}
\widehat{F^\lambda} &:&  I^\lambda \times \Delta 
& \rightarrow & \mathbb B_n\\
    
    & &  (v,\zeta) & \mapsto     & f_v^\lambda(\zeta)
\end{array}
$$ 
is smooth on $I^\lambda \times \Delta$. Moreover, if we fix a
small positive constant $\varepsilon_0$, then by shrinking $\lambda_0$
if necessary there is, for every $\lambda < \lambda_0$, a smooth
function $\widehat{G^\lambda}$ defined on $I^\lambda \times \Delta$,
satisfying $\|\widehat{G^\lambda(v,\zeta)}\| \leq \varepsilon_0 |\zeta|^2$ on
$I^\lambda \times \Delta$, such that for every $\lambda <
\lambda_0$ we have on $I^\lambda \times \Delta$~:

\begin{equation}\label{equation}
\widehat{F^ \lambda}(v,\zeta) = \zeta v + \widehat{G^\lambda}(v,\zeta).
\end{equation}

\vskip 0,1cm
Consider now the restriction of $\widehat{F^\lambda}$ to $I^\lambda
\times [0,1]$. This is a smooth map, still denoted by 
$\widehat{F^\lambda}$.  We
have the following~:
\begin{proposition}\label{PROPRO}
There exists $\lambda_1 \leq \lambda_0$ such that 
for every $\lambda < \lambda_1$ the family
$(\widehat{F^\lambda}(v,r))_{(v,r) \in I^\lambda \times [0,1[}$
is a real foliation of $\mathbb B_n \backslash \{0\}$. 
\end{proposition}
\proof 

\noindent{\it Step 1}.
For $r \neq 0$ we write $w:=rv$. Then $r=\|w\|$, $v=w/\|w\|$ 
and we denote by $\widetilde{F^\lambda}$ the function 
$\widetilde{F^\lambda}(w):=\widehat{F^\lambda}(v,r)$. 
For $\lambda < \lambda_0$, $\widetilde{F^\lambda}$ is a smooth map 
of the variable $w$ on
$\mathbb B_n \backslash\{0\}$, satisfying~:
$$
\widetilde{F^\lambda}(w) = w + \widetilde{G^\lambda}(w)
$$
where $\widetilde{G^\lambda}$ is a smooth map on $\mathbb B_n \backslash \{0\}$
with $\|\tilde{G}^\lambda(w)\| \leq \varepsilon_0\|w\|^2$ on $\mathbb B_n
\backslash\{0\}$. This implies that $\widetilde{F^\lambda}$ is a local
diffeomorphism at each point in $\mathbb B_n \backslash\{0\}$, and so that
$\widehat{F^\lambda}$ is a local diffeomorphism at each point in
$I^\lambda \times ]0,1[$. Moreover, the
condition $\|\widetilde{G^\lambda}(w)\| \leq \varepsilon_1\|w\|^2$ on
$\mathbb B_n \backslash\{0\}$ implies that $\widetilde{G^\lambda}$ is
differentiable at the origin with $d\widetilde{G^\lambda}(0) = 0$. Hence
by the implicit function theorem there exists $\lambda_1 < \lambda_0$
such that the map $\widetilde{F^\lambda}$ is a local diffeomorphism at the
origin for $\lambda < \lambda_1$. So there exists $0<r_1<1$ and a 
neighborhood $U$ of the origin in $\C^n$ such that
$\widehat{F^\lambda}$ is a diffeomorphism from $I^\lambda \times
]0,r_1[$ to $U \backslash \{0\}$, for $\lambda < \lambda_1$.

\vskip 0,1cm
\noindent{\it Step 2}. We show that $\widehat{F^\lambda}$ 
is injective on $I^\lambda \times ]0,1]$ for sufficiently small $\lambda$. 
Assume by contradiction that for every $n$ there exist 
$\lambda_n \in \mathbb R$, $r_n, r'_n \in ]0,1]$, $v^n, w^n \in I^{\lambda_n}$ 
such that:

$ \bullet \lim_{n \rightarrow \infty}\lambda_n =0, 
\ \lim_{n \rightarrow \infty}r_n=r,\ 
\lim_{n \rightarrow \infty}r'_n=r',
$

$\bullet
\lim_{n \rightarrow \infty}v^n=v \in \mathbb S^{2n-1}, \ 
\lim_{n \rightarrow \infty}w^n=w \in \mathbb S^{2n-1}
$

and satisfying
$$
\widehat{F^{\lambda_n}}(v^n,r_n) = \widehat{F^{\lambda_n}}(w^n,r'_n)
$$
for every $n$. Since $\widehat{F}$ is smooth with respect to $\lambda, v, r$,
it follows that $\widehat{F^0}(v,r) = \widehat{F^0}(w,r')$ and so
$v=w$ and $r=r'$. If $r < r_1$ then the contradiction follows from the fact 
that $\widehat{F^\lambda}$ is a diffeomorphism from 
$I^\lambda \times ]0,r_1[$ to $U \backslash \{0\}$. If $r \geq r_1$
then for every neighborhood $U_\infty$ of $rv$ in 
$\mathbb B_n \backslash \{0\}$, $r_nv^n \in U_\infty$ and 
$r_n'w^n \in U_\infty$ for sufficiently large $n$. Since we may choose 
$U_\infty$ such that $\widehat{F^\lambda}$ is a diffeomorphism from a
neighborhood of $(v,r)$ in $I^\lambda \times ]r_1,1]$ uniformly with respect
to $\lambda <<1$, we still obtain a contradiction.

\vskip 0,1cm
\noindent{\it Step 3}. We show that $\widehat{F^\lambda}$ is surjective
from $I^\lambda \times ]0,1[$ to $\mathbb B_n \backslash \{0\}$. 
It is sufficient to show that $\widehat{F^\lambda}$ is surjective from
$I^\lambda \times [r_1,1[$ to $\mathbb B_n \backslash U$.
Consider the nonempty set 
$E_\lambda=\{w \in \mathbb B_n \backslash U : w=\widehat{F^\lambda}(v,r) \ 
{\rm for \ some \ }(v,r) \in I^\lambda \times ]r_1,1[\}$. Since the jacobian
of $\widehat{F^\lambda}$ does not vanish for $\lambda=0$ and 
$\widehat{F^\lambda}$ is smooth with respect to $\lambda$, 
the set $E_\lambda$ is open for
sufficiently small $\lambda$.
Moreover it follows immediately from its definition that $E_\lambda$ 
is also closed in $\mathbb B_n \backslash U$. Thus 
$E_\lambda=\mathbb B_n \backslash U$.

These three steps prove the result. \qed

\vskip 0,1cm
We can construct now the map $\Psi_{J_\lambda}$ for $\lambda <
\lambda_1$. For every $z \in \mathbb B_n \backslash \{0\}$ consider the unique
couple $(v(z),r(z)) \in I^\lambda \times ]0,1[$ such that
$f_v^\lambda$ is the unique canonical disc passing through $z$ (its
existence and unicity are given by Proposition~\ref{PROPRO}) with
$f_{v(z)}^\lambda(0) = 0$, $df_{v(z)}^\lambda(0)(\partial / \partial
Re(\zeta)) = v(z)$ and $f_{v(z)}^\lambda(r(z)) = z$. The map $\Psi_{J_\lambda}$ is
defined by~:
$$
\begin{array}{llcll}
\Psi_{J_\lambda} &: & \bar{\mathbb B}_n \backslash \{0\} & \rightarrow & \C^n\\
     &  &  z                   & \mapsto     & r(z) v(z).
\end{array}
$$
\begin{definition}
The map $\Psi_{J_\lambda}$ is called the Riemann map associated with the almost
complex structure
$J_{\lambda}$.
\end{definition}
This map is an analogue of the circular representation of a strictly
convex domain introduced by L.Lempert~\cite{le81}. The term ``Riemann map''
was used by S. Semmes~\cite{se92} for a slightly different map where the
vector $v(z)$ is normalized (and so such a map takes values in the unit ball).
In this paper we work with the indicatrix since this is more convenient for
our applications.

The Riemann map $\Psi_{J_\lambda}$ has the following properties~:
\begin{proposition}\label{PROPRO2}
\noindent $(i)$ For every $(v,\zeta) \in I^\lambda \times
\Delta$ we have $(\Psi_{J_\lambda} \circ f_v^\lambda)(\zeta) =
\zeta v$ and so $\log\|(\Psi_{J_\lambda} \circ f_v^\lambda)(\zeta)\| =
\log|\zeta|$.

\noindent $(ii)$ There exist constants $0 < C' < C$ such that $C'
\|z\| \leq \|\Psi_{J_\lambda}(z)\| \leq C \|z\|$ on $\mathbb B_n$.
\end{proposition}

\vskip 0,1cm
\proof $(i)$ Let $\zeta = e^{i\theta}r \in \Delta(0,r_0)$ with $\theta
\in [0,2\pi[$. Then $f_v^\lambda(\zeta) = f_v^\lambda(e^{i\theta}r) =
f_{e^{i\theta}v}^\lambda(r)$. Hence we have $(\Psi_{J_\lambda} \circ
f_v^\lambda)(\zeta) = \Psi_{J_\lambda}( f_{e^{i\theta}v}^\lambda(r)) =
e^{i\theta}vr = \zeta v$.

$(ii)$ Let $z \in \mathbb B_n \backslash \{0\}$. Then according to
equation~(\ref{equation}) we have the inequality
$\|\Psi_{J_\lambda}(z)\|\ (1-\varepsilon_1\|\Psi_{J_\lambda}(z)\|) \leq
\|z\|^2 \leq \|\Psi_{J_\lambda}(z)\| \ (1+ \varepsilon_1
\|\Psi_{J_\lambda}(z)\|)$. Since $\|\Psi_{J_\lambda}(z)\| \leq 1$ we obtain
the desired inequality with $c'=1/1+\varepsilon_1$ and
$c=1/1-\varepsilon_1$. \qed

\vskip 0,1cm
From the above analysis we deduce the following basic properties of the Riemann map.

\begin{proposition}\label{PR1}

\begin{itemize}
\item[]
\item[(i)] The indicatrix $I^{\lambda}$ is a compact circled smooth
$J_{\lambda}$-strictly pseudoconvex hypersurface bounding a domain
denoted by $\Omega^{\lambda}$.
\item[(ii)] The Riemann map $\Psi_{J_\lambda}: \bar {\mathbb B}_n \backslash
\{ 0 \} \rightarrow \bar\Omega^{\lambda} \backslash \{ 0 \}$ is a smooth
diffeomorphism.
\item[(iii)] For every canonical disc $f_v^{\lambda}$ we have
$\Psi_{J_{\lambda}}\circ f^{\lambda}_v(\zeta) = v\zeta$.
\end{itemize}
\end{proposition}

\subsubsection{Local Riemann map}
We introduce the notion of local indicatrix 
$I_{v^0}^\lambda$ for $v^0 \in \mathbb S^{2n-1}$. We may localize the
notion of the Riemann map, introducing a similar associated with the
local indicatrix. Denote by $\Omega_{v^0}^\lambda$ the set 
$I_{v^0}^\lambda \times
[0,1[$. The arguments used in the proof of Proposition~\ref{PROPRO}
show that $\widehat{F^\lambda}(\Omega_{v^0}^\lambda)$ 
is foliated by stationary discs centered at the origin.
We may therefore define the Riemann map $\Psi_{J_\lambda,v^0}$ on 
$\widehat{F^\lambda}(\Omega_{v^0}^\lambda)$ by:
$$
\Psi_{J_\lambda,v^0}(z)=r(z)v(z)
$$
where $v(z)$ is the tangent vector at the origin of the unique stationary disc
$f_{v(z)}^\lambda$ passing through $z$ and $f_{v(z)}^\lambda(r(z)) = z$.

\begin{remark}\label{REM}
We point out that the Riemann map can be defined in any sufficiently small 
deformation of the unit ball and satisfies all the same properties. Moreover
one can easily generalize this construction to strictly convex domains in
$\mathbb C^n$ equipped with small almost complex deformations of the standard
structure.
\end{remark}

\subsubsection{Structure properties of the Riemann map}
Assume now that $M \subset \mathbb C^n$ and let $J_\lambda$ be an almost 
complex deformation of the standard structure on $M$. 
Let $i: T^*_{(1,0)}(M,J)\rightarrow T^*(M)$ be the canonical
identification. Let $D$ be a smoothly bounded
domain in $M$ with the boundary $\Gamma$. The conormal bundle
$\Sigma_J(\Gamma)$ of $\Gamma$ is a real subbundle of
$T^*_{(1,0)}(M,J)|_\Gamma$ whose fiber at $z\in \Gamma$ is defined by
$\Sigma_{z}(\Gamma) = \{ \phi \in T^*_{(1,0)}(M,J): Re \,\phi \vert
H_{(1,0)}^J(\Gamma) = 0 \}$.  Since the form $\partial_J \rho$ forms a
basis in $\Sigma_{J}(\Gamma)$, every $\phi \in \Sigma_J(\Gamma)$ has
the form $\phi = c \partial_J \rho$, $c \in \R$.

\begin{definition}
A continuous map $f:\bar\Delta \rightarrow (\bar
D,J)$, $J$-holomorphic on $\Delta$, is called
a {\it stationary} disc if there
exists a smooth map $\hat{f}= (f,g): \Delta \backslash \{ 0 \}
\rightarrow T^*_{(1,0)}(M,J)$, $\hat{f} \neq 0$ which is continuous on
$\bar \Delta \backslash \{ 0 \}$ and such that
\begin{itemize}
\item[(i)] $\zeta \mapsto \hat{f}(\zeta)$ satisfies the
$\bar\partial_{\tilde J_\lambda}$-equation on $\Delta \backslash \{0\}$,
\item[(ii)] $(i \circ (f,\zeta^{-1}g))(\partial \Delta) \subset
\Sigma_J(\Gamma)$.
\end{itemize}
\end{definition}

We call $\hat{f}$ a {\it lift} of $f$ to the conormal bundle of
$\Gamma$. Clearly the notion
of a stationary disc is {\it invariant} in the following sense: if
$\phi$ is a $\mathcal C^1$ diffeomorphism between $\bar{D}$ and 
$\bar{D}'$ and a $(J,J')$-biholomorphism from
$D$ to $D'$, then for every stationary disc $f$ in $(D,J)$ the
composition $\phi \circ f$ is a stationary discs in $(D',J')$.

Let now $D$ coincide with the unit ball $ {\mathbb B_n}$ equipped with
an almost complex deformation $J_{\lambda}$ of the standard
structure. Then it follows by definition that stationary discs in
$({\mathbb B_n},J_{\lambda})$ may be described as solutions of a
nonlinear boundary problem $({\mathcal BP}_{\lambda})$ associated with
$J_\lambda$. The above techniques give the existence and efficient
parametrization of the variety of
stationary discs in  $({\mathbb B_n},J_{\lambda})$ for $\lambda$ small
enough. This allows to apply the definition of the Riemann map and
gives its existence. We sum up our considerations in the following
Theorem, giving the main structural properties of the Riemann
map.

\begin{theorem}\label{TH1}
Let $J_\lambda$, $J'_{\lambda}$ be almost complex perturbations 
of the standard 
structure on $\bar{\mathbb B}_n$. 
The Riemann map $\Psi_{J_\lambda}$ exists 
for sufficiently small $\lambda$ and
satisfies the following properties:
\begin{itemize}
\item[(a)] $\Psi_{J_\lambda} : \bar{\mathbb B}_n \backslash \{ 0 \} 
\rightarrow \bar
\Omega^\lambda \backslash \{ 0 \}$ is a smooth diffeomorphism
\item[(b)] The restriction of $\Psi_{J_\lambda}$ on every stationary disc
through the origin is $(J_{st},J_{st})$ holomorphic (and even linear)
\item[(c)] $\Psi_{J_\lambda}$ commutes with biholomorphisms. More
precisely for sufficiently small $\lambda'$ and for every
$\CC^1$ diffeomorphism $\varphi$ of $\bar{\mathbb B}_n$, 
$(J_\lambda,J_{\lambda'})$-holomorphic in $\mathbb B_n$
and satisfying $\varphi(0) = 0$, we have
$$
\varphi = ( \Psi_{J_{\lambda'}})^{-1} \circ d\varphi_0 \circ \Psi_{J_\lambda}.
$$
\end{itemize}
\end{theorem}

\vskip 0,2cm
\noindent{\it Proof of Theorem~\ref{TH1}.} Conditions (a) and (b) are 
conditions (i) and (ii) of Proposition~\ref{PR1}. 

\noindent For condition~(c), let $\varphi: ({\mathbb B_n},J) \rightarrow
({\mathbb B_n},J')$ be a $(J,J')$-biholomorphism of class $\CC^1$ on
$\bar{\mathbb B}_n$ satisfying $\varphi(0) = 0$. We know that a disc
$f_v^J$ is a canonical disc for the almost complex structure $J$ if
and only if $\varphi(f_v^J)$ is a canonical disc for the almost
complex structure $J'$. Since $\varphi(f_v^J) =
f_{d\varphi_0(v)}^{J'}$ by definition, $\Psi_J(f_v^J)(\zeta) = \zeta
\,v$ and $\Psi_{J'}(f_{d\varphi_0(v)})(\zeta)=\zeta \,d\varphi_0(v)$ by
Proposition~\ref{PR1} (iii), condition (c) follows from the following
diagram (see Figure 3), which ends the proof of Theorem~\ref{TH1}~:

\bigskip
\begin{center}
\input{diagram3.pstex_t}
\end{center}
\vskip 0,5cm
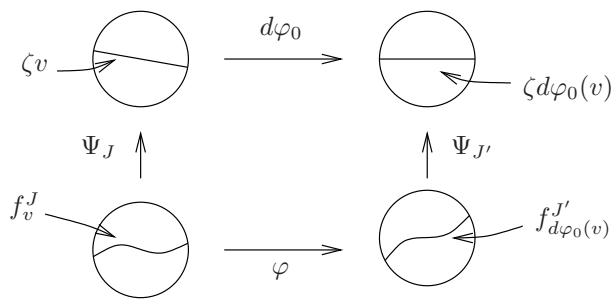
\centerline{Figure 3}

\qed
\vskip 0,5cm

Riemann maps are useful for the boundary study of biholomorphisms in
almost complex manifolds. We have
\begin{corollary}
If $\lambda, \lambda' <<1$ and $\varphi$ is a $\CC^1$ diffeomorphism of 
$\bar{\mathbb B}_n$, 
$(J_\lambda,J_{\lambda'})$-holomorphic in $\mathbb B_n$
satisfying $\varphi(0) = 0$, then $\varphi$ is of class $\CC^{\infty}$
on $\bar{\mathbb B}_n$.
\end{corollary}
\proof This follows immediately by Theorem~\ref{TH1} condition~(c) 
since the Riemann map is smooth up to the boundary. \qed

\subsubsection{Rigidity and local equivalence problem}

Condition $(c)$ of Theorem~\ref{TH1} implies the following partial 
generalization of Cartan's theorem for almost complex manifolds:
\begin{corollary}
If $\lambda <<1$ and if 
$\varphi$ is a $\CC^1$ diffeomorphism of $\bar{\mathbb B}_n$, 
$(J_\lambda,J_{\lambda})$-holomorphic in $\mathbb B_n$,
satisfying $\varphi(0) = 0$ and $d\varphi(0) = I$ then $\varphi$ is the 
identity.
\end{corollary}
This provides an efficient parametrization of the
isotropy group of the group of biholomorphisms of $({\mathbb B_n},J_\lambda)$.

\vskip 0,1cm
We can solve the local biholomorphic 
equivalence problem between almost complex manifolds in terms of the 
Riemann map similarly to~\cite{bl-du-ka87,le88} (see the paper
\cite{li50} by P. Libermann for a traditional approach to this problem 
based on Cartan's equivalence method for $G$-structures).
Let $I^\lambda$ (resp. $(I')^\lambda$) be the indicatrix of 
($\mathbb B_n,J_\lambda$) (resp. ($\mathbb B_n,J'_\lambda$)) bounding
the domain $\Omega^\lambda$ (resp. $(\Omega')^\lambda$) and let 
$\Psi_{J_\lambda}$ (resp. $\Psi_{J'_\lambda}$) be the associated Riemann map.
This induces the almost complex structure 
$J_\lambda^*:=d\Psi_{J_\lambda} \circ J_\lambda \circ d(\Psi_{J_\lambda})^{-1}$
 (resp.
$(J'_\lambda)^*:= 
d\Psi_{J'_\lambda} \circ J_\lambda \circ d(\Psi_{J'_\lambda})^{-1}$) 
on $\Omega^{\lambda}$ (resp. $(\Omega')^{\lambda}$). Then we have:

\begin{theorem}\label{TH3}
The following conditions are equivalent:

$(i)$ There exists a $\CC^\infty$ diffeomorphism $\varphi$ of 
$\bar{\mathbb B}_n$, $(J_\lambda,J_{\lambda}')$-holomorphic on $\mathbb B_n$
and satisfying $\varphi(0)=0$, 

$(ii)$ There exists a $J_{st}$-linear isomorphism
$L$ of $\C^n$, $(J_\lambda^*,(J'_\lambda)^*)$-holomorphic 
on $\Omega^\lambda$ and such that $L(\Omega^\lambda)=(\Omega')^\lambda$.
\end{theorem}

\noindent{\it Proof}. If $\varphi$ satisfies condition $(i)$, then
$L:=d\varphi_0$ satisfies condition $(ii)$, in view of 
the commutativity of the following diagram (see Figure 4) given by
Theorem~\ref{TH1}~:

\bigskip
\begin{center}
\input{riemann2.pstex_t}
\end{center}
\vskip 0,5cm
\centerline{Figure 4}

\vskip 0,5cm

\noindent 
Conversely if $L$ satisfies condition $(ii)$ 
then the map $\varphi:=(\Psi'_{J_\lambda})^{-1} \circ
L \circ \Psi_{J_\lambda}$ satisfies condition $(i)$. \qed
%In conclusion we point out that our main results concerning the existence 
%and the main properties of the Riemann map can be generalized to the case of
%a pointed strictly convex domain $(D,p)$ (where $p$ is a fixed point in $D$)
%and a small deformation $J_\lambda$ ($J_\lambda(p) = J_{st}$) of the standard 
%structure $J_{st}$, using the results of Section~3.2. 

\section{Kobayashi metric on almost complex manifolds}

In this section we give a lower estimate
on the Kobayashi-Royden infinitesimal metric on a strictly pseudoconvex
domain in an almost complex manifold. In particular, we prove that
every point in an almost complex manifold has a basis of complete hyperbolic
neighborhoods. These results were obtained in the paper \cite{ga-su}.

\subsection{Localization of the Kobayashi-Royden metric}

\subsubsection{Kobayashi-Royden infinitesimal pseudometric}
Let $(M,J)$ be an almost complex
manifold.
According to \cite{ni-wo}, for every $p \in M$ there is a
neighborhood $\mathcal V$ of $0$ in $T_pM$ such that for every $v \in
\mathcal V$ there exists $f \in \mathcal O_J(\Delta,M)$ satisfying
$f(0) = p,$ $df(0) (\frac{\partial}{\partial Re(\zeta)}) = v$.
This allows to define
the Kobayashi-Royden infinitesimal pseudometric $K_{(M,J)}$.
\begin{definition}\label{dd}
For $p \in M$ and $v \in T_pM$, $K_{(M,J)}(p,v)$ is the infimum of the
set of positive $\alpha$ such that there exists a $J$-holomorphic disc
$f:\Delta \rightarrow M$ satisfying $f(0) = p$ and $df(0)(\frac{\partial}
{\partial Re(\zeta)}) = v/\alpha$.
\end{definition}
The following statement is an obvious consequence of the above definition~:
\begin{proposition}\label{ppp}
Let $f:(M',J') \rightarrow (M,J)$ be a $(J',J)$-holomorphic map. 
Then $K_{(M,J)}(f(p'),df(p')(v')) \leq K_{(M',J')}(p',v')$ for every
$p'\in M', \ v' \in T_{p'}M'$.
\end{proposition}

We denote by $d_{(M,J)}^K$ the integrated pseudodistance of the
Kobayashi-Royden infinitesimal pseudometric. According to the almost
complex version of Royden's theorem \cite{kr99,iv-ro04}, the
infinitesimal pseudometric is an upper semicontinuous function on the
tangent bundle $TM$ of $M$ and $d^K_{(M,J)}$ coincides with the
usual Kobayashi pseudodistance on $(M,J)$ defined by means of
$J$-holomorphic discs.  Similarly to the case of the integrable
structure we have~:
\begin{definition}\label{dddd} 
$(i)$ Let $p \in M$. Then $M$ is locally hyperbolic at $p$ if 
there exists a neighborhood $U$ of $p$ and a positive constant $C$ such 
that for every $q \in U$, $v \in T_qM$~: $K_{(M,J)}(q,v) \geq C \|v\|$.

$(ii)$ $(M,J)$ is hyperbolic if it is locally hyperbolic
at every point.

$(iii)$ $(M,J)$ is complete hyperbolic if the Kobayashi 
ball $B_{(M,J)}^K(p,r):=\{q \in M : d_{(M,J)}^K(p,q) < r\}$ is relatively 
compact in $M$ for every $p \in M$, $r \geq 0$.
\end{definition}

\begin{lemma}\label{Lemlem}
Let $r < 1$ and let $\theta_r$ be a 
smooth nondecreasing function on
$\R^+$ such that $\theta_r(s)= s$ for $s \leq r/3$ and $\theta_r(s) =
1$ for $s \geq 2r/3$. Let $(M,J)$ be an almost complex manifold, and
let $p$ be a point of $M$. Then there exists a neighborhood $U$ of
$p$, positive constants $A = A(r)$, $B=B(r)$ and a diffeomorphism $z:U
\rightarrow \mathbb B$ such that $z(p) = 0$, $dz(p) \circ J(p) \circ
dz^{-1}(0) = J_{st}$ and the function ${\rm log}(\theta_r(\vert z
\vert^2)) + \theta_r(A\vert z \vert) + B\vert z \vert^2$ is
$J$-plurisubharmonic on $U$.
\end{lemma}

\noindent{\sl Proof of Lemma~\ref{Lemlem}}. 
Denote by $w$ the standard coordinates in $\C^n$. It follows
from Lemma~3 that there exist positive constants $A$ and $\lambda_0$
such that the function ${\rm log}(\vert w \vert^2) + A \vert w \vert$
is $J'$-plurisubharmonic on $\mathbb B$ for every almost complex
structure $J'$, defined in a neighborhood of $\bar{\B}$ in $\C^n$ and
such that $\|J'-J_{st}\|_{\CC^2(\bar{\mathbb B})} \leq \lambda_0$. This
means that the function $v(w) = \log(\theta_r(|w|^2)) + \theta_r(A|w|)$
is $J'$-plurisubharmonic on $B(0,r')=\{w \in \C^n : |w| < r'\}$ for
every such almost complex structure $J'$, where $r'=\inf(\sqrt{r/3},
r/3A)$. Decreasing $\lambda_0$ if necessary, we may assume that the
function $|w|^2$ is strictly $J'$-plurisubharmonic on $\B$. Then,
since $v$ is smooth on $\mathbb B \backslash B(0,r')$, there exists a
positive constant $B$ such that the function $v + B\vert w \vert^2$ is
$J'$-plurisubharmonic on $\mathbb B$ for
$\|J'-J_{st}\|_{\CC^2(\bar{\mathbb B})} \leq \lambda_0$. According to
Lemma \ref{suplem1} there exists a neighborhood $U$ of $p$ and a
diffeomorphism $z:U \rightarrow \mathbb B$ such that $\vert\vert
z_*(J) - J_{st} \vert\vert_{\mathcal C^2(\bar{\mathbb B})} \leq
\lambda_0$. Then the function $v \circ z = {\rm
log}(\theta_r(\vert z \vert^2)) + \theta_r(A \vert z \vert) +
B\vert z \vert^2$
is $J$-plurisubharmonic on $U$. \qed

\begin{proposition}\label{thm3}
(Localization principle) Let $D$ be a domain in an almost complex
manifold $(M,J)$, let $p \in \bar{D}$, let $U$ be a neighborhood of
$p$ in $M$ (not necessarily contained in $D$) and let $z:U \rightarrow
\mathbb B$ be the diffeomorphism given by Lemma~\ref{Lemlem}.  
Let $u$ be a $\mathcal C^2$ function on $\bar{D}$, negative and
$J$-plurisubharmonic on $D$. We assume that $-L \leq u < 0$ on $D \cap
U$ and that $u-c|z|^2$ is $J$-plurisubharmonic on $D \cap U$, where
$c$ and $L$ are positive constants. Then there exist a positive
constant $ s$ and a neighborhood $V \subset \subset U$ of $p$,
depending on $c$ and $L$ only, such that for $q \in D \cap V$ and $v
\in T_qM$ we have the following inequality~:

\begin{equation}\label{E2}
K_{(D,J)}(q,v) \geq s K_{(D \cap U,J)}(q,v).
\end{equation}
\end{proposition}

We note that a similar statement was 
obtained by F.Berteloot \cite{be95} in the integrable case. The proof
is based on N.Sibony's method \cite{si81}.

\vskip 0,2cm
\noindent{\sl Proof of Proposition~\ref{thm3}}.
 Fix a neighborhood $V$ of $p$, relatively compact in $U$.
For every $q \in V$ we consider a diffeomorphism $z_q$ from $U$ to
$\mathbb B$ such that $z_q(q) = 0$ and $(z_q)_*(J)(0) = J_{st}$. We
also may assume that the function $u -c|z_q|^2$ is
$J$-plurisubharmonic on $D \cap U$. 

According to Lemma
~\ref{Lemlem}, there exist uniform positive constants $A$ and $B$ such
that the function 
$$
{\rm log}(\theta_r(|z_q|^2))+ \theta_r(A|z_q|)+ B|z_q|^2
$$
is $J$-plurisubharmonic on $U$ for every $q \in V$. 
Set $\tau=2B/c$ and
define, for every point $q \in V$, the function $\Psi_{q}$ by~:
$$
\left\{
\begin{array}{lll}
\Psi_{q}(z) &=& \theta_r(|z_q|^2)\exp(\theta_r(A|z_q|)) 
\exp(\tau u(z))\ {\rm if} \  z \in D \cap U,\\
& & \\
\Psi_{q} &=& \exp(1+\tau u) \ {\rm on} \ D \backslash U.
\end{array}
\right.
$$

Then for every $0 < \varepsilon \leq B$, the function ${\rm
log}(\Psi_{q})-\varepsilon|z_q|^2$ is $J$-plurisubharmonic on $D \cap U$
and hence $\Psi_{q}$ is $J$-plurisubharmonic on $D \cap U$. Since
$\Psi_{q}$ coincides with $\exp(\tau u)$ outside $U$, it is globally
$J$-plurisubharmonic on $D$. 

Let $f:\Delta \rightarrow D$ be a $J$-holomorphic disc such that
$f(0)=q \in V$ and
$(\partial f/\partial Re(\zeta))(0) = v/\alpha$ where $v \in T_qM$ and $\alpha
>0$. For $\zeta$ sufficiently close to 0 we have
$$
f(\zeta) = q + df(0)(\zeta) +
\mathcal O(|\zeta|^2).
$$
Consider the function
$$
\varphi(\zeta) = \Psi_q(f(\zeta))/|\zeta|^2
$$
which is subharmonic on
$\Delta \backslash \{0\}$. Since
$$
\varphi(\zeta) = |z_q(f(\zeta))|^2/|\zeta|^2 \exp(A|z_q(f(\zeta))|) 
\exp(\tau u(f(\zeta)))
$$
for $\zeta$ close to 0 and 
$$
(z_q\circ f)(\zeta) = dz_q(q)(df(0)(\partial / \partial
  Re(\zeta)))\zeta + \mathcal O(|\zeta|^2)
$$
we obtain that $\limsup_{\zeta \rightarrow 0}\varphi(\zeta)
$ is finite. Moreover
$$
\limsup_{\zeta \rightarrow 0}\varphi(\zeta) \geq \|dz_q(q)(df(0)(\partial
/\partial Re(\zeta))\|^2\exp(-2B|u(q)|/c).
$$
Applying the maximum principle to a subharmonic extension of $\varphi$
on $\Delta$ we obtain the inequality 
$$
\|dz_q(q)(df(0)(\partial / \partial Re(\zeta)))\|^2 \leq \exp(1+2B|u(q)|/c).
$$
Furthermore there exists a positive constant $C$ such that
$$
\|df(0)(\partial / \partial Re(\zeta))\|^2 \leq C
\|dz_q(q)(df(0)(\partial / \partial Re(\zeta)))\|^2.
$$

Hence, by definition of the Kobayashi-Royden infinitesimal pseudometric, 
we obtain for every $q \in D \cap V$, $v \in T_qM$~:
\begin{eqnarray}
\label{localhyp}
K_{(D,J)}(q,v) \geq C^{-1/2} \left(\exp\left(-1-2B\frac{|u(q)|}{c}\right)
\right)^{1/2}\|v\|.
\end{eqnarray}
Consider now the Kobayashi ball $B_{(D,J)}(q,\alpha)=\{w \in D :
d_{(D,J)}^K(w,q)<\alpha\}$. It follows from Lemma~2.2 of \cite{ccs99} (whose
proof is identical in the almost complex setting) that restricting $V$
if necessary we can find a positive constant $s<1$,
independent of $q$, such that for every $J$-holomorphic disc
$f:\Delta \rightarrow D$ satisfying $f(0) \in D \cap V$ we have $f(s\Delta)
\subset D \cap U$ (see Figure 5). This gives the inequality (\ref{E2}). \qed

\bigskip
\begin{center}
\input{figure2.pstex_t}
\end{center}
\bigskip

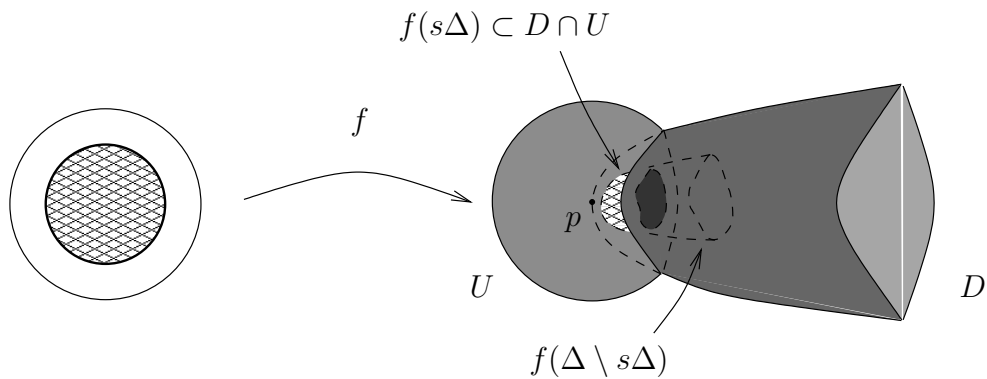
\centerline{Figure 5}
\bigskip

\subsection{Uniform estimates of the Kobayashi-Royden metric}

In the present section we refine the lower estimate of the
Kobayashi-Royden metric.

\begin{proposition}\label{addest}
Let $D$ be a domain in an almost complex
manifold $(M,J)$, let $p \in \bar{D}$, let $U$ be a neighborhood of
$p$ in $M$ (not necessarily contained in $D$) and let $z:U \rightarrow
\B$ be the diffeomorphism given by Lemma~\ref{Lemlem}.  
Let $u$ be a $\mathcal C^2$ function on $\bar{D}$, negative and
$J$-plurisubharmonic on $D$. We assume that $-L \leq u < 0$ on $D \cap
U$ and that $u-c|z|^2$ is $J$-plurisubharmonic on $D \cap U$, where
$c$ and $L$ are positive constants. Then there exists a neighborhood
$U'$ of $p$ and a constant $c'
> 0$, depending on $c$ and $L$ only, such that :

\begin{equation}\label{e1}
K_{(D,J)}(q,v) \geq c'\frac{\|v\|}{|u(q)|^{1/2}},
\end{equation}
for every $q \in D \cap U'$ and every $v \in T_qM$.
\end{proposition}

\noindent{\it Proof of proposition~\ref{addest}.}
We use the notations of the proof of Proposition~\ref{thm3}. 
Consider a positive constant $r$ that will be specified later and
let $\theta$ be a smooth non decreasing function
on $\R^+$ such that $\theta(x) =x$ for $x \leq 1/3$ and $\theta(x) =
1$ for $x \geq 2/3$. Restricting $U$ if necessary we may aasume that
the function 
$
\log(\theta(|z_q/r|^2)) + A|z_q| + B|z_q/r|^2
$
is $J$-plurisubharmonic on $D \cap U$, independently of $q$ and $r$.

Consider now the function 
$
\Psi_{q}(z) =
\theta (|z_q|^2/r^2)\exp(A|z_q|)\exp(\tau u)
$
where $\tau=1/|u(q)|$ and $r=(2B|u(q)|/c)^{1/2}$. 
Since the function $\tau u - 2B |z/r|^2$ is $J$-plurisubharmonic,
we may assume, shrinking $U$ if necessary, that the function
$\tau u - B |z_q/r|^2$ is $J$-plurisubharmonic on $D \cap U$ for every
$q \in V$. Hence the
function $\log(\Psi_q)$ is $J$-plurisubharmonic on $D \cap U$.
It follows from the estimate
(\ref{E2}) that there is a positive constant $s$ such that 
$D_{(D,J)}(q,v) \geq s K_{(D\cap U,J)}(q,v)$ 
for every $q \in V,\ v \in T_qM$. Let $q \in V$, let $v \in T_qM$ and let 
$f:\Delta \rightarrow D \cap U$ be a $J$-holomorphic disc such that
$f(0) = q$ and $df(0)(\partial / \partial Re(\zeta)) = v/\alpha$ where $\alpha
>0$.
Consider the function
$
\varphi(\zeta) = \Psi_q(f(\zeta))/|\zeta|^2
$
which is subharmonic on
$\Delta\backslash \{0\}$. As above we obtain that $\limsup_{\zeta
  \rightarrow 0}\varphi(\zeta)$
is finite and $
\limsup_{\zeta \rightarrow
0}\phi(\zeta) \geq \|v\|^2\exp(2)/(r^2\alpha^2)$.
There exists a positive constant $C'$, independent of $q$, such
that $|z_q| \leq C'$ on $U$.  Applying the maximum
principle to a subharmonic extension of $\varphi$ on $\Delta$, we obtain the
inequality
$$
\alpha \geq \sqrt{\frac{c}{2B\exp(1+AC')}}\|v\|^2/|u(q)|^{1/2}.
$$
This completes the proof. \qed

\subsubsection{Scaling and  estimates of the Kobayashi-Royden metric} 
In this Section we present a precise lower estimate on the Kobayashi-Royden
infinitesimal metric on a strictly pseudoconvex domain in
$(M,J)$.

\begin{theorem}\label{THM}
Let $M$ be a real $2n$-dimensional  manifold with an almost complex
structure $J$ and
let $D=\{\rho<0\}$ be a relatively compact domain in $(M,J)$.
We assume that $\rho$ is a $\CC^2$ defining function of $D$,
strictly $J$-plurisubharmonic in a neighborhood of $\bar{D}$. Then
there exists a positive constant $c$  such that~:

\begin{equation}\label{e3}
K_{(D,J)}(p,v) \geq c\left[\frac{|\partial_J\rho(p)(v - iJ(p)v)|^2}
{|\rho(p)|^2} + 
\frac{\|v\|^2}{|\rho(p)|}\right]^{1/2},
\end{equation}
for every $p \in D$ and every $v \in T_pM$.
\end{theorem}

We start with the small almost complex deformations of the standard
structure. In the second subsection, we consider the case of an
arbitrary almost complex structure, not necessarily close to the
standard one. We use non-isotropic dilations in special coordinates
``reducing'' an almost complex structure in order to represent a
strictly pseudoconvex hypersurface on an almost complex manifold 
as the Siegel sphere equipped with an arbitrary
small deformation of the standard structure. We stress
that such a representation cannot be obtained by the isotropic
dilations of Lemma 1 since the limit hypersurface is just a
hyperplane. 

\subsubsection{Small deformations of the standard structure}

We start the proof of Theorem~\ref{THM} with the following~:

\begin{proposition}\label{thm2}
Let $D=\{\rho < 0\}$ be a bounded domain in $\C^n$, where $\rho$ is a
 $\CC^2$ defining function of $D$, strictly $J_{st}$-plurisubharmonic in
 a neighborhood of $\bar{D}$. Then there exist positive constants $c$
 and $\lambda_0$ such that for every almost complex structure $J$
 defined in a neighborhood of $\bar{D}$ and such that
 $\|J-J_{st}\|_{\CC^2(\bar{D})} \leq \lambda_0$ estimate~(\ref{e3}) is
satisfied for every $p \in D,$ $v \in \C^n$.
\end{proposition}

\vskip 0,1cm
\noindent{\it Proof}. We note that according to Proposition~\ref{thm3}
 (see estimate (\ref{localhyp}))
it is sufficient to prove the inequality 
near $\partial D$. Suppose by contradiction that there exists a
sequence $(p^{\nu})$ of points in $D$ converging to a boundary point
$q$, a sequence $(v^{\nu})$ of unitary vectors and a sequence $(J_\nu)$ 
of almost complex structures defined in a neighborhood of $\bar{D}$, 
satisfying 
$\|J_\nu-J_{st}\|_{\CC^2(\bar{D})}\rightarrow_{\nu \rightarrow \infty} 0$, 
such that the quotient

\begin{equation}
\label{quot1}
K_{(D\cap U,J_\nu)}(p^{\nu},v^{\nu})\left[\frac{|\partial_{J_{\nu}} 
\rho(p^{\nu})(v^{\nu} - iJ_{\nu}(p^{\nu})v^{\nu})|^2}{|\rho(p^{\nu})|^2} 
+ \frac{\|v^{\nu}\|^2}{|\rho(p^{\nu})|}\right]^{-1/2}
\end{equation}
tends to $0$ as $\nu$ tends to $\infty$, where
$U$ is a neighborhood of $q$.
For sufficiently large $\nu$ denote by
$\delta_{\nu}$ the euclidean distance from $p^{\nu}$ to the
boundary of $D$ and by $q^{\nu} \in \partial D$ the unique point such
that $\vert p^{\nu} - q^{\nu} \vert = \delta_{\nu}$. Without loss of
generality we assume that $q = 0$, that $T_0(\partial D) = \{z:=('z,z_n) \in
\C^n : Re(z_n) = 0\}$ and that $J_\nu(q^\nu) = J_{st}$ for every $\nu$.

Consider a sequence of biholomorphic (for the standard structure)
transformations $T^{\nu}$ in a neigborhood of the origin, such that
$T^\nu(q^{\nu}) = 0$ and such that the image $D^{\nu} : =T^\nu(D)$
satisfies
$$
T_0(\partial D^{\nu})=\{ z \in \C^n : Re(z_n) = 0\}.
$$
We point out that the sequence $(T^{\nu})_\nu$ converges uniformly to the
identity map since $q^{\nu} \rightarrow q=0$ as $\nu \rightarrow \infty$
and hence that the sequence $((T^\nu)^{-1})_\nu$ is bounded. We
still denote by $J_\nu$ the direct image
$(T^\nu)_*(J_\nu)$. Let $U_1$ be a neighborhood of the origin such that
$\bar{U} \subset U_1$. For sufficiently large $\nu$ we have 
$T^\nu(U) \subset U_1$. We may assume that every domain $D^{\nu}$ is
defined on $U_1$ by
$$
D^\nu \cap U_1 = \{z \in U_1 : \rho^{\nu}(z) := 
Re(z_n) + |'z|^2 +\mathcal O(|z|^3) <0\},
$$ 
and that the
sequence $(\hat p^\nu = T^{\nu}(p^{\nu}) =(0',-\delta_\nu))_\nu$ is on
the real inward normal to $\partial D^{\nu}$ at 0. Of course, the
functions $\rho^{\nu}$ converge uniformly with all derivatives to the
defining function $\rho$ of $D$. In what follows we omit the hat and
write $p^{\nu}$ instead of $\hat{p}^{\nu}$.

Denote by $R$ the function
$$
{R}(z)=Re(z_n) + |'z|^2 + (Re(z_n) + \vert 'z \vert^2)^2.
$$
There is a neighborhood $V_0$ of the origin in $\mathbb C^n$
such that the function
${R}$ is strictly $J_{st}$-plurisubharmonic on $V_0$. Fix $\alpha > 0$
small enough 
such that the point $z^\alpha=('0,-\alpha)$ belongs to $V_0$.
Consider the dilation $\Lambda_\nu$ defined on
$\C^n$ by $\Lambda_\nu(z) =
({(\alpha / \delta_\nu)^{1/2}}'z,(\alpha/\delta_\nu)z_n)$. 
If we set $J^\nu :=\Lambda_\nu \circ J_\nu \circ (\Lambda_\nu)^{-1}$ 
then we have~: 

\begin{lemma}\label{Lemma1}
 $\lim_{\nu \rightarrow \infty}J^\nu = J_{st}$, uniformly on compact 
subsets of $\C^n$. 
\end{lemma}

\noindent{\it Proof}. Considering $J$ as a matrix valued function,
we may assume that the Taylor expansion of $J_\nu$ at the origin is given by 
$J_\nu = J_{st} + L_\nu(z) + \mathcal O(|z|^2)$ on $U$, uniformly with 
respect to $\nu$. 
Hence ${J}^\nu(z^0)(v) = J_{st}(v) +
{L}_\nu('z,(\delta_\nu/\alpha)^{1/2}z_n)(v)
+ \mathcal O(|(\delta_\nu|) 
\ \|v\|$.
Since $\lim_{\nu \rightarrow \infty}L_\nu = 0$ by assumption, 
we obtain the desired result. \qed

\vskip 0,1cm
Let $\tilde{\rho}^\nu:=(\alpha / \delta_\nu) \rho^\nu \circ 
\Lambda_\nu^{-1}$ and
$G^\nu:=\{z \in \Lambda_\nu(U_1) : \tilde{\rho}^\nu(z) < 0\}$.
Then the function $R^\nu:= \tilde{\rho}^\nu + (\tilde{\rho}^\nu)^2$
converges with all its derivatives to ${R}$, 
uniformly on compact subsets of $\mathbb C^n$.
Hence $R^\nu$ is strictly plurisubharmonic on $V_0$
and according to Lemma~\ref{Lemma1} there is a positive constant $C$ 
such that for sufficiently large $\nu$ the function $R^\nu - C|z|^2$ 
is strictly $J^\nu$-plurisubharmonic on $V_0$.
Since $\sup_{z \in G^\nu \cap \partial V_0} (R^\nu(z) -C|z|^2) =-C'<0$, 
the function 
$$
\tilde{R}^\nu:=\left\{
\begin{array}{lll}
R^\nu - C|z|^2 & {\rm on} & D^\nu \cap V_0\\
 & & \\
-\frac{C'}{2} & {\rm on} & D^\nu \backslash V_0
\end{array}
\right.
$$
is $J^\nu$-plurisubharmonic on $G^\nu$, strictly $J^\nu$-plurisubharmonic
on $G^\nu \cap V_0$. 
Since $z^\alpha$ belongs to $V_0$, it follows from the Proposition~\ref{thm3} 
(see estimate (\ref{localhyp})) that there exists a  positive constant
$C' > 0$ such that for sufficiently large $\nu$ we have~:

$$
K_{(G^\nu,J^\nu)}(z^\alpha,v) \geq  C'\|v\|
$$
for every $v \in \C^n$. 

Moreover for $v \in \C^n$ and for sufficiently large $\nu$ we have~:
\begin{eqnarray*}
& &   
K_{(D^\nu\cap U_1, J_\nu)}(p^\nu,v) =
K_{(G^\nu,J^\nu)}(z^\alpha,\Lambda_\nu(v)) \geq  C'\parallel \Lambda_\nu(v)
\parallel.
\end{eqnarray*}
This gives the inequality~:
$$
K_{(D^{\nu} \cap U,J_\nu)}(p^\nu,v) \geq C' \left(
\frac{\alpha |v_1|^2}{\delta_\nu} + \cdots + 
\frac{\alpha |v_{n-1}|^2}{\delta_\nu} + 
\frac{\alpha^2 |v_n|^2}{\delta_\nu^2}\right)^{1/2}.
$$
Since $C_1\delta_\nu$ is equivalent to $|\rho(p^\nu)|$ as $\nu
\rightarrow \infty$, we obtain that there
is a positive constant $C''$ such that
$$
K_{(D^{\nu} \cap U,J_\nu)}( p^\nu,v) \geq C'' \left(
\frac{\|v\|^2}{|\rho(p^\nu)|} + 
\frac{|v_n|^2}{|\rho(p^\nu)|^2}\right)^{1/2}.
$$
Since $J_{\nu}(0) = J_{st}$, we have  $|\partial\rho(p^\nu)(v - 
iJ_\nu(p^\nu)(v))|^2 = |\partial_{J_{st}}\rho(p^\nu)(v)|^2 +
\mathcal O(\delta_{\nu})\parallel v \parallel^2 = \vert v_n \vert^2 
 + \mathcal O(\delta_\nu)\parallel v \parallel^2$.
Hence there exists a positive constant $\tilde{C}$ such that 
$$
K_{(D^{\nu} \cap U,J_\nu)}( p^\nu,v) \geq \tilde{C} \left(
\frac{\|v\|^2}{|\rho(p^\nu)|} + 
\frac{|\partial_J\rho(p^\nu)(v-iJ_\nu(p^\nu)(v))|^2}
{|\rho(p^\nu)|^2}\right)^{1/2},
$$
contradicting the assumption on the quotient 
(\ref{quot1}). This proves the desired estimate. \qed

We  have the following corollary~:
\begin{corollary}\label{cor3}
Let $(M,J)$ be an almost complex manifold. Then every $p \in M$ has a 
basis of complete hyperbolic neighborhoods.
\end{corollary}

\proof Let $p \in M$. According to Example~1 there exist a
neighborhood $U$ of $p$ and a diffeomorphism $z:U \rightarrow \B$,
centered at $p$, such that the function $|z|^2$ is strictly
$J$-plurisubharmonic on $U$ and 
$\|z_\star(J)-J_{st}\|_{\mathcal C^2(U)} \leq \lambda_0$. 
Hence the open ball $\{x \in \mathbb C^n :
\|x\|<1/2\}$ equipped with the structure  $z_\star(J)$ satisfies the
hypothesis of Theorem \ref{thm2}. Now the estimate on the
Kobaysahi-Royden metric given by this theorem implies that this ball
is  complete hyperbolic by the standard integration argument. 
\qed  

\subsubsection{Arbitrary  almost complex structures}

We turn now to the proof of Theorem~\ref{THM} on an
arbitrary  strictly pseudoconvex domain in an almost
complex manifold $(M,J)$ ($J$ is not supposed to
be a small deformation of the standard structure).
In view of Proposition~\ref{thm3} it suffices to prove
the statement in a neighborhood $U$ of  a boundary point  $q \in \partial D$.
Considering local coordinates $z$ centered at $q$, we may assume that
$D \cap U$ is a domain in $\mathbb C^n$ and
$0 \in \partial D$, $J(0) = J_{st}$.
The idea of the proof is to reduce the situation to the case of a small
deformation of the standard structure considered in Proposition~\ref{thm2}.
In the case of real dimension four Theorem~\ref{THM} is a direct corollary of
Proposition~\ref{thm2}. In the case of arbitrary dimension the proof of
Theorem~\ref{THM} requires a slight modification of Proposition~\ref{thm2}.
So we treat this case seperately.

\subsubsection{Case where $dim M = 4$}
According to \cite{si} Corollary~3.1.2,
there exist a neighborhood $U$ of $q$ in $M$ and complex coordinates
$z=(z_1,z_2) : U \rightarrow \mathbb B_2 \subset \mathbb C^2$, $z(0) =
0$ such that $z_*(J)(0) = J_{st}$ and moreover, a map $f: \Delta
\rightarrow \mathbb B$ is $J':= z_*(J)$-holomorphic if it satisfies the
equations 

\begin{eqnarray}
\label{Jhol1}
\frac{\partial f_j}{\partial \bar \zeta} =
A_j(f_1,f_2)\overline{\left ( \frac{\partial f_j}
{\partial \zeta}\right ) }, j=1,2
\end{eqnarray} 
where $A_j(z) =  O(\vert
z \vert)$, $j=1,2$.

As pointed out before, one can obtain such coordinates by considering
two transversal
foliations of $\mathbb B$ by $J'$-holomorphic curves
and then taking these curves into the lines $z_j = const$ by a local
diffeomorphism (see Figure 1). The direct image of the almost complex structure $J$
under such a diffeomorphism has a diagonal matrix $ J'(z_1,z_2) =
(a_{jk}(z))_{jk}$ with $a_{12}=a_{21}=0$ and $a_{jj}=i+\alpha_{jj}$
where $\alpha_{jj}(z)=\mathcal O(|z|)$ for $j=1,2$.
We point out that the lines $z_j = const$ are
$J$-holomorphic after a suitable parametrization (which, in general,
is not linear). 

\vskip 0,1cm
In what follows we omit the prime and denote this structure again by
$J$. We may assume that the complex tangent space $T_0(\partial D)
\cap J(0) T_0(\partial D) = T_0(\partial D) \cap i T_0(\partial D)$ is
given by $\{ z_2 = 0 \}$.
In particular, we have the following expansion for the defining
function $\rho$ of $D$ on $U$~:
$\rho(z,\bar{z}) = 2 Re(z_2) + 2Re K(z) + H(z) + \mathcal O(\vert z
\vert^3)$, where
$K(z)  = \sum k_{\nu\mu} z_{\nu}{z}_{\mu}$, $k_{\nu\mu} =
k_{\mu\nu}$ and 
$H(z) = \sum h_{\nu\mu} z_{\nu}\bar z_{\mu}$, $h_{\nu\mu} =
\bar h_{\mu\nu}$.

\vskip 0,1cm
Consider the non-isotropic dilations $\Lambda_{\delta}: (z_1,z_2) \mapsto
(\delta^{-1/2}z_1,\delta^{-1}z_2) = (w_1,w_2)$ with $\delta > 0$. 
If $J$ has the above
diagonal form in the coordinates $(z_1,z_2)$ in $\mathbb C^2$, then
its direct image  $J_{\delta}= (\Lambda_{\delta})_*(J)$ has the form
$J_{\delta}(w_1,w_2) =(a_{jk}(\delta^{1/2}w_1,\delta w_2))_{jk}$
and so $J_{\delta}$ tends to $J_{st}$ in the $\mathcal C^2$ norm as $\delta
\rightarrow 0$. On the other hand, $\partial D$ is, in the coordinates
$w$, the zero set of the function 
$\rho_{\delta}= \delta^{-1}(\rho \circ \Lambda_{\delta}^{-1})$.
As $\delta \rightarrow 0$, the function $\rho_{\delta}$ tends to 
the function $2 Re w_2 + 2 Re K(w_1,0) + H(w_1,0)$ which defines a
$J_{st}$- strictly pseudoconvex domain by Lemma~\ref{PP}.
So we may apply Proposition~\ref{thm2}.
This proves Theorem~\ref{THM} in dimension 4.

\subsubsection{Case where $dim M = 2n$.} 
In this case we cannot apply directly Proposition \ref{thm2} since
$J$ can not be deformed by the non-isotropic dilations to the standard
structure. Instead we use the invariance of the Levi form with respect
to the non-isotropic dilations.

We suppose that in a neighborhhod of the origin we have $J = J_{st} +
{\mathcal O}(\vert z \vert)$.
We also may assume that in these coordinates the defining function
$\rho$ of $D$ has the form $\rho = 2Re z_n + 2ReK(z) + H(z) + 
\mathcal O(\vert z \vert^3)$, where $K$ and $H$ are defined similarly
to the 4-dimensional case and $\rho$ is strictly $J$-plurisubharmonic
at the origin. We use the notation $z = ('z,z_n)$.

Consider the non-isotropic dilations $\Lambda_{\delta} : ('z,z_n)
\mapsto (w',w_n) = {(\delta^{-1/2}}'z,\delta^{-1}z_n)$ and set
$J_{\delta} = (\Lambda_{\delta})_*(J)$. Then $J_{\delta}$ tends to the
almost complex structure $J_0(z)= J_{st} + L('z,0)$ where 
$L('z,0) = (L_{kj}('z,0))_{kj}$
denotes a matrix with $L_{kj} = 0$ for $k = 1,...,n-1$, $j = 1,...,n$,
$L_{nn} = 0$ and $L_{nj}('z,0)$, $j=1,...,n-1$ being (real) linear
forms in $'z$.  

Let $\rho_{\delta} = \delta^{-1}(\rho \circ \Lambda_{\delta}^{-1})$.
As $\delta \rightarrow 0$, the function $\rho_{\delta}$ tends to 
the function $\tilde{\rho} = 2 Re z_n + 2 Re K('z,0) + H('z,0)$ in the
$\mathcal C^2$ norm. By the invariance of the Levi form we have
${\mathcal L}^J(\rho)(0)(\Lambda_{\delta}^{-1}(v)) = {\mathcal
  L}^{J_\delta}(\rho \circ \Lambda_{\delta}^{-1})(0)(v)$. Since $\rho$ is
strictly $J$-plurisubharmonic, multiplying by $\delta^{-1}$ and
passing to the limit at the right side  as $\delta \longrightarrow 0$ , 
we obtain that
${\mathcal L}^{J_0}(\tilde \rho)(0)(v) \geq 0$ for any $v$. Now let $v =
('v,0)$. Then $\Lambda_{\delta}^{-1}(v) = \delta^{1/2}v$ and so 
${\mathcal L}^J(\rho)(0)(v) = {\mathcal
  L}^{J_\delta}(\rho_{\delta})(0)(v)$. Passing to the limit as $\delta$
tends to zero, we obtain that  ${\mathcal L}^{J_0}(\tilde \rho)(0)(v) > 0$
for any $v = ('v,0)$ with $'v \neq 0$.

Consider now the function $R=\tilde{\rho} + \tilde{\rho}^2$. 
Then ${\mathcal L}^{J_0}(R)(0)(v) = {\mathcal L}^{J_0}(\tilde \rho)(0)(v)
+ 2 v_n \overline v_n$, so $R$ is strictly $J_0$-plurisubharmonic in a
neighborhood of the origin.
Thus the functions $R^{\nu}$ used in the proof of
Proposition~\ref{thm2} are strictly $J^\nu$-plurisubharmonic and their Levi
forms are bounded from below by a positive constant independent of $\nu$.
This allows to use Proposition~\ref{thm3} and the proof can be
proceeded quite similarly to the proof of Proposition~\ref{thm2}
without any changes.  \qed
    
\subsubsection{Upper estimate of the Kobayashi-Royden metric}

In this subsection we prove the following~:
\begin{theorem}\label{THEOREM}
Let $M$ be a real $2n$-dimensional  manifold with an almost complex
structure $J$ and
let $D=\{\rho<0\}$ be a relatively compact domain in $(M,J)$.
We assume that $\rho$ is a $\CC^2$ defining function of $D$,
strictly $J$-plurisubharmonic in a neighborhood of $\bar{D}$. Then
there exists a positive constant $c$  such that~:

\begin{equation}\label{E3}
K_{(D,J)}(p,v) \leq c\left[\frac{|\partial_J\rho(p)(v - iJ(p)v)|^2}
{|\rho(p)|^2} + 
\frac{\|v\|^2}{|\rho(p)|}\right]^{1/2},
\end{equation}
for every $p \in D$ and every $v \in T_pM$.
\end{theorem}
Since the estimates are purely local we may
assume that $D \subset \mathbb C^n$.
Let $p \in D$ be sufficiently close to $\partial D$ and let $\|v\| =
1$.
 We choose coordinates $z=(z^1,\dots,z^n)$ on $\mathbb C^n$ such that
 $0 \in \partial D$, $p=(0',-\delta)$ where $0 < \delta < 1$,
$dist(p,\partial D) = dist(0,p)$, and
 $D' :=\{z \in B(0,2) / Re(z_n) + |z|^2 <0\} \subset D \cap B(0,2)$.
Moreover we may assume that the map $f_v:\Delta \rightarrow \mathbb
C^n$ defined by $f_v(\zeta) = \zeta v$ is $J$-holomorphic. We denote by
$d$ the distance from $p$ to $\partial D'$ along the
line $\mathbb C v$. Since $f_v(d\Delta)$ is
contained in $D'$ we have the following inequalities~:
$$
K_{(D,J)}(p,v) \leq K_{(D',J)}(p,v) \leq \frac{1}{d}.
$$

Moreover by the strict convexity of $D'$ we have~:
$$
d \geq \frac{1}{2}\left(\delta ^2 \|v_n\|^2 + \delta
  \|v'\|^2|\right)^{1/2}.
$$ 
On the other hand
$$
\frac{1}{(\delta^2 \|v_n\|^2 + \delta \|v'\|^2|)^{1/2}} \leq
\left(\frac{\|v_n\|^2}{\delta ^2} +
\frac{\|v'\|^2}{\delta}\right)^{1/2}
$$
which implies the desired statement. \qed

\subsection{ Boundary continuity and localization of biholomorphisms}
 
In this section we give some technical results necessary for the
proof of the Fefferman mapping theorem. They are also of independent
interest.

\subsubsection{Boundary continuity of diffeomorphisms}
Using estimates of the
Kobayashi-Royden metric together with the boundary distance preserving
property, we obtain, by means of classical arguments
(see, for instance, K.Diederich-J.E.Fornaess \cite{df77}), the following

\begin{theorem}\label{Reg}
Let $D$ and $D'$ be two  smoothly relatively compact strictly pseudoconvex
domains in almost complex manifolds $(M,J)$ and $(M',J')$ respectively. Let
$f:  D  \rightarrow  D'$ be  a smooth
diffeomorphism biholomorphic with respect to $J$ and $J'$. 
Then $f$ extends as a $1/2$-H{\"o}lder homeomorphism
between the  closures of $D$ and $D'$.
\end{theorem}

We recall the following estimates
of the Kobayashi-Royden infinitesimal metric obtained previously~:
\begin{lemma}
\label{lowest1}
Let $D$ be a relatively compact strictly pseudoconvex domain in an
almost complex manifold $(M,J)$. Then there
exists a constant $C > 0$ such that

\begin{eqnarray*}
(1/C)\| v \| / dist(p,\partial D)^{1/2} \leq K_{(D,J)}(p,v) 
\leq C\| v \|/dist(p,\partial D)
\end{eqnarray*}for every $p \in D$ and
$v \in T_pM$.
\end{lemma}

\noindent{\it Proof of Theorem~\ref{Reg}}. For any $p \in D$ and any
tangent vector $v$ at $p$ we have by Lemma~\ref{lowest1}~:

\begin{eqnarray*}
C_1\frac{\| df_p(v) \|}{dist(f(p),\partial D')^{1/2}} \leq 
K_{(D',J')}(f(p),df_p(v)) = K_{(D,J)}(p,v) \leq C_2\frac{\| v
\|}{dist(p,\partial D)}
\end{eqnarray*}
which implies, by Proposition~\ref{equiv}, the estimate
$$\vert\vert\vert df_p \vert\vert\vert \leq C\frac{\| v
\|}{dist(p,\partial D)^{1/2}}.$$
This gives the desired statement. \qed

\vskip 0,2cm
\noindent Theorem~\ref{Reg} allows to reduce the proof of Fefferman's
theorem to a {\it local situation}. 
Indeed, let $p$ be a boundary point of $D$ and
$f(p) = p' \in \partial D'$. It suffices to prove that $f$ extends
smoothly to a neighborhood of $p$ on $\partial D$. Consider
coordinates $z$ and $z'$ defined in small neighborhoods $U$ of $p$ and $U'$ of
$p'$ respectively, with $U' \cap D' = f(D \cap U)$ (this is possible since
$f$ extends as a homeomorphism at $p$). We obtain the following situation.
If $\Gamma = z(\partial D \cap U)$ and $\Gamma' = z'(\partial D' \cap U')$
then the map $z' \circ f \circ z^{-1}$ is defined on $z(D \cap U)$ in $\C^2$,
continuous up to the hypersurface $\Gamma$ with $f(\Gamma) \subset \Gamma'$.
Furthermore the map $z' \circ f \circ z^{-1}$ is a diffeomorphism between
$z(D \cap U)$ and $z'(D' \cap U')$ and the hypersurfaces $\Gamma$ and
$\Gamma'$ are strictly pseudoconvex for the structures $z_*(J)$ and
$(z')_*(J')$ respectively. Finally, we may choose $z$ and $z'$ such that
$z_*(J)$ and $z'_*(J')$ are represented by diagonal matrix functions in the
coordinates $z$ and $z'$. 
As we proved in Lemma~\ref{PP}, $\Gamma$ (resp. $\Gamma'$) is also strictly
$J_{st}$-psdeudoconvex at the origin. We call such
coordinates $z$ (resp. $z'$) {\it canonical coordinates} at $p$
(resp. at $p'$). Using the non-isotropic
dilation as in Section 2.5, we may assume that the norms
$\| z_*(J) - J_{st}\|_{\CC^2}$ and
$\| z'_*(J') - J_{st}\|_{\CC^2}$ are as small as needed.
This localization is crucially used in the sequel and we write $J$ (resp.
$J'$) instead of $z_*(J)$ (resp. $z'_*(J')$); we 
identify $f$ with $z' \circ f \circ z^{-1}$.

\subsubsection{H{\"o}lder extension of holomorphic discs}
We study the boundary continuity of pseudoholomorphic discs
attached to smooth totally real submanifolds in almost complex manifolds.

Recall that in the case of the integrable structure every smooth totally
real submanifold  $E$ (of maximal dimension) is the zero set of a positive
strictly plurisubharmonic function of class $\CC^2$. This remains true
in the almost complex case. Indeed, we can choose coordinates
$z$ in a neighborhood $U$ of $p \in E$ such that $z(p) = 0$, $z_*(J) =
J_{st} + O(\vert z \vert)$ on $U$ and 
$z(E \cap U) = \{w=(x,y) \in z(U) : r_j(w) = x^j +o(\vert
w \vert) = 0 \}$. The function $\rho = \sum_{j=1}^n r_j^2$ is strictly
$J_{st}$-plurisubharmonic on $z(U)$ and so remains strictly
$z_*(J)$-plurisubharmonic, restricting $U$ if necessary.
Covering $E$ by such neighborhoods, we conclude by mean of the partition of
unity.

Let $f : \Delta \rightarrow M$ be a $J$-holomorphic disc and let
$\gamma$ be an open arc on the unit circle $\partial \Delta$.
As usual we denote by $C(f,\gamma)$ the cluster set of $f$ on $\gamma$;
this consists of points $p \in M$ such that $p=\lim_{k \rightarrow
\infty}f(\zeta_k)$ for a sequence $(\zeta_k)_k$ in $\Delta$ converging
to a point in $\gamma$.

\begin{theorem}
\label{Regth1}
Let $G$ be a relatively compact domain in an almost complex manifold $(M,J)$
and let $\rho$ be a strictly $J$-plurisubharmonic
function of class $\CC^2$ on $\bar{G}$.
Let $f:\Delta \rightarrow G$ be a
$J$-holomorphic disc such that $\rho \circ f \geq 0$ on $\Delta$.
Suppose that $\gamma$ is an open non-empty arc on
$\partial \Delta$ such that the cluster set
$C(f,\gamma)$ is contained in the zero set of $\rho$.
Then $f$ extends as a H{\"o}lder 1/2-continuous map on $\Delta \cup
\gamma$.
\end{theorem}

We begin the proof by the following well-known assertion
(see, for instance, \cite{BeL}).

\begin{lemma}
\label{dlem3.1}
Let $\phi$ be a positive subharmonic function in $\Delta$ such that
the measures $\mu_r(e^{i\theta}) := \phi(re^{i\theta})d\theta$ converge in
the weak-star topology to
a measure $\mu$ on $\partial \Delta$ as $r \rightarrow 1$. Suppose that
$\mu$ vanishes on an open arc $\gamma \subset \partial \Delta$. Then for
every compact subset $K \subset \Delta \cup \gamma$ there exists a constant
$C>0$ such that
$\phi(\zeta) \leq C(1 - \vert \zeta \vert)$ for any
$\zeta \in K \cup \Delta$.
\end{lemma}
 
Now fix a point $a \in \gamma$, a constant $\delta > 0$ small
enough so that the intersection $\gamma \cap (a + \delta
\bar\Delta )$ is compact in $\gamma$; we denote by
$\Omega_{\delta}$ the intersection $\Delta \cap (a + \delta\Delta
)$. By Lemma~\ref{dlem3.1}, there exists a constant $C > 0$ such that,
for any $\zeta$ in $\Omega_{\delta}$, we have

\begin{eqnarray}
\label{dd4}
\rho \circ f(\zeta) \leq C (1 - \vert \zeta \vert ).
\end{eqnarray}

Let $(\zeta_k)_k$ be a sequence of points in $\Delta$ converging to $a$
with $\lim_{k \rightarrow \infty}f(\zeta_k) = p$.
By assumption, the function $\rho$ is strictly $J$-plurisubharmonic in a
neighborhood $U$ of $p$; hence there is a constant
$\varepsilon > 0$ such that the function $\rho - \varepsilon \vert z
\vert^2$ is $J$-plurisubharmonic on $U$.

\begin{lemma}
\label{dlem3.2}
There exists a constant $A > 0$ with the following property~: If
$\zeta$ is an arbitrary point of $\Omega_{\delta/2}$ such that
$f(\zeta)$ is in $G \cap z^{-1}(\B)$, then
$\vert \vert \vert df_\zeta \vert \vert \vert
\leq A(1- \vert \zeta \vert)^{-1/2}$.
\end{lemma}
\noindent{\it Proof of Lemma~\ref{dlem3.2}.}
Set $d = 1 - \vert \zeta \vert$; then the disc $\zeta +
d\Delta$ is contained in $\Omega_{\delta}$. Define the domain $G_d =
\{ w \in G: \rho(w) < 2Cd \}$. Then it follows by (\ref{dd4}) that the
image $f(\zeta + d\Delta)$ is contained in $G_d$, where the
$J$-plurisubharmonic function $u_d = \rho - 2Cd$ is negative. 
Moreover we have lower estimates on the Kobayashi-Royden
infinitesimal pseudometric given by Proposition~\ref{addest}.

Hence there exists a positive constant $M$
(independent of $d$) such that $K_{(G_d,J)}(w,\eta) \geq M \vert \eta
\vert \vert u_d(w) \vert^{-1/2}$, for any $w$ in $G \cap z^{-1}(\B)$ and
any $\eta \in T_{w}\Omega$. On another hand, we have $K_{\zeta +
d\Delta}(\zeta,\tau ) = \vert \tau \vert /d$ for any $\tau$ in
$T_{\zeta}\Delta$ indentified with $\C$. By the decreasing property
of the Kobayashi-Royden metric, for any $\tau$ we have 

\begin{eqnarray*}
M \| df_\zeta(\tau)\| \ \vert u_d(f(\zeta))\vert^{-1/2} \leq
K_{(G_d,J)}(f(\zeta),df_\zeta(\tau)) \leq K_{\zeta + d\Delta}(\zeta,\tau) =
\vert \tau \vert/d.
\end{eqnarray*}
Therefore, $\vert \vert \vert df_\zeta\vert\vert\vert \leq M^{-1}\vert
u_d(f(\zeta))\vert^{1/2}/d$. As $-2Cd \leq u_d(f(\zeta)) < 0$, this
implies the desired statement in Lemma~\ref{dlem3.2}
with $A = M^{-1}(2C)^{1/2}$. \qed

\vskip 0,1cm
\noindent{\it Proof of Theorem~\ref{Regth1}}.
Lemma~\ref{dlem3.2} implies that $f$ extends  as a 1/2-H{\"o}lder map
to a neighborhood of the point $a$
in view of an integration argument inspired by the classical
Hardy-Littlewood theorem. 
This proves Theorem~\ref{Regth1}. \qed

\section{Nonisotropic scaling of almost complex manifolds with boundary}

We consider a biholomorphism $f$ between two relatively compact,
strictly pseudoconvex domains $D$ and $D'$ in almost complex manifolds
$(M,J)$ and $(M',J')$.
The aim of this Section is to provide a precise information about the boundary
behavior of the tangent map of $f$.
For convenience of the reader, we begin this Section with the
case of four real dimensional almost complex manifolds (subsections
4.1 and 4.2). The case of higher dimension will be treated in
Subsection 4.3.

\subsection{Localization and boundary behavior of the tangent map} 

Our considerations being purely local, we can assume that $D$ and $D'$ are
domains in $\C^2$, $\Gamma$ and $\Gamma'$ are open
smooth pieces of their boundaries containing the origin, the almost complex
structure $J$ (resp. $J'$) is defined in a neighborhood of $D$
(resp. $D'$), $f$ is a $(J,J')$ biholomorphism from $D$ to $D'$,
continuous up to $\Gamma$, $f(\Gamma) = \Gamma'$, $f(0) = 0$.
The matrix $J$ (resp. $J'$) is diagonal on $D$ (resp. $D'$).

\vskip 0,1cm
Consider a basis $(\omega_1,\omega_2)$ of $(1,0)$ differential forms
(for the structure $J$) in a neighborhood of the origin. Since $J$ is
diagonal, we may choose $\omega_j = dz^j - B_{j}(z)d\bar z^j$, $j=1,2$.
Denote by $Y=(Y_1,Y_2)$ the corresponding dual basis
of $(1,0)$ vector fields. Then $Y_j = \partial /\partial z^j -
\beta_j(z)\partial/\bar\partial z^j$, $j=1,2$. Here $\beta_j(0) =
\beta_k(0) = 0$. The basis $Y(0)$ simply coincides with the canonical (1,0)
basis of $\C^2$.
In particular $Y_1(0)$ is a basis vector of the holomorphic tangent space
$H^J_0(\partial D)$ and $Y_2(0)$ is ``normal'' to $\partial D$.
Consider now for $t \geq 0$ the translation $\partial D -
t$ of the boundary of $D$ near the origin. Consider, in a neighborhood of the
origin, a $(1,0)$ vector field $X_1$ (for $J$) such that $X_1(0) = Y_1(0)$
and $X_1(z)$ generates the complex tangent space $H^J_z(\partial D - t)$ at
every point $z \in \partial D - t$, $0 \leq t <<1$.
Setting $X_2 = Y_2$, we obtain a basis of vector fields
$X = (X_1,X_2)$ on $D$ (restricting $D$ if necessary). 
Any complex tangent vector $v \in T_z^{(1,0)}(D,J)$ at
point  $z \in D$ admits the unique
decomposition $v = v_t + v_n$ where $v_t = \alpha_1
X_1(z)$ (the tangent component) and $v_n = \alpha_2 X_2(z)$ (the normal
component). Identifying $T_z^{(1,0)}(D,J)$ with $T_zD$ we may
consider the decomposition $v=v_t + v_n$ for $v \in T_z(D)$.
Finally we consider this decomposition for points $z$ in a neighborhood of
the boundary.

 We fix a (1,0) basis vector fields $X$
(resp. $X'$) on $D$ (resp. $D')$ as above.

\begin{proposition}
\label{matrix}
The matrix $A = (A_{kj})_{k,j= 1,2}$ of the differential $df_z$ with respect
to the bases
$X(z)$ and $X'(f(z))$ satisfies the following estimates~:
$A_{11} = O(1)$, $A_{12}= O(dist(z,\partial D)^{-1/2})$,
$A_{21}= O(dist(z,\partial D)^{1/2})$ and $A_{22}= O(1)$.   
\end{proposition}
\vskip 0,1cm
\noindent{\it Proof of Proposition~\ref{matrix}}.
Consider the case where $v = v_t$. It follows from
the estimates of the Kobayashi-Royden metric obtained previously that~:
$$
\begin{array}{llcll}
\displaystyle \frac{1}{C}\left (
\frac{\| (df_z(v_t))_t \|}{dist(f(z),\partial
D')^{1/2}} + \frac{\|(df_z(v_t))_n
\|}{dist(f(z),\partial D')} \right ) & \leq & 
K_{(D',J')}(f(z),df_z(v_t)) && \\
& = & K_{(D,J)}(z,v_t) & \leq & \displaystyle C
\frac{\| v_t \|}{dist(z,\partial D)^{1/2}}.
\end{array}
$$
This implies that 
$\|(df_z(v_t))_t \| \leq C^{5/2} \| v_t \|$
and
$\vert\vert (df_z(v_t))_n \vert\vert \leq C^{3}dist(z,\partial D)^{1/2}
\| v_t \|$, by the boundary distance
preserving property given in Proposition~\ref{equiv}.
We obtain the estimates for the normal component in a similar way.
\qed

\subsection{Non isotropic dilations}

Our goal now is to prove Fefferman's mapping theorem without the
assumption of $\CC^1$-smoothness of $f$ up to the boundary. This requires
an application of the estimates of the Kobayashi-Royden metric given
above and the scaling method due to S.Pinchuk; we adapt this to the
almost complex case.

We reduced the problem to the following local
situation. Let $D$ and $D'$ be domains in $\C^2$, $\Gamma$ and
$\Gamma'$ be open $\CC^{\infty}$-smooth pieces of their boundaries,
containing the origin. We assume that an almost complex structure $J$
is defined and $\CC^{\infty}$-smooth in a neighborhood of the closure
$\bar D$, $J(0) = J_{st}$ and $J$ has a diagonal form in a
neighborhood of the origin: $J(z) = diag(a_{11}(z),a_{22}(z))$.
Similarly, we assume that $J'$
is diagonal in a neighborhood of the origin, $J'(z) =
diag(a_{11}'(z),a_{22}'(z))$ and $J'(0) = J_{st}$. The hypersurface
$\Gamma$ (resp. $\Gamma'$) is supposed to be strictly $J$-pseudoconvex
(resp. strictly $J'$-pseudoconvex). Finally, we assume that $f: D
\rightarrow D'$ is a $(J,J')$-biholomorphic map, $1/2$-H{\"o}lder
homeomorphism between  $D \cup \Gamma$ and $D' \cup \Gamma'$, such that
$f(\Gamma) = \Gamma'$ and $f(0) = 0$. Finally,
$\Gamma$ may be defined in a neighborhood of the origin
by the equation $\rho(z) = 0$ where $\rho(z) = 2Re
z^2 + 2Re K(z) + H(z) + o(\vert z \vert^2)$ and $K(z) = \sum
K_{\mu\nu}z^{\mu\nu}$, $H(z) = \sum h_{\mu\nu}z^{\mu}\bar
z^{\nu}$, $k_{\mu\nu} = k_{\nu\mu}$, $h_{\mu\nu} = \bar
h_{\nu\mu}$. The crucial point is that $H(z^1,0)$ is a positive
Hermitian form on $\C$, meaning that in these coordinates $\Gamma$ is
strictly pseudoconvex at the origin with respect to the standard
structure of $\C^2$ (see Lemma~\ref{PP} for the proof). Of course,
$\Gamma'$ admits a similar local representation. In what follows we
assume that we are in this setting. 

Let $(p^k)$ be a sequence of points in $D$  converging to $0$ and let
$\Sigma := \{ z \in \C^2: 2Re z^2 + 2Re K(z^1,0) + H(z^1,0) < 0\}$,
$\Sigma' := \{ z \in \C^2: 2Re z^2 + 2Re K'(z^1,0) + H'(z^1,0) < 0\}$.
The scaling procedure associates with the pair $(f,(p^k)_k)$
a biholomorphism $\phi$ (with respect to the standard structure $J_{st}$)
between $\Sigma$ and $\Sigma'$. Since $\phi$ is obtained as a limit of a
sequence of biholomorphic maps conjugated with $f$, some of their properties
are related and this can be used to study boundary properties of
$f$ and to prove that its cotangent lift is continuous up to the conormal
bundle $\Sigma(\partial D)$.

\subsubsection{Fixing suitable local coordinates and dilations.}
For any boundary point $t \in
\partial D$ we consider the change of variables $\alpha^t$ defined by 
$$
(z^1)^* = \frac{\partial \rho}{\partial \bar z^2}(t)(z^1 - t^1)
- \frac{\partial \rho}{\partial \bar z^1}(t)(z^2 - t^2),\
(z^2)^* = \sum_{j=1}^2 \frac{\partial \rho}{\partial z^j}(t)(z^j -
t^j).
$$
Then $\alpha^t$ maps $t$ to $0$. The real
normal at $0$ to $\Gamma$ is mapped by $\alpha^t$ to the line $\{
z^1 = 0, y_2 = 0 \}$. 
For every $k$, we denote by $t^k$ the projection of
$p^k$ onto $\partial D$ and by $\alpha^k$ the change of variables
$\alpha^t$ with $ t = t^k$.  Set $\delta_k = dist(p^k,\Gamma)$.
Then $\alpha^k(p^k) = (0,-\delta_k)$ and $\alpha^k(D)=\{2Re z^2 + O(\vert z
\vert^2) < 0\}$ near the origin. Since the sequence $(\alpha^k)_k$ converges
to the identity map, the sequence $(\alpha^k)_*(J)$ of almost
complex structures tends to $J$ as $k \rightarrow \infty$. Moreover
there is a sequence $(L^k)$ of linear automorphisms of $\R^4$
such that $(L^k \circ \alpha^k)_*(J)(0) = J_{st}$.
Then $(L^k \circ \alpha^k)(p^k) = (o(\delta_k),-\delta_k')$ with
$\delta_k' \sim \delta_k$ and
$(L^k \circ \alpha^k)(D) = \{Re (z^2 + \tau_k z^1) + O(\vert z \vert^2) < 0\}$
near the origin, with $\tau_k = o(1)$.
Hence there is sequence $(M^k)$ of
$\C$-linear transformations of $\C^2$, converging to the identity, 
such that $(T^k: = M^k \circ L^k
\circ \alpha^k)$ is a sequence of linear transformations converging to the
identity, and $D^k:= T^k(D)$ is defined near the origin by
$D^k=\{\rho_k(z) = Re z^2 + O(\vert z \vert^2) < 0\}$.
Finally $\tilde p_k = T^k(p^k)= (o(\delta_k),\delta_k'' + io(\delta_k))$ with
$\delta_k''\sim \delta_k$. 
We also denote by $\Gamma^k = \{\rho_k = 0 \}$ the
image of $\Gamma$ under $T^k$.
Furthermore, the sequence of almost complex structures
$(J_k:= (T^k)_*(J))$ converges to $J$ as $k \rightarrow \infty$
and $J_k(0) = J_{st}$.

We proceed quite similarly for the target domain $D'$. 
For $s \in \Gamma'$ we define the transformation $\beta^s$ by
$$
(z^1)^* = \frac{\partial \rho'}{\partial \bar z^2}(s)(z^1 - s^1)
- \frac{\partial \rho'}{\partial \bar z^1}(s)(z^2 - s^2),
(z^2)^* = \sum_{j=1}^2 \frac{\partial \rho'}{\partial z^j}(s)(z^j -
s^j).
$$

Let $s^k$ be the projection of $q^k:=f(p^k)$ onto $\Gamma'$ and
let $\beta^k$ be the corresponding map $\beta^s$ with $s = s^k$. 
The sequence $(q^k)$ converges
to  $0 = f(0)$ so $\beta^k$ tends to the identity. Considering linear
transformations $(L')^k$ and $(M')^k$, we obtain a sequence $(T'^k)$ of
linear transformations converging to the identity and satisfying the
following properties. The domain 
$(D^k)':= T'^k(D')$ is defined near the origin by
$(D^k)'=\{\rho_k'(z) := Re z^2 + O(\vert z \vert^2) < 0\}$,
$\Gamma_k' = \{ \rho_k' = 0 \}$ and $\tilde q_k = T'^k(q^k) =
(o(\varepsilon_k),\varepsilon_k''+ io(\varepsilon_k))$
with $\varepsilon_k'' \sim \varepsilon_k$, where
$\varepsilon_k = dist(q^k,\Gamma')$. 
The sequence of almost complex structures $(J_k':= (T'^k)_*(J'))$
converges to $J'$ as $k \rightarrow \infty$ and $J_k'(0) = J_{st}$.  

Finally, the map $ f^k:= T'^k \circ f \circ (T^k)^{-1}$ satisfies
$f^k(\tilde p_k) = \tilde q_k$  and is
a biholomorphism between the domains $D^k$ and $(D')^k$ 
with respect to the almost
complex structures $J_k$ and $J_k'$. 

Consider now the non isotropic dilations
$\phi_k : (z^1,z^2) \mapsto (\delta_k^{1/2}z^1,\delta_kz^2)$ and
$\psi_k(z^1,z^2)=
(\varepsilon_k^{1/2}z^1,\varepsilon_kz^2)$ and set $\hat f^k =
(\psi_k)^{-1} \circ f^k \circ \phi_k$.
Then the map $\hat f^k$ is biholomorphic with respect to the almost complex
structures $\hat J_k:=((\phi_k)^{-1})_*(J_k)$ and
$\hat J'_k:= (\psi_k^{-1})_*(J'_k)$.
Moreover if $\hat D^k:=\phi_k^{-1}(D^k)$ and
$(\hat{D'})^k:=\psi_k^{-1}((D')^k)$ then
$\hat D^k = \{ z \in \phi_k^{-1}(U): \hat \rho_k(z) < 0\}$
where
$$
\begin{array}{lll}
\hat \rho_k(z) &: =& \delta_k^{-1}\rho(\phi_k(z))\\
 & = & 2Re z^2 + \delta_k^{-1}[2
Re K(\delta_k^{1/2}z^1,\delta_kz^2) + H(\delta_k^{1/2}z^1,\delta_kz^2)
+  o(\vert (\delta_k^{1/2}z^1,\delta_kz^2) \vert^2).
\end{array}
$$
and $(\hat D')^k=\{ z \in \phi_k^{-1}(U): \hat \rho'_k(z) < 0\}$
where
$$
\begin{array}{lll}
\hat \rho'_k(z) &: =& \varepsilon_k^{-1}\rho'(\psi_k(z))\\
 & = & 2Re z^2 +
\varepsilon_k^{-1}[2 Re K'(\varepsilon_k^{1/2}z^1,\varepsilon_kz^2) +
H'(\varepsilon_k^{1/2}z^1,\varepsilon_kz^2)
+  o(\vert (\varepsilon_k^{1/2}z^1,\varepsilon_kz^2) \vert^2).
\end{array}
$$
Since $U$
is a fixed neighborhood of the origin, the pullbacks $\phi_k^{-1}(U)$
tend to $\C^2$ and the functions $\hat\rho_k$ tend
to $\hat \rho(z) = 2Re z^2 + 2Re K(z^1,0) + H(z^1,0)$ in the $\CC^2$ norm
on any compact subset of $\C^2$. Similarly, since $U'$
is a fixed neighborhood of the origin, the pullbacks $\psi_k^{-1}(U')$
tend to $\C^2$ and the functions $\hat\rho_k'$ tend
to $\hat \rho'(z) = 2Re z^2 + 2Re K'(z^1,0) + H'(z^1,0)$ in the $\CC^2$ norm
on any compact subset of $\C^2$. If $\Sigma :=
\{ z \in \C^2: \hat \rho(z) < 0 \}$ and $\Sigma' := \{ z \in \C^2:
\hat \rho'(z) < 0 \}$ then the sequence of points $\hat p^k =
\phi_k^{-1}(\tilde p_k) \in \hat D^k$ converges to the point $(0,-1) \in
\Sigma$ and the sequence of points $\hat q^k =
\psi^{-1}_k(\tilde q^k) \in \hat{D'}^k$ converges to $(0,-1) \in
\Sigma'$. Finally $\hat{f}^k(\hat p^k) = \hat q^k$.

\subsubsection{Convergence of the dilated families.} We begin with the
following
 
\begin{lemma}\label{convseq2}
The sequences $(\hat J'_k)$ and $(\hat J_k)$ of almost complex structures
converge to the standard structure uniformly (with all partial
derivatives of any order) on compact subsets of $\C^2$.
\end{lemma}

\noindent{\it Proof of Lemma~\ref{convseq2}.}
Denote by $a_{\nu\mu}^k(z)$ the elements of the matrix
$J_k$. Since $J_k \rightarrow J$ and $J$ is diagonal, we have $a_{\nu\mu}^k
\rightarrow a_{\nu\mu}$ for $\nu = \mu$ and $a_{\nu\mu}^k
\rightarrow 0$ for $\nu \neq \mu$. Moreover, since $J_k(0) =
J_{st}$, $a_{\nu\mu}^k(0) = i$ for $\nu = \mu$ and $a_{\nu\mu}(0) = 0$
for $\nu \neq \mu$.
 The elements $\hat
a_{\nu\mu}^k$ of the matrix
$\hat J_k$ are given by: $\hat a_{\nu\mu}^k(z^1,z^2) = a_{\nu
\mu}^k(\delta_k^{1/2}z^1,\delta_k z^2)$ for $\nu = \mu$, $\hat
a_{12}^k(z^1,z^2) = \delta_k^{1/2}a(\delta_k^{1/2}z^1,\delta_k z^2)$ and
$\hat a_{21}^k(z^1,z^2) = \delta_k^{-1/2}a_{21}^k(\delta_k^{1/2}z^1,\delta_k
z^2)$. This implies the desired result. \qed

\vskip 0,1cm
The next statement is crucial.
\begin{proposition}
\label{scaling}
The sequence $(\hat f^k )$ (together with all derivatives) is a
relatively compact family (with
respect to the compact open topology) on
$\Sigma$; every cluster point $\hat f$ is
a biholomorphism (with respect to $J_{st}$) between $\Sigma$ 
and $\Sigma'$, satisfying $\hat f(0,-1) = (0,-1)$ and
$(\partial \hat f^2/\partial z^2)(0,-1) = 1$. 
\end{proposition}
\noindent{\it Proof of Proposition~\ref{scaling}.}
{\it Step 1: convergence.} Our proof is based on the method
developed by F.Berteloot-G.Coeur{\'e}~\cite{BerCo}.
Consider a domain
$G \subset \C^2$ of the form $G = \{ z \in W: \lambda(z) = 2Re z^2 +
2Re K(z) + H(z) + o(\vert z \vert^2) < 0 \}$ where $W$ is a
neighborhood of the origin. We assume that an almost complex
structure $J$ is diagonal on $W$ and that the hypersurface 
$\{ \lambda = 0 \}$ is strictly $J$-pseudoconvex at any point.
Given $a \in \C^2$ and $\delta > 0$ 
denote by $Q(a,\delta)$ the non-isotropic ball
$Q(a,\delta ) = \{ z: \vert z^1 - a_1 \vert < \delta^{1/2}, \vert z^2
- a_2 \vert < \delta \}$. Denote also by $d_{\delta}$ the non-isotropic
dilation $d_{\delta}(z^1,z^2) = (\delta^{-1/2}z^1,\delta^{-1}z^2)$. 
\begin{lemma}\label{Control}
There exist positive constants $\delta_0, C, r$ satisfying the
following property : for
every $\delta \leq \delta_0$ and for every $J$-holomorphic disc $g:\Delta
\rightarrow G$ such that $g(0) \in Q(0,\delta)$ we have the
inclusion $g(r\Delta) \subset Q(0,C\delta)$.
\end{lemma}
\noindent{\it Proof of Lemma~\ref{Control}.}
Assume by contradiction that there exist positive sequences $\delta_k
\rightarrow 0$, $C_k \rightarrow +\infty$, a sequence $\zeta_k \in
\Delta$, $\zeta_k \rightarrow 0$ and a sequence $g_k: \Delta
\rightarrow G$ of
$J$-holomorphic discs such that $g_k(0) \in Q(0,\delta_k)$ and
$g_k(\zeta_k) \not\in Q(0,C_k\delta_k)$. Denote by $d_k$ the
dilations $d_{\delta}$ with $\delta = \delta_k$ and consider the
composition $h_k = d_k \circ g_k$ defined on $\Delta$. The
dilated domains $G_k:= d_k(G)$ are defined by $\{ z \in d_k(W):
\lambda_k(z):= \delta_k^{-1}\lambda \circ d_k^{-1}(z) < 0 \}$ and the
sequence $(\lambda_k)$ converges uniformly on compact subsets of
$\C^2$ to $\hat \lambda : z \mapsto 2Re z^2 + 2Re K(z) + H(z^1,0)$. Since $J$
is diagonal, the sequence  of structures $J_k:=(d_k)_*(J)$ converges to
$J_{st}$ in the $\CC^2$ norm on compact subsets of $\C^2$.  

The discs $h_k$ are $J_k$-holomorphic and the sequence $(h_k(0))$ is
contained in $Q(0,1)$; passing to a subsequence we may assume that this
converges to a point $p \in \overline{Q(0,1)}$. On the other hand, the
function $\hat \lambda + A\hat \lambda^2$ is
strictly $J_{st}$-plurisubharmonic on $Q(0,5)$ for a suitable constant
$A > 0$. Since the structures $J_k$ tend
to $J_{st}$, the functions $\lambda_k + A\lambda_k^2$
are strictly $J_k$-plurisubharmonic on $Q(0,4)$ for every $k$ large
enough and their Levi forms admit a uniform
lower bound with respect to $k$.
By Proposition~\ref{addest} the Kobayashi-Royden infinitesimal
pseudometric on $G_k$ admits the following lower bound~:
$K_{G_k}(z,v) \geq C \vert v \vert$ for any $z \in G_k \cap Q(0,3)$,
$v \in \C^2$,
with a positive constant $C$ independent of $k$. Therefore, there exists a
constant $C' > 0$ such that $\vert \vert \vert (dh_k)_\zeta \vert \vert \vert
\leq C'$ for any $\zeta \in (1/2)\Delta$ satisfying
$h_k(\zeta) \in G_k \cap Q(0,3)$.
On the other hand, the sequence $(\vert h_k(\zeta_k) \vert)$ tends to $+
\infty$. Denote by $[0,\zeta_k]$ the segment 
(in $\C$) joining the origin and $\zeta_k$ and let 
$\zeta_k' \in [0,\zeta_k]$ be the point the closest to the origin such
that $h_k([0,\zeta_k']) \subset G_k \cap
\overline{Q(0,2)}$ and $h_k(\zeta_k') \in \partial Q(0,2)$. Since $h_k(0)
\in Q(0,1)$, we have  $\vert h_k(0) - h_k(\zeta_k') \vert \geq C''$
for some constant $C'' > 0$. Let $\zeta_k' = r_k e^{i\theta_k}$, $r_k \in
]0,1[$. Then 
$$
\vert h_k(0) - h_k(\zeta_k') \vert \leq \int_{0}^{r_k}
\vert \vert \vert (dh_k)(te^{i\theta_k}) \vert \vert \vert dt \leq C'r_k
\rightarrow 0.
$$
This contradiction proves Lemma~\ref{Control}. \qed

\vskip 0,1cm
The statement of Lemma~\ref{Control} remains true if we replace the unit
disc $\Delta$ by the unit ball $\B_2$ in $\C^2$ equipped with an almost
complex structure $\tilde J$ close enough (in the $\CC^2$ norm) to
$J_{st}$. For the proof it is sufficient to foliate $\B_2$ by $\tilde
J$-holomorphic curves through the origin (in view of a smooth
dependence on  small perturbations of $J_{st}$ 
such a foliation is a small
perturbation of the foliation by complex lines through the origin and
apply Lemma~\ref{Control} to the foliation. 

\vskip 0,1cm
As a corollary we have the following
\begin{lemma}
\label{conv}
Let $(M,\tilde J)$ be an almost complex manifold and let $F^k: M
\rightarrow G$ be a sequence of $(\tilde J,J)$-holomorphic maps.
Assume that for some point $p^0 \in M$ we have $F^k(p) = (0,-\delta_k)$,
$\delta_k \rightarrow 0$, and that the sequence $(F^k)$ converges
to $0$ uniformly on compact subsets of $M$. 
Consider the rescaled maps $d_k \circ
F^k$. Then for any compact subset $K \subset M$ the sequence of norms
$(\| d_k \circ F^k \|_{\mathcal C^0(K)})$ is bounded.
\end{lemma}
\noindent{\sl Proof of Lemma~\ref{conv}}. It is sufficient to consider
a covering of a compact subset of $M$ by sufficiently small balls,
similarly to \cite{BerCo}, p.84.
Indeed, consider a covering of $K$ by the balls $p^j + r\B$, $j=0,...,N$
where $r$ is given by Lemma~\ref{Control} and $p^{j+ 1} \in p^j + r\B$ for
any $j$. For $k$ large enough, we obtain
that $F^k(p^0 + r\B) \subset Q(0,2C\delta_k)$, and $F^k(p^1 + r\B)
\subset Q(0,4C^2\delta_k)$. Continuing this process we obtain that
$F^k(p^N + r\mathbb B) \subset Q(0,2^NC^N\delta_k)$.
This proves Lemma~\ref{conv}. \qed

\vskip 0,1cm
Now we return to the proof of Proposition~\ref{scaling}. Lemma
\ref{conv} implies that the sequence $(\hat f^k)$ is bounded (in the $\CC^0$
norm) on any
compact subset $K$ of $\Sigma$. Covering $K$ by small bidiscs,
consider two transversal foliations by $J$-holomorphic curves on every
bidisc. Since the restriction of $\hat f^k$ on every such curve is
uniformly bounded in the $\CC^0$-norm, 
it follows by the elliptic
estimates that this is bounded in $\CC^l$ norm for every $l$ (see
\cite{si}). Since the bounds are uniform with respect to curves,
the sequence $(\hat f^k)$ is bounded in every
$\CC^l$-norm. So the family $(\hat f^k)$ is relatively compact. 

{\it Step 2: Holomorphicity of the limit maps.} Let $(\hat f^{k_s})$ be a
subsequence converging to a smooth map $\hat f$. 
Since $f^{k_s}$ satisfies the holomorphicity condition
$\hat J'_{k_s} \circ d\hat f^{k_s} = d\hat f^{k_s} \circ
J_{k_s}$, since $\hat J_{k_s}$ and $\hat J'_{k_s}$ converge to $J_{st}$,
we obtain, passing to the limit in the holomorphicity condition, that $\hat f$
is holomorphic with respect to $J_{st}$.

{\it Step 3: Biholomorphicity of $\hat f$.}
Since $\hat f(0,-1) = (0,-1) \in \Sigma'$ and
$\Sigma'$ is defined by a plurisubharmonic function, it follows by the
maximum principle that $\hat f(\Sigma) \subset \Sigma'$ (and not just a
subset of $\bar\Sigma'$). Applying a similar argument to the
sequence $(\hat f^k)^{-1}$ of inverse map, we obtain that this converges
(after extraction of a subsequence) to the inverse of $\hat f$. 

Finally the domain $\Sigma$ (resp. $\Sigma'$) is
biholomorphic to ${\mathbb H}$ by means of the transformation $(z^1,z^2)
\mapsto (z^1,z^2 + K(z^1,0))$ (resp. $(z^1,z^2) \mapsto (z^1,z^2 +
K'(z^1,0))$). Since a biholomorphism of ${\mathbb H}$ fixing the point
$(0,-1)$ has the form
$(e^{i\theta}z^1,z^2)$ (see, for instance, \cite{co}), $\hat f$ is conjugated
to this transformation by the above quadratic biholomorphisms of
$\C^2$. Hence~: 

\begin{eqnarray}
\label{derivative}
\frac{\partial \hat f^2}{\partial z^2}(0,-1) = 1.
\end{eqnarray}

\vskip 0,1cm
\noindent This property will be used in the next Section. \qed

\subsubsection{Boundary behavior of the tangent map} 

We suppose that we are in the local situation described at the
beginning of the previous section. Here we prove two statements
concerning the boundary behavior of the tangent map of $f$ near
$\Gamma$. They are obvious if $f$ is of class $\CC^{1}$ up to
$\Gamma$. In the general situation, their proofs
require the scaling method of the previous section.
Let $p \in \Gamma$. After a local change of coordinates $z$ we
may assume that $p = 0$, $J(0) = J_{st}$ and $J$ is assumed to be diagonal.
In the $z$ coordinates, we 
consider a base $X$ of (1,0) (with respect to $J$) vector fields
defined previously. Recall that $X_2 = \partial
/\partial z^2 + a(z) \partial/\bar\partial z^2$, $a(0) = 0$,
$X_1(0) = \partial/\partial z^1$ and at every point $z^0$, $X_1(z^0)$
generates the holomorphic tangent space $H_z^J(\partial D - t)$, $t \geq 0$.
If we return to the initial coordinates and move the point $p \in \Gamma$,
we obtain for every $p$ a basis $X_p$ of $(1,0)$ vector fields, defined in a
neighborhood of $p$. Similarly, we define the basis $X'_q$ for
$q \in \partial D'$. 

The elements of the matrix of the tangent map
$df_z$ in the bases $X_p(z)$ and $X'_{f(p)}(z)$ are denoted by
$A_{js}(p,z)$. According to Proposition~\ref{matrix} the function
$A_{22}(p,\cdot)$ is upper bounded on $D$. 

\begin{proposition}
\label{REALITY}
We have:
\begin{itemize}
\item[(a)] Every cluster point of the function $z \mapsto A_{22}(p,z)$
(in the notation of Proposition~\ref{matrix}) is real when $z$ tends to a
point $p \in \partial D$.
\item[(b)] For $z \in D$, let $p \in \Gamma$ such that
$|z-p| = dist(z,\Gamma)$. There exists a constant $A$, independent of
$z \in D$, such that $\vert A_{22}(p,z) \vert \geq A$.
\end{itemize}
\end{proposition}

The proof of these statements use the above scaling
construction. So we use the notations of the previous section. 

\vskip 0,1cm
\noindent{\it Proof of Proposition~\ref{REALITY}}.
(a) Suppose that there exists a sequence of points $(p^k)$ converging
to a boundary point $p$ such that $A_{22}(p,\cdot)$ tends to a complex number
$a$. Applying the above scaling construction,
we obtain a sequence of maps $(\hat f^k)_k$.
Consider the two basis $\hat X^k:= 
\delta_k^{1/2}((\phi_k^{-1}) \circ T^k)(X_1),
\delta_k((\phi_k^{-1})\circ T^k)(X_2))$ and $(\hat
X')^k:= (\varepsilon_k^{-1/2}((\psi_k^{-1}) \circ T'^k)(X'_1),
\varepsilon_k^{-1}((\psi_k^{-1})\circ T'^k)(X'_2))$. These vector
fields tend to the standard (1,0) vector field base of $\C^2$ as $k$
tends to $\infty$. Denote by $\hat A^k_{js}$ the elements of the
matrix of $d\hat f^k(0,-1)$. Then $A^k_{22} \rightarrow (\partial
\hat f^2/\partial z^2)(0,-1) = 1$, according to (\ref{derivative}). On the
other hand, $A^k_{22} = \varepsilon_k^{-1}\delta_k A_{22}$ and tends
to $a$ by the boundary distance preserving property (Proposition~\ref{equiv}).
This gives the statement.

(b) Suppose that there is a sequence of points $(p^k)$ converging
to the boundary such that $A_{22}$ tends to $0$. Repeating precisely
the argument of (a), we obtain that $(\partial \hat f^2/\partial
z^2)(0,-1) = 0$; this contradicts~(\ref{derivative}). \qed

\vskip 0,1cm
In order to establish the next proposition, it is convenient to associate
a wedge with the totally real part of the conormal bundle
$\Sigma_J(\partial D)$ of $\partial D$ as edge.
Consider in $\R^{4} \times \R^{4}$ the set $S = \{ (z,L):
dist((z,L),\Sigma_J(\partial D)) \leq dist(z,\partial D), z \in D \}$.
Then, in a neighborhood $U$ of any totally real point of
$\Sigma_J(\partial D)$, the set S contains a wedge $W_U$ with
$\Sigma_J(\partial D) \cap U$ as totally real edge.

\begin{proposition}
\label{cluster2}
Let $K$ be a compact subset of the totally real part of the conormal
bundle $\Sigma_J(\partial D)$. Then the cluster set of the cotangent lift
$\tilde f $ of $f$ on the conormal bundle 
$\Sigma(\partial D)$, when $(z,L)$ tends to $\Sigma_J(\partial D)$
along the wedge $W_U$, is relatively compactly contained 
in the totally real part of $\Sigma(\partial D')$.
\end{proposition}
\noindent{\sl Proof of Proposition~\ref{cluster2}}.
Let $(z^k,L^k)$ be a sequence in $W_U$
converging to $(0,\partial_J\rho(0)) = (0,dz^2)$. Set $g = f^{-1}$.
We shall prove that the
sequence of linear forms $Q^k :=
{}^tdg(w^k)L^k$, where $w^k = f(z^k)$, converges to a linear form which up to
a {\it real} factor (in view of Part (a) of
Proposition \ref{REALITY})
coincides with $\partial_J \rho'(0)= dz^2$
(we recall that ${}^t$ denotes the transposed map).
It is sufficient to prove that the first component of $Q^k$ with
respect to the dual basis $(\omega_1,\omega_2)$ of $X$ tends to $0$
and 
the second one is
bounded below from the origin as $k$ tend to infinity. The map $X$ being
of class $\CC^1$ we can replace $X(0)$ by $X(w^k)$.
Since $(z^k,L^k) \in W_U$, we have $L^k
= \omega_2(z^k) + O(\delta_k)$, where $\delta_k$ is the distance from
$z^k$ to 
the boundary. Since $\vert\vert\vert dg_{w^k} \vert\vert\vert =
0(\delta_k^{-1/2})$, we
have $Q^k = {}^tdg_{w^k}(\omega_2(z^k)) + O(\delta_k^{1/2})$.
By Proposition~\ref{matrix}, the
components of ${}^tdg_{w^k}(\omega_2(z^k))$ with respect to the basis
$(\omega_1(z^k),\omega_2(z^k))$ are the elements of the 
second line of the matrix 
$dg_{w^k}$ with respect to the basis $X'(w^k)$ and $X(z^k)$. So its first 
component is $0(\delta_k^{1/2})$ and tends to $0$ as $k$ tends to
infinity. Finally the component $A_{22}^k$ is bounded below from the
origin by Part (b) of Proposition~\ref{REALITY}. \qed

%\vskip 0,1cm
%\noindent{\it Proof of Theorem \ref{MainTheorem}}. In view of
%Proposition~\ref{cluster}, we may apply Proposition~\ref{wedges} to the
%cotangent lift $\tilde f$ of $f$.
%This gives the statement of Theorem~\ref{MainTheorem}. \qed

\subsection{Lifts of biholomorphisms to the cotangent bundle~: the
  case of arbitrary dimension}

Our further considerations rely deeply on the estimates
of the Kobayashi-Royden infinitesimal pseudometric given in Section 3.
\vskip 0,1cm
Let $D$ (resp. $D'$) be a strictly pseudoconvex domain in an almost complex
manifold $(M,J)$ (resp. $(M',J')$) and let $f$ be a $(J,J')$-biholomorphism
from $D$ to $D'$. Fix a point $p \in \partial D$ and a sequence $(p^k)_k$
in $D$ converging to $p$. After extraction we may assume that the sequence
$(f(p^k))_k$ converges to a point $p'$ in $\partial D'$. According to the
Hopf lemma, $f$ has the boundary distance property. Namely, there is a
positive constant $C$ such that
\begin{equation}\label{bdp}
(1/A) \ dist(f(p^k), \partial D') \leq dist(p^k, \partial D) \leq
A \ dist(f(p^k), \partial D'),
\end{equation}
where $A$ is independent of $k$ (see~\cite{co-ga-su}).

Since all our considerations are local we set $p=p'=0 \in \C^n$.
We may assume that $J(0) = J_{st}$ and $J'(0) = J_{st}$.
Let $U$  (resp. $V$) be a neighborhood of the origin in $\C^n$ such that
$D \cap U = \{z \in U : \rho(z,\bar{z}) : =
z_n + \bar{z}_n + Re(K(z)) + H(z,\bar{z}) + \cdots < 0\}$ 
(resp. $D' \cap V = \{w \in V : \rho'(w,\bar{w}) : =
w_n + \bar{w}_n + Re(K'(w)) + H'(w,\bar{w}) + \cdots < 0\}$) 
where 
$K(z)  = \sum k_{\nu\mu} z^{\nu}{z}^{\mu}$, $k_{\nu\mu} =
k_{\mu\nu}$,
$H(z) = \sum h_{\nu\mu} z^{\nu}\bar z^{\mu}$, $h_{\nu\mu} =
\bar h_{\mu\nu}$ and $\rho$ is a strictly $J$-plurisubharmonic function on
$U$ (resp. $K'(z)  = \sum k'_{\nu\mu} w^{\nu}{w}^{\mu}$, $k'_{\nu\mu} =
k'_{\mu\nu}$,
$H'(w) = \sum h'_{\nu\mu} w^{\nu}\bar w^{\mu}$, $h'_{\nu\mu} =
\bar h'_{\mu\nu}$ and $\rho'$ is a strictly $J'$-plurisubharmonic function on
$V$).

\subsubsection{Asymptotic behaviour of the tangent map of $f$}
We wish to understand the limit behaviour (when $k \rightarrow \infty$) of
$df(p^k)$. Consider the vector fields
$$
v^j:=(\partial \rho/\partial x^n)\partial / \partial x^j
- (\partial \rho/\partial x^j)\partial / \partial x^n
$$
for $j=1,\dots,n-1$, and
$$
v^n:=(\partial \rho/\partial x^n)\partial /
\partial y^n - (\partial \rho/\partial y^n)\partial / \partial x^n.
$$
Restricting $U$ if necessary, the vector fields $X^1,\dots,X^{n-1}$
defined by $X^j:=v^j-iJ(v^j)$ form a basis 
of the $J$-holomorphic tangent space to $\{\rho = \rho(z)\}$
at any $z \in U$. Moreover, if $X^n:=v^n-iJv^n$ then the family
$X:=(X^1,\dots,X^n)$ forms a basis of $(1,0)$ vector fields on $U$.
Similarly we define a basis
$X':=(X'^1,\dots,X'^n)$ of $(1,0)$ vector fields on $V$ such that
$(X'^1(w),\dots,X'^{n-1}(w))$ defines a basis
of the $J'$-holomorphic tangent space to $\{\rho' = \rho'(w)\}$
at any $w \in V$. 
We denote by $A(p^k):=(A(p^k)_{j,l})_{1 \leq j,l \leq n}$ the matrix of the
map $df(p^k)$ in the basis $X(p^k)$ and $X(f(p^k))$.

\begin{remark}\label{precise}
In sake of completeness we should write $X_0$ and $X'_0$ to emphasize that
the structure was normalized by the condition $J(0) = J_{st}$ and
$A(0,p^k)$ for $A(p^k)$. The same construction is valid for any
boundary point of $D$.
The corresponding notations will be used in Proposition~\ref{reality}.
\end{remark}
\begin{proposition}\label{tangent}
The matrix $A(p^k)$ satisfies the following estimates~:
$$
A(p^k)=\left(
\begin{array}{ccc}
O_{n-1,n-1}(1) & & O_{n-1,1}(dist(p^k,\partial D)^{-1/2})\\
 & & \\
O_{1,n-1}(dist(p^k,\partial D)^{1/2}) & & O_{1,1}(1)
\end{array}
\right).
$$
\end{proposition}

The matrix notation means that the following estimates are satisfied~:
$A(p^k)_{j,l} = O(1)$ for $1 \leq j,l \leq n-1$,
$A(p^k)_{j,n} = O(dist(p^k,\partial D)^{-1/2})$ for $1 \leq j \leq n-1$,
$A(p^k)_{n,l} = O(dist(p^k,\partial D)^{1/2})$ for $1 \leq l \leq n-1$ and
$A(p^k)_{n,n} = O(1)$. 
\vskip 0,1cm

The proof of Proposition~\ref{tangent} was given previously
in dimension two (Proposition~\ref{matrix}) but is valid without any
modification in any dimension. We note
that the asymptotic behaviour of $A(p^k)$ depends only on the distance
from the point to $\partial D$, not on the choice of the sequence
$(p^k)_k$.

\subsubsection{Scaling process and model domains}

The following construction is similar to the two dimensional case.
For every $k$ denote by $q^k$ the projection of $p^k$ to $\partial D$ and
consider the change of variables $\alpha^k$ defined by 
$$
\left\{
\begin{array}{ccccc}
(z^j)^* & = & \displaystyle \frac{\partial \rho}{\partial \bar z^n}(q^k)
(z^j - (q^k)^j)
- \displaystyle \frac{\partial \rho}{\partial \bar z^j}(q^k)(z^n - (q^k)^n),
& &
{\rm for} \ 1 \leq j \leq n-1,\\
(z^n)^* & = & \sum_{j=1}^n \displaystyle \frac{\partial \rho}{\partial z^j}
(q^k)(z^j - (q^k)^j).
\end{array}
\right.
$$
If $\delta_k := dist(p^k,\partial D)$ then $\alpha^k(p^k) = (0,-\delta_k)$
and $\alpha^k(D)=\{2Re z^n + O(\vert z
\vert^2) < 0\}$ near the origin. Moreover, the sequence $(\alpha^k)_*(J)$
converges to $J$ as $k \rightarrow \infty$, since the sequence
$(\alpha^k)_k$ converges
to the identity map. Let $(L^k)_k$ be a sequence of linear automorphisms of
$\R^{2n}$
such that $(T^k: = L^k
\circ \alpha^k)_k$ converges to the
identity, and $D^k:= T^k(D)$ is defined near the origin by
$D^k=\{\rho_k(z) = Re z^n + O(\vert z \vert^2) < 0\}$.
The sequence of almost complex structures
$(J_k:= (T^k)_*(J))_k$ converges to $J$ as $k \rightarrow \infty$
and $J_k(0) = J_{st}$.
Furthermore $\tilde p_k := T^k(p^k)$ satisfies
$\tilde p_k = (o(\delta_k),\delta_k'' + io(\delta_k))$
with
$\delta_k''\sim \delta_k$. 

We proceed similarly on $D'$. 
Denote by $s^k$ the projection of $f(p^k)$ onto $\partial D'$ and
define the transformation $\beta^k$ by
$$
\left\{
\begin{array}{ccccc}
(w^j)^* & = & \displaystyle \frac{\partial \rho'}{\partial \bar w^n}(s^k)
(w^j - (s^k)^j)
- \displaystyle \frac{\partial \rho'}{\partial \bar w^j}(s^k)(w^n - (s^k)^n),
& & {\rm for} \ 1 \leq j \leq n-1,\\
(w^n)^* & = & \sum_{j=1}^n \displaystyle \frac{\partial \rho'}{\partial w^j}
(s^k)(w^j - (s^k)^j).
\end{array}
\right.
$$
We define a sequence $(T'^k)_k$ of
linear transformations converging to the identity and satisfying the
following properties. The domain 
$(D^k)':= T'^k(D')$ is defined near the origin by
$(D^k)'=\{\rho_k'(w) := Re w^n + O(\vert w \vert^2) < 0\}$,
and $\tilde f(p_k) = T'^k(f(p_k)) =
(o(\varepsilon_k),\varepsilon_k''+ io(\varepsilon_k))$
with $\varepsilon_k'' \sim \varepsilon_k$, where
$\varepsilon_k = dist(f(p_k),\partial D')$. 
The sequence of almost complex structures $(J_k':= (T'^k)_*(J'))_k$
converges to $J'$ as $k \rightarrow \infty$ and $J_k'(0) = J_{st}$.  

Finally, the map $ f^k:= T'^k \circ f \circ (T^k)^{-1}$ satisfies
$f^k(\tilde p_k) = \tilde f(p_k)$  and is
a $(J_k,J'_k)$-biholomorphism between the domains $D^k$ and $(D')^k$. 

Let $\phi_k : ('z,z^n) \mapsto (\delta_k^{1/2} {'z},\delta_kz^n)$ and
$\psi_k('w,w^n)=
(\varepsilon_k^{1/2} \ 'w,\varepsilon_kw^n)$ and set $\hat f^k =
(\psi_k)^{-1} \circ f^k \circ \phi_k$.
The map $\hat f^k$ is $(\hat J_k,\hat J'_k)$-biholomorphic, where
$\hat J_k:=((\phi_k)^{-1})_*(J_k)$ and
$\hat J'_k:= (\psi_k^{-1})_*(J'_k)$.
If $\hat D^k:=\phi_k^{-1}(D^k)$ and
$(\hat{D'})^k:=\psi_k^{-1}((D')^k)$ then
$\hat D^k = \{ z \in \phi_k^{-1}(U): \hat \rho_k(z) < 0\}$
where
$$
\begin{array}{lll}
\hat \rho_k(z) & : = & \delta_k^{-1}\rho(\phi_k(z))\\
 & = & 2Re z^n + \delta_k^{-1}[2
Re K(\delta_k^{1/2}{'z},\delta_kz^n) + H(\delta_k^{1/2}{'z},\delta_kz^n)
+  o(\vert (\delta_k^{1/2}{'z},\delta_kz^n) \vert^2).
\end{array}
$$
and $(\hat D')^k=\{w \in \phi_k^{-1}(V): \hat \rho'_k(z) < 0\}$
where
$$
\begin{array}{lll}
\hat \rho'_k(w) & : = &\varepsilon_k^{-1}\rho'(\psi_k(w))\\
 & = & 2Re w^n +
\varepsilon_k^{-1}[2 Re K'(\varepsilon_k^{1/2}\ 'w,\varepsilon_kw^n) +
H'(\varepsilon_k^{1/2}\ 'w,\varepsilon_kw^n)
+  o(\vert (\varepsilon_k^{1/2}\ 'w,\varepsilon_kw^n) \vert^2).
\end{array}
$$
Since $U$
is a neighborhood of the origin, the pullbacks $\phi_k^{-1}(U)$
converge to $\C^n$ and the functions $\hat\rho_k$ converge
to $\hat \rho(z) = 2Re z^n + 2Re K({'z},0) + H({'z},0)$ in the $\CC^2$ norm
on compact subsets of $\C^n$. Similarly, since $V$
is a neighborhood of the origin, the pullbacks $\psi_k^{-1}(U')$
converge to $\C^n$ and the functions $\hat\rho_k'$ converge
to $\hat \rho'(w) = 2Re w^n + 2Re K'({'w},0) + H'({'w},0)$ in the $\CC^2$ norm
on compact subsets of $\C^n$. If $\Sigma :=
\{z \in \C^n: \hat \rho(z) < 0 \}$ and $\Sigma' := \{w \in \C^n:
\hat \rho'(w) < 0 \}$ the sequence of points $\hat p_k =
\phi_k^{-1}(\tilde p_k) \in \hat D^k$ converges to the point $(0,-1) \in
\Sigma$ and the sequence of points $\hat f(p_k) =
\psi^{-1}_k(\tilde f(p_k)) \in \hat{D'}^k$ converges to $(0,-1) \in
\Sigma'$. Finally $\hat{f}^k(\hat p_k) = \hat f(p_k)$.

\noindent The limit behaviour of the dilated objects is given by the following
proposition (see Figure~6).

\begin{proposition}\label{convseq}
$(i)$ The sequences $(\hat J_k)$ and $(\hat J'_k)$ of almost complex
structures converge to model structures $J_0$ and $J'_0$
uniformly (with all partial derivatives of any order) on compact subsets of
$\C^n$.

\vskip 0,1cm
$(ii)$ $(\Sigma,J_0)$ and $(\Sigma',J'_0)$ are model domains.

\vskip 0,1cm
$(iii)$ The sequence $(\hat f^k)$ (together with all derivatives) is a
relatively compact family (with respect to the compact open topology) on
$\Sigma$; every cluster point $\hat f$ is
a $(J_0,J'_0)$-biholomorphism between $\Sigma$ 
and $\Sigma'$, satisfying $\hat f(0,-1) = (0,-1)$ and
$\hat f^n('0,z^n) = z^n$ on $\Sigma$.
\end{proposition}

\medskip
\begin{center}
\input{figure3.pstex_t}
\end{center}
\medskip
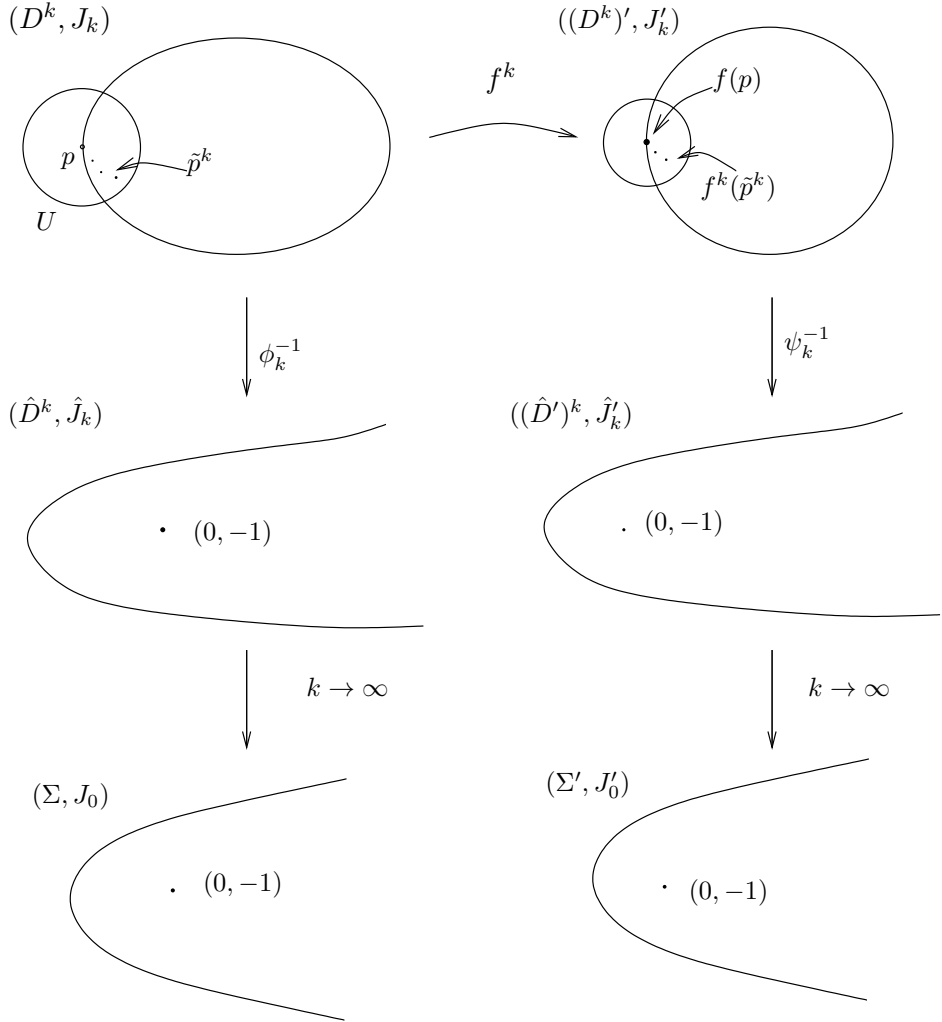
\centerline{Figure 6}
\bigskip

\noindent{\it Proof of Proposition~\ref{convseq}.}
We start with the proof of $(i)$. We focus on structures $\hat{J}_k$.
Consider $J=J_{st} + L(z) +
  O(|z|^2)$ as a matrix
valued function, where $L$ is a real linear matrix.
The Taylor expansion of $J_k$ at
the origin is given by  $J_k = J_{st} + L^k(z) +   O(|z|^2)$
on $U$, uniformly with respect to $k$. Here $L^k$ is a real linear
matrix converging to $L$ at infinity. 
Write $\hat{J}_k = J_{st} + \hat{L}^k +   O(\delta_k)$.
If $L^k=(L^k_{j,l})_{j,l}$ then
$\hat{L}^k_{j,l}= L^k_{j,l}(\phi_k(z))$ for $1 \leq j \leq n-1,\ 1 \leq
l \leq n$, $\hat{L}^k_{n,l}=\delta_k^{-1/2}L^k_{n,l}(\phi_k(z))$
for $1 \leq l \leq n-1$ and $\hat{L}^k_{n,n}=L^k_{n,n}(\phi_k(z))$.
This gives the conclusion. 

\vskip 0,1cm
Proof of $(ii)$. We focus on $(\Sigma,J_0)$.
By the invariance of the Levi form we
have ${\mathcal L}^{J_k}(\rho_k)(0)(\phi_k(v)) = {\mathcal
L}^{\hat J_k}(\rho_k \circ \phi_k)(0)(v)$.
Write $J_0 = J_{st} + L^\infty$.
Since $\rho_k$ is strictly $J_k$-plurisubharmonic uniformly with respect to
$k$ ($\rho_k$ converges to $\rho$ and $J_k$ converges to $J$),
multiplying by $\delta_k^{-1}$ and
passing to the limit at the right side as $k \rightarrow \infty$, 
we obtain that
${\mathcal L}^{J_0}(\hat \rho)(0)(v) \geq 0$ for any $v$. Now let $v =
(v',0) \in T_0(\partial \Sigma)$. Then
$\phi_k(v) = \delta_k^{1/2}v$ and so 
${\mathcal L}^J_k(\rho)(0)(v) = {\mathcal
  L}^{\hat J_k}(\rho_k)(0)(v)$. Passing to the limit as $k$
tends to infinity, we obtain that
${\mathcal L}^{J_0}(\hat \rho)(0)(v) > 0$
for any $v = (v',0)$ with $v' \neq 0$. 

\vskip 0,1cm
Proof of $(iii)$. The proof of the existence
and of the biholomorphicity of $\hat{f}$ is the same as in dimension two.
We prove the identity on $\hat f^n$.
Let $t$ be a real positive number. Then we have~:
\begin{lemma}\label{infty}
$\lim_{t \rightarrow \infty} \hat{\rho}'(\hat{f}('0,-t)) =  \infty$.
\end{lemma}
\noindent{\it Proof of Lemma~\ref{infty}}.
According to the boundary distance property~(\ref{bdp}) we have
$$
|\rho'(f \circ (T^k)^{-1} \circ \phi_k)('0,-t)| \geq C \ 
dist(T_k^{-1}('0,-\delta_k t)).
$$
Then
$$
|\hat{\rho}'_k(\hat{f}^k('0,-t))| \geq C \varepsilon_k^{-1}\delta_k \ t.
$$
Since $\hat{\rho}'_k$ converges to $\hat{\rho}'$ uniformly on
compact subsets of $\Sigma'$ and $\varepsilon_k \simeq \delta_k$ (by
the boundary distance property~(\ref{bdp})) we obtain~:
$$
|\hat{\rho}'(\hat{f}('0,-t))| \geq Ct.
$$
This proves Lemma~\ref{infty}. \qed

\vskip 0,1cm
We turn back to the proof of part $(iii)$ of Proposition~\ref{convseq}.
Assume first that $J$ (and similarly $J'$) are not integrable
(see Proposition~\ref{prop-hyp}). Consider a $J$-complex hypersurface
$A \times \C$ in $\C^n$ where $A$ is a $J_{st}$ complex hypersurface in
$\C^{n-1}$.
Since $f((A \times \C) \cap \mathbb H_{P_1}) = (A' \times \C) \cap
\mathbb H_{P_2}$ where $A'$ is a $J_{st}$ complex hypersurface in
$\C^{n-1}$, it follows that the restriction of $\hat f^n$ to $\{'z=\ '0,
Re(z^n) < 0\}$ is a $J_{st}$ automorphism of $\{'z={'0}, Re(z^n) < 0\}$.
Let $\phi : \zeta \mapsto (\zeta -1)/(\zeta + 1)$. The function
$\hat{g}:=\phi^{-1} \circ \hat{f}^n \circ \phi$ is a $J_{st}$ automorphism
of the unit disc in $\C$. In view of Lemma~\ref{infty} this satisfies
$\hat{g}(0) = 0$ and $\hat{g}(1) = 1$. Hence $\hat{g} \equiv id$ and
$\hat{f}^n('0,z^n) = z^n$ on $\Sigma$.

Assume now that $J$ and $J'$ are integrable.
Let $F$ (resp.
$F'$) be the diffeomorphism from $\Sigma$ to $\mathbb H_{P}$
(resp. from $\Sigma$ to $\mathbb H_{P'}$) given in the proof of
Proposition~\ref{prop-hyp}. The diffeomorphism $g:=F' \circ f \circ
F^{-1}$ is a $J_{st}$-biholomorphism from $\mathbb H_{P}$ to $\mathbb
H_{P'}$ satisfying $g('0,-1) = ('0,-1)$. Since $(\Sigma,J)$ and
$(\Sigma', J')$ are model domains, the domains $\mathbb H_{P}$ and
$\mathbb H_{P'}$ are strictly $J_{st}$-pseudoconvex. In particular, since
$P$ and $P'$ are homogeneous of degree two, there are linear complex maps
$L,\ L'$ in $\C^{n-1}$ such that the map $G$ (resp. $G'$) defined by
$G('z,z_n)=(L('z),z_n)$ (resp. $G'('z,z_n)=(L'('z),z_n)$) is a
biholomorphism from $\mathbb H_{P}$ (resp. $\mathbb H_{P'}$) to
$\mathbb H$. The map $G' \circ g \circ G^{-1}$ is an automorphism of
$\mathbb H$ satisfying $G' \circ g \circ G^{-1}('0,-1) = ('0,-1)$.
Let $\Phi$ be the $J_{st}$ biholomorphism from $\mathbb H$ to
the unit ball $\mathbb B_n$ of $\C^n$ defined by
$\Phi('z,z^n) = (\sqrt{2}'z/1-z^n,(1+z^n)/(1-z^n))$.
Let $\hat{g} :=\Phi^{-1} \circ g
\circ \Phi$. In view of lemma~\ref{bdp} this satisfies
$\hat{g}(0) = 0$ and $\hat{g}('0,1)=('0,1)$. Hence $\hat{g}^n \equiv id$
and $\hat{f}^n('z,z^n) = z^n$ for every $z$ in $\Sigma$. \qed

\vskip 0,1cm
According to part $(ii)$ of Proposition~\ref{convseq}
and restricting $U$ if necessary, one may view
$D \cap U$ as a strictly $J_0$-pseudoconvex domain
in $\C^n$ and $J$ as a small deformation of ${J}_0$ in a
neighborhood of $\bar{D} \cap U$. The same holds for $D' \cap V$.

\vskip 0,1cm
For $p \in \partial D$ and $z \in D$ let $X_p(z)$ and $X'_{f(p)}(f(z))$
be the basis of $(1,0)$ vector fields defined above.
The elements of the matrix of $df_z$ in
the bases $X_p(z)$ and $X'_{f(p)}(f(z))$ are denoted by
$A_{js}(p,z)$. According to Proposition~\ref{tangent} the function
$A_{n,n}(p,\cdot)$ is upper bounded on $D$. 

\begin{proposition}
\label{reality}
We have:
\begin{itemize}
\item[(a)] Every cluster point of the function $z \mapsto A_{n,n}(p,z)$
is real when $z$ tends to $p \in \partial D$.
\item[(b)] For $z \in D$, let $p \in \partial D$ such that
$|z-p| = dist(z,\partial D)$. There exists a constant $A$, independent of
$z \in D$, such that $\vert A_{n,n}(p,z) \vert \geq A$.
\end{itemize}
\end{proposition}

\vskip 0,1cm
\noindent{\it Proof of Proposition~\ref{reality}}.
(a) Suppose that there exists a sequence of points $(p^k)$ converging
to a boundary point $p$ such that $A_{n,n}(p,\cdot)$ tends to a complex number
$a$. Applying the above scaling construction,
we obtain a sequence of maps $(\hat f^k)_k$.
For $k \geq 0$ consider the dilated vector fields
$$
Y^j_k:=\delta_k^{1/2}((\phi_k^{-1}) \circ T^k)(X^j(p^k))
$$
for $j=1,\dots,n-1$, and
$$
Y^n_k:=\delta_k((\phi_k^{-1})\circ T^k)(X_n(p^k)).
$$
Similarly we define
$$
Y'^j_k:=\varepsilon_k^{-1/2}((\psi_k^{-1}) \circ T'^k)
(X'^j(f(p^k)))
$$
for $j=1,\dots,n-1$, and
$$
Y'^n_k:=\varepsilon_k^{-1}((\psi_k^{-1})\circ T'^k)(X'_n(f(p^k))).
$$
For every $k$, the $n$-tuple
$Y^k:= (Y^1_k,\dots,Y^n_k)$ is a basis of $(1,0)$ vector fields for
the dilated structure $\hat{J}^k$. In view of Proposition~\ref{convseq}
the sequence $(Y^k)_k$
converges to a basis
of $(1,0)$ vector fields of $\C^n$ (with respect to $J_0$) as $k$
tends to $\infty$. Similarly, the $n$-tuple
$Y'^k := (Y'^1_k,\dots,Y'^n_k)$ is a basis of $(1,0)$ vector fields for
the dilated structure $\hat{J}'^k$ and $(Y'^k)_k$
converges to a basis of $(1,0)$ vector fields
of $\C^n$ (with respect to $J'_0$) as $k$ tends to $\infty$.
In particular the last components $Y^n_k$ and $Y'^n_k$
converge to the $(1,0)$ vector field $\partial / \partial z^n$.
Denote by $\hat A^k_{js}$ the elements of the
matrix of $d\hat f^k(0,-1)$. Then $A^k_{n,n}$ converges to $(\partial
\hat f^n/\partial z^n)(0,-1) = 1$, according to Proposition~\ref{convseq}.
On the other hand, $A^k_{n,n} = \varepsilon_k^{-1}\delta_k A_{n,n}$
converges to $a$ by the boundary distance preserving property~(\ref{bdp}).
This gives the statement.

(b) Suppose that there is a sequence of points $(p^k)$ converging
to the boundary such that $A_{n,n}$ tends to $0$. Repeating precisely
the argument of (a), we obtain that $(\partial \hat f^n/\partial
z^n)(0,-1) = 0$; this contradicts part $(iii)$ of Proposition~\ref{convseq}.
\qed

\vskip 0,1cm
\begin{proposition}\label{PPP}
The cluster set of the cotangent lift $f^*$ on
$\Sigma(\partial D)$ is contained in $\Sigma(\partial D')$.
\end{proposition}

\vskip 0,1cm

\vskip 0,1cm
\noindent{\it Proof of Proposition~\ref{PPP}}.

{\it Step one.}
We first reduce the problem to the following
local situation. Let $D$ and $D'$ be domains in $\C^n$, $\Gamma$ and
$\Gamma'$ be open $\CC^{\infty}$-smooth pieces of their boundaries,
containing the origin. We assume that an almost complex structure $J$
is defined and $\CC^{\infty}$-smooth in a neighborhood of the closure
$\bar D$, $J(0) = J_{st}$.
Similarly, we assume that $J'(0) = J_{st}$. The hypersurface
$\Gamma$ (resp. $\Gamma'$) is supposed to be strictly $J$-pseudoconvex
(resp. strictly $J'$-pseudoconvex). Finally, we assume that $f: D
\rightarrow D'$ is a $(J,J')$-biholomorphic map. It follows from the estimates
of the Kobayashi-Royden infinitesimal pseudometric given in \cite{ga-su}
that $f$ extends as a $1/2$-H{\"o}lder
homeomorphism between  $D \cup \Gamma$ and $D' \cup \Gamma'$, such that
$f(\Gamma) = \Gamma'$ and $f(0) = 0$. Finally
$\Gamma$ is defined in a neighborhood of the origin
by the equation $\rho(z) = 0$ where $\rho(z) = 2Re
z^n + 2Re K(z) + H(z) + o(\vert z \vert^2)$ and $K(z) = \sum
K_{\mu\nu}z^{\mu\nu}$, $H(z) = \sum h_{\mu\nu}z^{\mu}\bar
z^{\nu}$, $k_{\mu\nu} = k_{\nu\mu}$, $h_{\mu\nu} = \bar
h_{\nu\mu}$. As we noticed at the end of Section~3 the hypersurface
$\Gamma$ is strictly $\hat{J}$-pseudoconvex at the origin. The hypersurface
$\Gamma'$ admits a similar local representation. In what follows we
assume that we are in this setting. 

Let $\Sigma := \{ z \in \C^n: 2Re z^n + 2Re K('z,0) + H('z,0) < 0\}$,
$\Sigma' := \{ z \in \C^n: 2Re z^n + 2Re K'('z,0) + H'('z,0) < 0\}$.
If $(p^k)$ is a sequence of points in $D$  converging to $0$, then according
to Proposition~\ref{convseq}, the
scaling procedure associates with the pair $(f,(p^k)_k)$ two linear almost
complex structures ${J}_0$ and ${J}'_0$, both defined on $\C^n$,
and a $(J_0,{J}'_0)$-biholomorphism $\hat{f}$ between $\Sigma$ and
$\Sigma'$. Moreover $(\Sigma,J_0)$ and $(\Sigma',J'_0)$ are model 
domains. To prove that the cluster set of the cotangent lift of $f$ at a point
in $N(\Gamma)$ is contained in $N(\Gamma')$, it is sufficient to prove that
$(\partial \hat{f}^n / \partial z^n)('0,-1) \in \mathbb R \backslash \{0\}$.

\vskip 0,1cm
\noindent{\it Step two.} The proof of Proposition~\ref{PPP} is given
by the following Lemma.

\begin{lemma}
\label{cluster}
Let $K$ be a compact subset of the totally real part of the conormal
bundle $\Sigma_J(\partial D)$. Then the cluster set of the cotangent lift
$f^*$ of $f$ on the conormal bundle 
$\Sigma(\partial D)$, when $(z,L)$ tends to $\Sigma_J(\partial D)$
along the wedge $W_U$, is relatively compactly contained 
in the totally real part of $\Sigma(\partial D')$.
\end{lemma}
We recall that the totally real part of $\Sigma(\partial D')$ is
the complement of the zero section in $\Sigma(\partial D')$.

\vskip 0,1cm
\noindent{\sl Proof of Lemma~\ref{cluster}}.  Let $(z^k,L^k)$ be
a sequence in $W_U$ converging to $(0,\partial_J\rho(0)) =
(0,dz^n)$.  We shall prove that the sequence of
linear forms $Q^k := {}^tdf^{-1}(w^k)L^k$, where $w^k = f(z^k)$, converges
to a linear form which up to a {\it real} factor (in view of Part (a)
of Proposition \ref{reality}) coincides with $\partial_{J} \rho(0)=
dz^n$ (we recall that ${}^t$ denotes the transposed map).  It is
sufficient to prove that the $(n-1)$ first component of $Q^k$ with respect to
the dual basis $(\omega_1,\dots,\omega_n)$ of $X$ converge to $0$ and the
last one is bounded below from the origin as $k$ goes to infinity.
The map $X$ being of class $\CC^1$ we can replace $X(0)$ by $X(w^k)$.
Since $(z^k,L^k) \in W_U$, we have $L^k = \omega_n(z^k) +
O(\delta_k)$, where $\delta_k$ is the distance from $z^k$ to the
boundary. Since $\vert\vert\vert df^{-1}_{w^k} \vert\vert\vert =
0(\delta_k^{-1/2})$, we have $Q^k = {}^tdf^{-1}_{w^k}(\omega_n(z^k)) +
O(\delta_k^{1/2})$.  By Proposition~\ref{convseq}, the components of
${}^tdf^{-1}_{w^k}(\omega_n(z^k))$ with respect to the basis
$(\omega_1(z^k),\dots,\omega_n(z^k))$ are the elements of the last line of
the matrix $df^{-1}_{w^k}$ with respect to the basis $X'(w^k)$ and
$X(z^k)$. So its $(n-1)$ first components are $0(\delta_k^{1/2})$ and
converge to $0$ as $k$ tends to infinity. Finally the component $A_{n,n}^k$
is bounded below from the origin by Part (b) of
Proposition~\ref{reality}. \qed

\subsubsection{Compactness principle}
In this section we prove the following

\begin{theorem}\label{wr}
Let $(M,J)$ be an almost complex manifold, not equivalent to a model
domain. Let $D=\{r<0\}$ be a relatively compact domain in a smooth
manifold $N$ and let $(f^\nu)_\nu$ be a sequence of diffeomorphisms
from $M$ to $D$. Assume that

$(i)$ the sequence $(J_\nu:=f^\nu_*(J))_\nu$ extends smoothly up to
$\bar{D}$ and is compact in the $C^2$ convergence on $\bar{D}$,

$(ii)$ the Levi forms of $\partial D$ , $\mathcal L^{J_\nu}(\partial
D)$ are uniformly bounded from below (with respect to $\nu$) by a
positive constant.

Then the sequence $(f^\nu)_\nu$ is compact in the compact-open
topology on $M$.
\end{theorem}

We proceed by contradiction.
Assume that there is a compact $K_0$ in $M$, points $p^\nu \in M$ and a point
$q \in \partial D$ such that $\lim_{\nu \rightarrow \infty}f^\nu(p^\nu) = q$.

\begin{lemma}\label{met-kob}
For every relatively compact neighborhood $V$ of $q$ there is $\nu_0$
such that for $\nu \geq \nu_0$ we have~:
$\lim_{x \rightarrow q}inf_{q' \in D \cap \partial V}d^K_{(D,J_\nu)}=\infty$.
\end{lemma}

\noindent{\it Proof of Lemma~\ref{met-kob}}. Restricting $U$ if
necessary, we may assume that the function $\rho + C \rho ^2$ is a
strictly $J_\nu$-plurisubharmonic function in a neighborhood of
$\bar{D} \cap U$, for sufficiently large $\nu$.
Moreover, using Proposition~B, we can focus on $K_{D \cap U}$. Smoothing
$D \cap U$, we may assume that the hypothesis of Proposition~A are satisfied
on $D \cap U$, uniformly for sufficiently large $\nu$.
In particular, the inequality~(\ref{e3}) is satisfied on
$D \cap U$, with a positive constant $c$ independent of $\nu$.
The result follows by a direct integration of this inequality.
\qed
 
\vskip 0,1cm
The following Lemma is a corollary of Lemma~\ref{met-kob}.
\begin{lemma}\label{lem3.3.1}
For every $K \subset \subset M$ we have : 
$\lim_{\nu \rightarrow \infty}f^\nu(K) =q$.
\end{lemma}

\noindent{\it Proof of Lemma \ref{lem3.3.1}}. Let $K \subset \subset M$ 
be such that $x^0 \in
K$. Since the function $x \mapsto d_D^K(x^0,x)$ is bounded from above by a
constant $C$ on $K$, it follows from the decreasing property of the Kobayashi 
pseudodistance that

\begin{equation}\label{eq2}
d_{(D,J_\nu)}^K(f^\nu(x^0),f^\nu(x)) \leq C
\end{equation}
for every $\nu$ and every 
$x \in K$. It follows from Lemma~\ref{met-kob} that for
every $V \subset \subset U$, containing $p$, we have :
\begin{equation}\label{eq3}
\lim_{\nu \rightarrow \infty}d_{(D,J_\nu)}^K 
(f^\nu(x^0),D \cap \partial V) = +\infty.
\end{equation} 
Then from conditions (\ref{eq2}) and (\ref{eq3}) we deduce that 
$f^\nu(K) \subset V$ for every sufficiently large $\nu$. 
This gives the statement. \qed 

\vskip 0,1cm
Fix now a point $p \in M$ and denote by $p^\nu$ the point $f^\nu(p)$.
We may assume that the sequence $(J_\nu:=f^\nu_*(J))_\nu$ converges
to an almost complex structure $J'$ on $\bar{D}$ and according to
Lemma~\ref{lem3.3.1} we may assume that
$\lim_{\nu \rightarrow \infty}p^\nu = q$.
We apply Subsection~4.3 to the domain $D$ and the sequence $(q^\nu)_\nu$.
We denote by $T^\nu$ the linear transformation
$T^\nu:=M^\nu \circ L^\nu \circ \alpha^\nu$, as in Subsection~4.3, and
we consider $D^\nu:=T^\nu(D)$, and $J^\nu:=T^\nu_*(J_\nu)$.
If $\phi_\nu$ is the nonisotropic dilation $\phi_\nu:('z,z^n) \mapsto
(\delta_\nu^{1/2}\ 'z,\delta_\nu z^n)$ then we set
$\hat{f}^\nu:=\phi_\nu^{-1} \circ T^\nu \circ f$ and
$\hat{J}^\nu:=(\phi_\nu^{-1})_*(J^\nu)$. We also consider
$\hat{\rho}_\nu:=\delta_\nu^{-1} \circ \rho \circ \phi_\nu$
and $\hat{D}^\nu:=\{\hat{\rho}_\nu < 0\}$.
As proved in Subsection~4.3, the sequence $(\hat{D}^\nu)_\nu$ converges,
in the local Hausdorff convergence, to a domain
$\Sigma:=\{z \in C^n:\hat \rho(z) := 2Re z^n + 2Re K({'z},0) + H({'z},0)<0\}$,
where $K$ and $H$ are homogeneous of degree two.
According to Proposition~\ref{convseq} we have~:

$(i)$ The sequence $(\hat{J}^\nu)$ converges to a model almost complex
structure $J_0$, uniformly (with all partial derivatives of any order)
on compact subsets of $\C^n$,

$(ii)$ $(\Sigma,J_0)$ is a model domain,

$(iii)$ the sequence $(\hat{f}^\nu)_\nu$ converges to a $(J,J_0)$
holomorphic map $F$ from $M$ to $\Sigma$.

\vskip 0,1cm
To prove Theorem~\ref{wr}, it remains to prove that $F$ is a diffeomorphism
from $M$ to $\Sigma$. 
We first notice that according to condition $(ii)$ of Theorem~\ref{wr}
and Lemma~\ref{met-kob}, the domain $D$ is
complete $J_\nu$-hyperbolic. In particular, since $f^\nu$ is a $(J,J_\nu)$
biholomorphism from $M$ to $D$, the manifold $M$ is complete $J$-hyperbolic.
Consequently, for every compact subset $L$ of $M$, there is a positive
constant $C$ such that for every $z \in L$ and every $v \in T_zM$ we have
$K_{(M,J)}(z,v) \geq C\|v\|$.
Consider the map $\hat{g}^\nu:=(\hat{f}^\nu)^{-1}$.
This is a $(\hat{J}^\nu,J)$ biholomorphism from $\hat{D}^\nu$ to $M$.
Let $K$ be a compact set in $\Sigma$. We may consider $\hat{g}^\nu(K)$
for sufficiently large $\nu$. By the decreasing property of the Kobayashi
distance, there is a compact subset $L$ in $M$ such that
$\hat{g}^\nu(K) \subset L$ for sufficiently large $\nu$. Then for every
$w \in K$ and for every $v \in T_w\Sigma$ we obtain, by the decreasing of the
Kobayashi-Royden infinitesimal pseudometric~:
$$
\|df^\nu(w)(v)\| \leq (1/C) \|v\|,
$$
uniformly for sufficiently large $\nu$.
According to Ascoli Theorem, we may extract from
$(\hat{g}^\nu)_\nu$ a subsequence, converging to a map $G$ from
$\Sigma$ to $M$. Finally, on any compact subset $K$ of $M$, by
the equality $\hat{g}^\nu \circ \hat{f}^\nu = id$ we obtain $F \circ G = id$.
This gives the result. \qed

\vskip 0,1cm
As a corollary of Theorem~\ref{wr} we obtain the following almost complex
version of the Wong-Rosay Theorem in real dimension four~:
\begin{corollary}\label{wr-2}
Let $(M,J)$ (resp. $(M',J')$) be an almost complex manifold of real dimension
four. Let $D$ (resp. $D'$) be a relatively compact domain in $M$ (resp. $N$).
Consider a sequence $(f^\nu)_\nu$ of diffeomorphisms from $D$ to
$D'$ such that the sequence $(J_\nu:=f^\nu_*(J))_\nu$ extends to $\bar{D}'$
and converges to $J'$ in the $C^2$ convergence on $\bar{D}'$.

Assume that there is a point $p\in D$ and a point $q \in \partial D'$ such
that $\lim_{\nu \rightarrow \infty}f^\nu(p) = q$ and such that
$D'$ is strictly $J'$-pseudoconvex at $q$.
Then there is a $(J,J_{st})$-biholomorphism from $M$ to the unit ball $\B^2$
in $\C^2$. 
\end{corollary}

\noindent{\it Proof of Corollary~\ref{wr-2}}.
The proof of Corollary~\ref{wr-2} follows exactly the same lines as
the proof of Theorem~\ref{wr}. \qed

\section{Elliptic regularity on almost complex manifolds with boundary}

This section is devoted to one of the main technical steps of our
construction. We prove that a pseudoholomorphic disc attached (in the
sense of the cluster set) to a smooth totally real submanifold in an almost
complex manifold, extends  smoothly up to the boundary. In the case of
the integrable structure, various versions of this statement have been
obtained by several authors. In the almost complex case, similar
assertions have been established by H.Hofer \cite{ho}, J.-C.Sikorav
\cite{si}, S.Ivashkovich-V.Shevchishin \cite{iv-sh}, E.Chirka
\cite{ch1}, D.McDuff-D.Salamon~\cite{mc-sa}
under stronger assumptions on the initial boundary regularity of the disc
(at least the continuity is required). Our proof consists of two
steps. First, we show that a disc extends as a $1/2$-H{\"o}lder
continuous map up to the boundary. The proof is based on special
estimates of the Kobayashi-Royden metric in ``Grauert tube'' type
domains. The second step is the reflection principle adapted to the
almost complex category; here we follow the construction of E.Chirka
\cite{ch1}.

%\subsubsection{Elliptic regularity and the reflection}

\subsection{Reflection principle and regularity of analytic discs}

We prove the following~:
\begin{theorem}\label{reflection}
Let $N$ be a smooth $\mathcal C^\infty$ totally real submanifold in $(M,J)$
and let $\varphi : \Delta^+ \rightarrow M$ be $J$-holomorphic, where
$\Delta^+:=\{\zeta \in \Delta : Im(\zeta) >0\}$.
Assume that the cluster set of $\varphi$ on the real interval $]-1,1[$ is
contained in $N$. Then $\varphi$ is of class $\mathcal C^\infty$ on
$\Delta^+ \cup ]-1,1[$.
\end{theorem}

In case $N$ has a weaker regularity then the exact regularity of $\varphi$,
related to that of $N$, can be derived directly from the following proof of
Theorem~\ref{reflection}.

\vskip 0,1cm
\noindent{\it Proof of Theorem~\ref{reflection}.}
\noindent{\it Step one}. It follows by
Theorem~\ref{Regth1} that $\varphi$ extends as a H{\"o}lder 1/2-continuous
map on $\Delta^+ \cup
]-1,1[$.

\vskip 0,1cm
\noindent{\it Step two : The disc $\varphi$ is of class $\mathcal C^{1+1/2}$.}
 The following construction of the reflection principle
for pseudoholomorphic discs is due to Chirka \cite{ch1}. For reader's
convenience we give the details.
Let $a\in ]-1,1[$. Our consideration being local at $a$, we may assume that
$N=\R^n \subset \C^n$, $a=0$ and $J$ is a smooth almost complex structure
defined in the unit ball $\B_n$ in $\C^n$.

After a complex linear change of coordinates we may assume that
$J = J_{st} + O(\vert z \vert)$ and $N$ is given by $x + ih(x)$ where
$x \in \R^n$ and $dh(0) = 0$. If $\Phi$ is the local diffeomorphism
$x \mapsto x$, $y \mapsto y - h(x)$ then $\Phi(N) = \R^n$ and the direct
image of $J$ by $\Phi$, still denoted by $J$, keeps the form $J_{st} +
O(\vert z \vert)$. Then $J$ has a basis of $(1,0)$-forms given in the
coordinates $z$ by $dz^j + \sum_k a_{jk}d\bar z^k$; using the
matrix notation we write it in the form $\omega = dz + A(z)d\bar z$ where
the matrix function $A(z)$ vanishes at the origin. Writing
$\omega = (I + A)dx + i(I - A)dy$ where $I$ denotes the identity
matrix, we can take as a basis of $(1,0)$ forms~: $\omega' = dx +
i(I +A)^{-1}(I - A)dy = dx + iBdy$. Here the matrix function $B$ satisfies
$B(0) = I$. Since $B$ is smooth, its restriction $B_{\vert \R^n}$ on $\R^n$
admits a smooth extension $\hat B$ on the unit ball such that
$\hat B - B_{\vert \R^n} = O(\vert y \vert^k)$ for any positive integer $k$.
Consider the diffeomorphism $z^* = x + i\hat B(z) y$.
In the $z^*$-coordinates the submanifold $N$ still coincides with $\R^n$
and $\omega' = dx + iBdy = dz^* + i(B - \hat B)dy - i(d\hat B)y = dz^* +
\alpha$, where the coefficients of the form $\alpha$ vanish up to
the first order on $\R^n$. Therefore there is a basis of $(1,0)$-forms
(with respect to the image of $J$ under the coordinate diffeomorphism
$z \mapsto z^*$) of the form $dz^* + A(z^*)d\bar z^*$,
where $A$ vanishes to first order on $\R^n$ and
$\| A \|_{\mathcal C^1(\bar{\B}_n)} < < 1$.

Consider the continuous map $\psi$ defined on $\Delta$ by
$$
\left\{
\begin{array}{cccc}
\psi &=& \varphi &{\rm on}\ \Delta^+\\
& & & \\
\psi(\zeta) &=&\overline{\varphi(\bar{\zeta})} &{\rm for}\ \zeta \in
\Delta^- :=\{\zeta \in \Delta / Im(\zeta) < 0\}.
\end{array}
\right.
$$

Since the map $\varphi$ satisfies
\begin{equation}\label{holo}
\bar \partial \varphi + A(\varphi)\overline{\partial \varphi} = 0
\end{equation}
on $\Delta^+$, the map $\psi$ satisfies the equation
$$
\bar\partial\psi(\zeta) + 
\overline{A(\varphi(\bar\zeta))}\
\overline{\partial\psi(\zeta)} = 0
$$
for $\zeta \in \Delta^-$.

Hence $\psi$ is a solution on $\Delta$ of the elliptic equation
\begin{equation}\label{elliptic}
\bar \partial \psi + \lambda(\cdot)\overline{\partial \psi} = 0
\end{equation}
where $\lambda$ is defined by $\lambda(\zeta) =
A(\varphi(\zeta))$ for
$\zeta \in \Delta^+ \cup ]-1,1[$ and $\lambda(\zeta) =
\overline{A(\varphi(\bar\zeta))}$ for
$\zeta \in \Delta^-$.
According to Step
one, the map $\lambda$ is H{\"o}lder $1/2$ continuous on $\Delta$
and vanishes on $]-1,1[$.
This implies that $\psi$ is of class $\mathcal C^{1+1/2}$ on $\Delta$
by equation~(\ref{elliptic}) (see~\cite{si,ve}).

\vskip 0,1cm
\noindent{\it Step three : Geometric bootstrap.} See Figure 7. Let $v=(1,0)$ in
$\R^2$ and
consider the disc $\varphi^c$ defined on $\Delta^+$ by
$$
\varphi^c(\zeta) = (\varphi(\zeta),d\varphi(\zeta)(v)).
$$

\vskip 0,1cm
The cluster set $C(\varphi^c,]-1,1[)$ is contained in the smooth
submanifold $TN$ of $TM$.

\begin{lemma}\label{tot-real}
If $N$ is a totally real submanifold in an almost complex manifold
$(M,J)$ then $TN$ is a totally real submanifold in $(TM,J^c)$.
\end{lemma}

\noindent{\it Proof of Lemma~\ref{tot-real}.} 
Let $X \in T(TN) \cap J^c(T(TN))$. If $X=(u,v)$ in the trivialisation
$T(TM) = TM \oplus TM$ then $u \in TN \cap J(TN)$,
implying that $u=0$. Hence $v  \in TN \cap J(TN)$,
implying that $v=0$. Finally, $X=0$. \qed

\vskip 0,1cm
Applying Step two to $\varphi^c$ and $TN$ we prove that the first derivative
of $\varphi$ with respect to $x$ ($x+iy$ are the standard coordinates on $\C$)
is of class $\mathcal C^{1+1/2}$ on $\Delta^+ \cup ]-1,1[$. 
The $J$-holomorphicity equation~(\ref{holo}) may be written as
$$
\frac{\partial \varphi}{\partial y} = J(\varphi)
\frac{\partial \varphi}{\partial x}
$$
on $\Delta^+ \cup ]-1,1[$.
Hence $\partial \varphi/\partial y$
is of class $\mathcal C^{1+1/2}$ on $\Delta^+ \cup ]-1,1[$, meaning that
$\varphi$ is of class $\mathcal C^{2+1/2}$ on $\Delta^+ \cup ]-1,1[$.
We prove now that $\varphi$ is of class $\mathcal C^{3+1/2}$ on
$\Delta^+ \cup ]-1,1[$. The reader will conclude, repeating the same
argument that $\varphi$ is of class $\mathcal C^\infty$ on
$\Delta^+ \cup ]-1,1[$.

\bigskip
\begin{center}
\input{bootstrap.pstex_t}
\end{center}
\bigskip
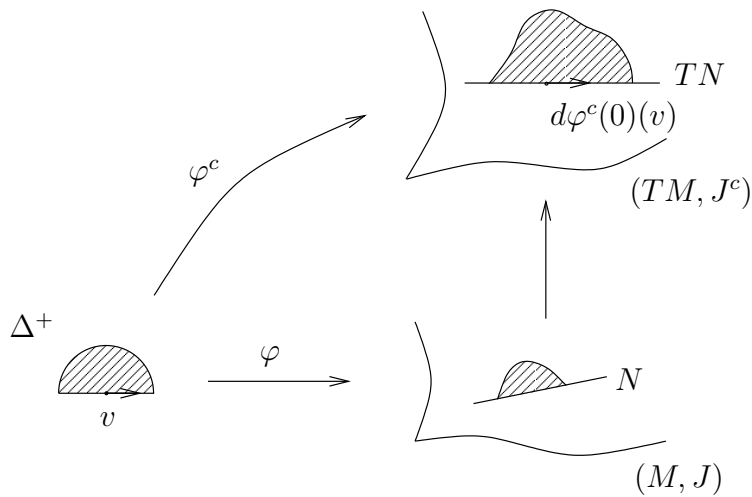
\centerline{Figure 7}
\bigskip

Replace now the data $(M,J)$ and $\varphi$ by $(TM,J^c)$ and
$\varphi^c$ in Step three. The map $^2\varphi^c$ defined on $\Delta^+$ by
$^2\varphi^c(\zeta) = (\varphi^c(\zeta), d\varphi^c(\zeta)(v))$
is $^2J^c$-holomorphic on $\Delta^+$ ($^2J^c$ is the complete lift of $J^c$
to the second tangent bundle $T(TM)$. According to Step two, its
first derivative $\partial (^2\varphi^c)/\partial x$ is of class $C^{1+1/2}$
on $\Delta^+ \cup ]-1,1[$. This means that the second derivatives
$\displaystyle \frac{\partial^2 \varphi}{\partial x^2}$ and
$\displaystyle \frac{\partial^2 \varphi}{\partial x \partial y}$
are $C^{1+1/2}$ on $\Delta^+ \cup ]-1,1[$. Differentiating
equation~(\ref{holo}) with respect to $y$, we prove that
$\displaystyle \frac{\partial^2 \varphi}{\partial y^2}$ is $C^{1+1/2}$ on
$\Delta^+ \cup ]-1,1[$ and so that $\varphi$ is $C^{3+1/2}$ on
$\Delta^+ \cup ]-1,1[$. \qed

\subsection{Behavior of pseudoholomorphic maps near totally real submanifolds}

Let $\Omega$ be a domain in an almost complex manifold $(M,J)$ and $E 
\subset \Omega$ be a smooth $n$-dimensional
totally real submanifold defined as the set of common zeros of the
functions $r_j$, $j=1,...,n$ smooth on $\Omega$. We suppose that
$\bar\partial_J r_1 \wedge ...\wedge \bar\partial_J r_n \neq
0$ on $\Omega$. Consider the ``wedge''
$W(\Omega,E)=\{ z \in \Omega: r_j(z) < 0, j= 1,...,n \}$ with ``edge''
$E$. For $\delta > 0$ we denote by $W_{\delta}(\Omega,E)$ the
``shrinked'' wedge 
$\{ z \in \Omega : r_j(z) - \delta \sum_{k \neq j} r_k <
0, j = 1,..., n \}$.
The main goal of this Section is to prove the following 

\begin{proposition}
\label{Wedges}
Let $W(\Omega,E)$ be a wedge in $\Omega \subset (M,J)$ with a totally real
n-dimensional edge $E$ of class $\CC^{\infty}$ and let $f:W(\Omega,E)
\rightarrow (M',J')$ be a  $(J,J')$-holomorphic map. Suppose that the
cluster set $C(f,E)$ is (compactly) contained in a
$\mathcal C^\infty$ totally real submanifold $E'$ of $M'$.
Then for any $\delta > 0$ the map $f$ extends to
$W_{\delta}(\Omega,E) \cup E$ as a $\CC^{\infty}$-map. 
\end{proposition}

We previously established this statement for a single
$J$-holomorphic disc. The general case also 
relies on the ellipticity of the
$\bar\partial$-operator.  It requires an additional
 technique of attaching pseudoholomorphic discs to a totally real
manifold which could be of independent interest.

Now we prove Proposition~\ref{Wedges}.
Let $(h_t)_t$ be the family of $J$-holomorphic discs, smoothly depending on
the parameter $t \in \R^{2n}$, defined in Lemma~\ref{lem-discs}.
It follows from Lemma~\ref{dlem3.2}, applied to the holomorphic disc $f \circ
h_t$, uniformly with respect to $t$, that there is a constant $C$
such that $\vert \vert \vert df(z) \vert \vert \vert \leq C dist(z,E)^{-1/2}$
for any $z \in W_{\delta}(\Omega,E)$.
This implies that $f$ extends as a
H{\"o}lder $1/2$-continuous map on $W_{\delta}(\Omega,E) \cup E$. 

It follows now from Theorem~\ref{reflection} that every
composition $f \circ h_t$ is smooth up to $\partial \Delta^+$. Moreover,
the direct examination of our argument shows that the $\CC^k$ norm of
the discs $f \circ h_t$ are uniformly
bounded, for any $k$. Recall the separate smoothness principle
(Proposition 3.1, \cite{tu}):

\begin{proposition}\label{separate}
Let $F_j$, $1 \leq j \leq n$, be $\CC^{\alpha}$ ($\alpha > 1$
 noninteger)
 smooth foliations in a domain $\Omega \subset \R^n$ such
that for every point $p \in \Omega$ the tangent vectors to the curves
$\gamma_j \in F_j$ passing through $p$ are linearly independent. Let
$f$ be a function on $\Omega$ such that the restrictions $f
_{\vert{\gamma_j}}$, $1 \leq j \leq n$, are of class
$\CC^{\alpha-1}$ and are uniformly bounded in the $\CC^{\alpha-1}$
norm. Then $f$ is of class $\CC^{\alpha-1}$.
\end{proposition}
Using Lemma~\ref{lem-discs} we construct $n$ transversal
foliations of $E$ by boundaries of Bishop's discs. Since the restriction
of $f$ on every such curve satisfies the hypothesis of
Proposition~\ref{separate}, $f$ is smooth up to $E$. This proves
Proposition~\ref{Wedges}. \qed

\vskip 0,1cm
Let $\Gamma$ and $\Gamma'$ be two totally real maximal submanifolds in almost
complex manifodls $(M,J)$ and $(M',J')$. Let $W(\Gamma,M)$ be a wedge in
$M$ with edge $\Gamma$. 
\begin{proposition}\label{wed-reg}
If $F :W(\Gamma,M) \rightarrow M'$ is $(J,J')$-holomorphic and if the
cluster set of $\Gamma$ is contained in $\Gamma'$ then $F$ extends as a
$\mathcal C^\infty$ map up to $\Gamma$.
\end{proposition}

\noindent{\it Proof of Proposition~\ref{wed-reg}.}
In view of Proposition~\ref{reflection} the proof is classical
(see \cite{co-ga-su}). \qed

\vskip 0,1cm
As a direct application of Proposition~\ref{wed-reg} we obtain the following
partial version of Fefferman's Theorem :
\begin{corollary}\label{feff1}
Let $D$ and $D'$ be two smooth relatively compact domains in real manifolds.
Assume that $D$ admits an almost complex structure $J$ smooth on $\bar D$ and
 such that $(D,J)$ is strictly pseudoconvex. Let $f$ be a smooth
diffeomorphism $f:  D  \rightarrow  D'$, extending as a $\mathcal C^1$
diffeomorphism (still called $f$) between $\bar{D}$ and $\bar{D}'$.
Then $f$ is a smooth $\mathcal C^\infty$
diffeomorphism between $\bar D$ and $\bar D'$  if and only
if the direct image $f_*(J)$ of $J$ under
$f$ extends smoothly on $ \bar D'$ and $(D', f_*(J))$ 
is strictly pseudoconvex.
\end{corollary}

\noindent{\it Proof of Corollary~\ref{feff1}.}
The cotangent lift $f^*$
of $f$ to the cotangent bundle over $D$, locally defined by
$f^*:=(f,^t(df)^{-1})$, is a $(\tilde{J},\tilde{J}')$-biholomorphism
from $T^*D$ to $T^*D'$, where $J':=f_*(J)$.
According to Proposition~\ref{prop-tot-real}, the conormal
bundle $\Sigma(\partial D)$ (resp. $\Sigma(\partial D')$) is a totally real
submanifold in $T^*M$ (resp. $T^*M'$).
We consider $\Sigma(\partial D)$ as the edge of a wedge
$W(\Sigma(\partial D),M)$ contained in $TD$. Then we may apply
Proposition~\ref{wed-reg} to $F=f^*$ to conclude. \qed

\subsection{Fefferman's mapping Theorem}

Here we present one of the main results of our paper. This was
obtained in the paper \cite{ga-su2}.

\begin{theorem}\label{theo-fefferman}
Let $D$ and $D'$ be two smooth relatively compact domains in real
manifolds.  Assume that $D$ admits an almost complex structure $J$
smooth on $\bar D$ and such that $(D,J)$ is strictly
pseudoconvex. Then a smooth diffeomorphism $f: D \rightarrow D'$
extends to a smooth diffeomorphism between $\bar D$ and $\bar D'$ if
and only if the direct image $f_*(J)$ of $J$ under $f$ extends
smoothly on $ \bar D'$ and $(D', f_*(J))$ is strictly pseudoconvex.
\end{theorem}

Theorem~\ref{theo-fefferman} is a consequence of
Proposition~\ref{PPP}.
We recall that according to Proposition~\ref{prop-tot-real}
the conormal bundle $\Sigma_J(\partial D)$ of $\partial D$ is a
totally real submanifold in the cotangent bundle $T^*M$.
Consider the set
$$
S = \{(z,L) \in \R^{2n} \times \R^{2n} :
dist((z,L),\Sigma_J(\partial D)) \leq dist(z,\partial D), z \in D \}.
$$

In a neighborhood $U$ of any totally real point of
$\Sigma_J(\partial D)$, the set S contains a wedge $W_U$ with
$\Sigma_J(\partial D) \cap U$ as totally real edge.

Then in view of Proposition~\ref{wed-reg} we obtain the following
Proposition~:
\begin{proposition}
\label{wedges}
There is a wedge $W_{U'}$ contained in the
wedge $W_U$ such that the map $f^*$ extends to $W_{U'} \cup \Sigma(\partial D)$
as a $\CC^{\infty}$-map. 
\end{proposition}
Proposition~\ref{wedges} implies immediately that $f$ extends
as a smooth $\mathcal C^\infty$ diffeomorphism from $\bar{D}$ to
$\bar{D'}$ (see Figure 8).

\bigskip
\begin{center}
\input{fefferman.pstex_t}
\end{center}
\bigskip
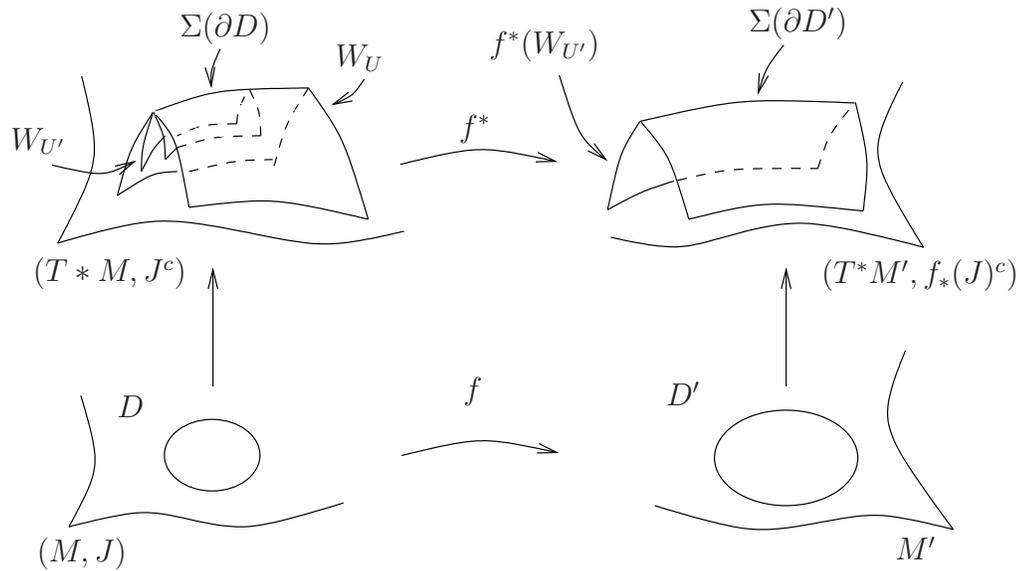
\centerline{Figure 8}

\bigskip
In this survey, we presented an overview of different results
dealing with local analysis in almost complex manifolds, establishing
some bases of the geometry of nonintegrable structures. We
point out that there are many open questions concerning for instance the contact
geometry (the contact properties of the Riemann map,...), the
study of Monge-Amp{\`e}re equations, or the links
between almost complex analysis and symplectic topology. Our approach
here may be considered as a necessary first step to study such
questions.

\end{document}

%% file: figure4.pstex_t
\begin{picture}(0,0)%
\includegraphics{figure4.pstex}%
\end{picture}%
\setlength{\unitlength}{3947sp}%
\begingroup\makeatletter\ifx\SetFigFont\undefined%
\gdef\SetFigFont#1#2#3#4#5{%
  \reset@font\fontsize{#1}{#2pt}%
  \fontfamily{#3}\fontseries{#4}\fontshape{#5}%
  \selectfont}%
\fi\endgroup%
\begin{picture}(5494,1859)(364,-963)
\put(545,-811){\makebox(0,0)[lb]{\smash{\SetFigFont{12}{14.4}{\familydefault}{\mddefault}{\updefault}{\color[rgb]{0,0,0}$(M,J)$}%
}}}
\put(713,689){\makebox(0,0)[lb]{\smash{\SetFigFont{12}{14.4}{\familydefault}{\mddefault}{\updefault}{\color[rgb]{0,0,0}$U \in \mathcal V(q)$}%
}}}
\put(3080,206){\makebox(0,0)[lb]{\smash{\SetFigFont{12}{14.4}{\familydefault}{\mddefault}{\updefault}{\color[rgb]{0,0,0}$z$}%
}}}
\put(1497,-72){\makebox(0,0)[lb]{\smash{\SetFigFont{12}{14.4}{\familydefault}{\mddefault}{\updefault}{\color[rgb]{0,0,0}$q$}%
}}}
\put(4658,-125){\makebox(0,0)[lb]{\smash{\SetFigFont{12}{14.4}{\familydefault}{\mddefault}{\updefault}{\color[rgb]{0,0,0}$0$}%
}}}
\put(5296,-552){\makebox(0,0)[lb]{\smash{\SetFigFont{12}{14.4}{\familydefault}{\mddefault}{\updefault}{\color[rgb]{0,0,0}$\mathbb B_2 \subset \mathbb C^2$}%
}}}
\put(3575,-905){\makebox(0,0)[lb]{\smash{\SetFigFont{12}{14.4}{\familydefault}{\mddefault}{\updefault}{\color[rgb]{0,0,0}$\{z^2 = c_2\}$}%
}}}
\put(5858,547){\makebox(0,0)[lb]{\smash{\SetFigFont{12}{14.4}{\familydefault}{\mddefault}{\updefault}{\color[rgb]{0,0,0}$\{z^1 = c_1\}$}%
}}}
\put(5551,-61){\makebox(0,0)[lb]{\smash{\SetFigFont{12}{14.4}{\familydefault}{\mddefault}{\updefault}{\color[rgb]{0,0,0}$z^1$}%
}}}
\put(4651,764){\makebox(0,0)[lb]{\smash{\SetFigFont{12}{14.4}{\familydefault}{\mddefault}{\updefault}{\color[rgb]{0,0,0}$z^2$}%
}}}
\end{picture}

%% file: figure1.pstex_t
\begin{picture}(0,0)%
\includegraphics{figure1.pstex}%
\end{picture}%
\setlength{\unitlength}{3947sp}%
\begingroup\makeatletter\ifx\SetFigFont\undefined%
\gdef\SetFigFont#1#2#3#4#5{%
  \reset@font\fontsize{#1}{#2pt}%
  \fontfamily{#3}\fontseries{#4}\fontshape{#5}%
  \selectfont}%
\fi\endgroup%
\begin{picture}(3898,3040)(306,-3393)
\put(1315,-2412){\makebox(0,0)[lb]{\smash{\SetFigFont{12}{14.4}{\familydefault}{\mddefault}{\updefault}{\color[rgb]{0,0,0}$f_{p,0}(r\Delta)$}%
}}}
\put(1395,-2028){\makebox(0,0)[lb]{\smash{\SetFigFont{12}{14.4}{\familydefault}{\mddefault}{\updefault}{\color[rgb]{0,0,0}$f_{p,t}(r\Delta)$}%
}}}
\put(306,-2184){\makebox(0,0)[lb]{\smash{\SetFigFont{12}{14.4}{\familydefault}{\mddefault}{\updefault}{\color[rgb]{0,0,0}$K$}%
}}}
\put(3188,-618){\makebox(0,0)[lb]{\smash{\SetFigFont{12}{14.4}{\familydefault}{\mddefault}{\updefault}{\color[rgb]{0,0,0}$\partial M_{11}$}%
}}}
\put(4204,-1724){\makebox(0,0)[lb]{\smash{\SetFigFont{12}{14.4}{\familydefault}{\mddefault}{\updefault}{\color[rgb]{0,0,0}$\{\rho=0\}$}%
}}}
\put(3027,-3334){\makebox(0,0)[lb]{\smash{\SetFigFont{12}{14.4}{\familydefault}{\mddefault}{\updefault}{\color[rgb]{0,0,0}$E=\bar{M}_{11} \cap \{\rho=0\}$}%
}}}
\end{picture}

%% file: diagram3.pstex_t
\begin{picture}(0,0)%
\includegraphics{diagram3.pstex}%
\end{picture}%
\setlength{\unitlength}{3947sp}%
\begingroup\makeatletter\ifx\SetFigFont\undefined%
\gdef\SetFigFont#1#2#3#4#5{%
  \reset@font\fontsize{#1}{#2pt}%
  \fontfamily{#3}\fontseries{#4}\fontshape{#5}%
  \selectfont}%
\fi\endgroup%
\begin{picture}(3300,1814)(3076,-3968)
\put(3526,-3061){\makebox(0,0)[lb]{\smash{\SetFigFont{10}{12.0}{\familydefault}{\mddefault}{\updefault}{$\Psi_J$}%
}}}
\put(3076,-3436){\makebox(0,0)[lb]{\smash{\SetFigFont{10}{12.0}{\familydefault}{\mddefault}{\updefault}{$f_v^ J$}%
}}}
\put(4651,-2311){\makebox(0,0)[lb]{\smash{\SetFigFont{10}{12.0}{\familydefault}{\mddefault}{\updefault}{$d\varphi_0$}%
}}}
\put(5851,-3061){\makebox(0,0)[lb]{\smash{\SetFigFont{10}{12.0}{\familydefault}{\mddefault}{\updefault}{$\Psi_{J'}$}%
}}}
\put(6376,-3511){\makebox(0,0)[lb]{\smash{\SetFigFont{10}{12.0}{\familydefault}{\mddefault}{\updefault}{$f_{d\varphi_0(v)}^ {J'}$}%
}}}
\put(4726,-3811){\makebox(0,0)[lb]{\smash{\SetFigFont{10}{12.0}{\familydefault}{\mddefault}{\updefault}{$\varphi$}%
}}}
\put(6301,-2686){\makebox(0,0)[lb]{\smash{\SetFigFont{10}{12.0}{\familydefault}{\mddefault}{\updefault}{$\zeta d{\varphi_0(v)}$}%
}}}
\put(3151,-2536){\makebox(0,0)[lb]{\smash{\SetFigFont{10}{12.0}{\familydefault}{\mddefault}{\updefault}{$\zeta v$}%
}}}
\end{picture}

%% file: riemann2.pstex_t
\begin{picture}(0,0)%
\includegraphics{riemann2.pstex}%
\end{picture}%
\setlength{\unitlength}{3947sp}%
\begingroup\makeatletter\ifx\SetFigFont\undefined%
\gdef\SetFigFont#1#2#3#4#5{%
  \reset@font\fontsize{#1}{#2pt}%
  \fontfamily{#3}\fontseries{#4}\fontshape{#5}%
  \selectfont}%
\fi\endgroup%
\begin{picture}(3000,1864)(1801,-1319)
\put(4201,239){\makebox(0,0)[lb]{\smash{\SetFigFont{12}{14.4}{\familydefault}{\mddefault}{\updefault}{\color[rgb]{0,0,0}$((\Omega')^\lambda,(J')_\lambda^*)$}%
}}}
\put(4801,-511){\makebox(0,0)[lb]{\smash{\SetFigFont{12}{14.4}{\familydefault}{\mddefault}{\updefault}{\color[rgb]{0,0,0}$\Psi'_{J_\lambda}$}%
}}}
\put(1801,239){\makebox(0,0)[lb]{\smash{\SetFigFont{12}{14.4}{\familydefault}{\mddefault}{\updefault}{\color[rgb]{0,0,0}$(\Omega^\lambda,J_\lambda^*)$}%
}}}
\put(1801,-1261){\makebox(0,0)[lb]{\smash{\SetFigFont{12}{14.4}{\familydefault}{\mddefault}{\updefault}{\color[rgb]{0,0,0}$(\mathbb B_n,J_\lambda)$}%
}}}
\put(3226,-1111){\makebox(0,0)[lb]{\smash{\SetFigFont{12}{14.4}{\familydefault}{\mddefault}{\updefault}{\color[rgb]{0,0,0}$\varphi$}%
}}}
\put(4351,-1261){\makebox(0,0)[lb]{\smash{\SetFigFont{12}{14.4}{\familydefault}{\mddefault}{\updefault}{\color[rgb]{0,0,0}$(\mathbb B_n,J'_\lambda)$}%
}}}
\put(2326,-511){\makebox(0,0)[lb]{\smash{\SetFigFont{12}{14.4}{\familydefault}{\mddefault}{\updefault}{\color[rgb]{0,0,0}$\Psi_{J_\lambda}$}%
}}}
\put(2926,389){\makebox(0,0)[lb]{\smash{\SetFigFont{12}{14.4}{\familydefault}{\mddefault}{\updefault}{\color[rgb]{0,0,0}$L=d\varphi_0$}%
}}}
\end{picture}

%% file: figure2.pstex_t
\begin{picture}(0,0)%
\includegraphics{figure2.pstex}%
\end{picture}%
\setlength{\unitlength}{3947sp}%
\begingroup\makeatletter\ifx\SetFigFont\undefined%
\gdef\SetFigFont#1#2#3#4#5{%
  \reset@font\fontsize{#1}{#2pt}%
  \fontfamily{#3}\fontseries{#4}\fontshape{#5}%
  \selectfont}%
\fi\endgroup%
\begin{picture}(5973,2331)(4003,-3886)
\put(7501,-3886){\makebox(0,0)[lb]{\smash{\SetFigFont{14}{16.8}{\familydefault}{\mddefault}{\updefault}{\color[rgb]{0,0,0} }%
}}}
\put(9976,-3361){\makebox(0,0)[lb]{\smash{\SetFigFont{12}{14.4}{\familydefault}{\mddefault}{\updefault}{\color[rgb]{0,0,0}$D$}%
}}}
\put(6451,-1711){\makebox(0,0)[lb]{\smash{\SetFigFont{12}{14.4}{\familydefault}{\mddefault}{\updefault}{\color[rgb]{0,0,0}$f(s \Delta) \subset D \cap U$}%
}}}
\put(7276,-3811){\makebox(0,0)[lb]{\smash{\SetFigFont{12}{14.4}{\familydefault}{\mddefault}{\updefault}{\color[rgb]{0,0,0}$f(\Delta \setminus s \Delta)$}%
}}}
\put(6151,-2311){\makebox(0,0)[lb]{\smash{\SetFigFont{12}{14.4}{\familydefault}{\mddefault}{\updefault}{\color[rgb]{0,0,0}$f$}%
}}}
\put(7501,-2911){\makebox(0,0)[lb]{\smash{\SetFigFont{12}{14.4}{\familydefault}{\mddefault}{\updefault}{\color[rgb]{0,0,0}$p$}%
}}}
\put(6901,-3361){\makebox(0,0)[lb]{\smash{\SetFigFont{12}{14.4}{\familydefault}{\mddefault}{\updefault}{\color[rgb]{0,0,0}$U$}%
}}}
\end{picture}

%% file: figure3.pstex_t
\begin{picture}(0,0)%
\includegraphics{figure3.pstex}%
\end{picture}%
\setlength{\unitlength}{3947sp}%
\begingroup\makeatletter\ifx\SetFigFont\undefined%
\gdef\SetFigFont#1#2#3#4#5{%
  \reset@font\fontsize{#1}{#2pt}%
  \fontfamily{#3}\fontseries{#4}\fontshape{#5}%
  \selectfont}%
\fi\endgroup%
\begin{picture}(5867,6390)(1426,-6836)
\put(5422,-3742){\makebox(0,0)[lb]{\smash{\SetFigFont{10}{12.0}{\familydefault}{\mddefault}{\updefault}{\color[rgb]{0,0,0}$(0,-1)$}%
}}}
\put(2651,-6009){\makebox(0,0)[lb]{\smash{\SetFigFont{10}{12.0}{\familydefault}{\mddefault}{\updefault}{\color[rgb]{0,0,0}$(0,-1)$}%
}}}
\put(2588,-3805){\makebox(0,0)[lb]{\smash{\SetFigFont{10}{12.0}{\familydefault}{\mddefault}{\updefault}{\color[rgb]{0,0,0}$(0,-1)$}%
}}}
\put(5701,-6061){\makebox(0,0)[lb]{\smash{\SetFigFont{10}{12.0}{\familydefault}{\mddefault}{\updefault}{\color[rgb]{0,0,0}$(0,-1)$}%
}}}
\put(5851,-961){\makebox(0,0)[lb]{\smash{\SetFigFont{10}{12.0}{\familydefault}{\mddefault}{\updefault}{\color[rgb]{0,0,0}$f(p)$}%
}}}
\put(5776,-1636){\makebox(0,0)[lb]{\smash{\SetFigFont{10}{12.0}{\familydefault}{\mddefault}{\updefault}{\color[rgb]{0,0,0}$f^k(\tilde{p}^k)$}%
}}}
\put(1762,-1440){\makebox(0,0)[lb]{\smash{\SetFigFont{10}{12.0}{\familydefault}{\mddefault}{\updefault}{\color[rgb]{0,0,0}$p$}%
}}}
\put(1603,-1860){\makebox(0,0)[lb]{\smash{\SetFigFont{10}{12.0}{\familydefault}{\mddefault}{\updefault}{\color[rgb]{0,0,0}$U$}%
}}}
\put(2551,-1486){\makebox(0,0)[lb]{\smash{\SetFigFont{10}{12.0}{\familydefault}{\mddefault}{\updefault}{\color[rgb]{0,0,0}$\tilde{p}^k$}%
}}}
\put(1426,-586){\makebox(0,0)[lb]{\smash{\SetFigFont{11}{13.2}{\familydefault}{\mddefault}{\updefault}{\color[rgb]{0,0,0}$(D^k,J_k)$}%
}}}
\put(1576,-5461){\makebox(0,0)[lb]{\smash{\SetFigFont{10}{12.0}{\familydefault}{\mddefault}{\updefault}{\color[rgb]{0,0,0}$(\Sigma,J_0)$}%
}}}
\put(4801,-5386){\makebox(0,0)[lb]{\smash{\SetFigFont{10}{12.0}{\familydefault}{\mddefault}{\updefault}{\color[rgb]{0,0,0}$(\Sigma',J'_0)$}%
}}}
\put(4876,-586){\makebox(0,0)[lb]{\smash{\SetFigFont{10}{12.0}{\familydefault}{\mddefault}{\updefault}{\color[rgb]{0,0,0}$((D^k)',J'_k)$}%
}}}
\put(6301,-2611){\makebox(0,0)[lb]{\smash{\SetFigFont{10}{12.0}{\familydefault}{\mddefault}{\updefault}{\color[rgb]{0,0,0}$\psi_k^{-1}$}%
}}}
\put(3001,-2686){\makebox(0,0)[lb]{\smash{\SetFigFont{10}{12.0}{\familydefault}{\mddefault}{\updefault}{\color[rgb]{0,0,0}$\phi_k^{-1}$}%
}}}
\put(3301,-4786){\makebox(0,0)[lb]{\smash{\SetFigFont{10}{12.0}{\familydefault}{\mddefault}{\updefault}{\color[rgb]{0,0,0}$k \rightarrow \infty$}%
}}}
\put(6451,-4786){\makebox(0,0)[lb]{\smash{\SetFigFont{10}{12.0}{\familydefault}{\mddefault}{\updefault}{\color[rgb]{0,0,0}$k \rightarrow \infty$}%
}}}
\put(4426,-961){\makebox(0,0)[lb]{\smash{\SetFigFont{11}{13.2}{\familydefault}{\mddefault}{\updefault}{\color[rgb]{0,0,0}$f^k$}%
}}}
\put(4576,-3061){\makebox(0,0)[lb]{\smash{\SetFigFont{10}{12.0}{\familydefault}{\mddefault}{\updefault}{\color[rgb]{0,0,0}$((\hat{D}')^k,\hat{J}'_k)$}%
}}}
\put(1426,-3061){\makebox(0,0)[lb]{\smash{\SetFigFont{10}{12.0}{\familydefault}{\mddefault}{\updefault}{\color[rgb]{0,0,0}$(\hat{D}^k,\hat{J}_k)$}%
}}}
\end{picture}

%% file: bootstrap.pstex_t
\begin{picture}(0,0)%
\includegraphics{bootstrap.pstex}%
\end{picture}%
\setlength{\unitlength}{3947sp}%
\begingroup\makeatletter\ifx\SetFigFont\undefined%
\gdef\SetFigFont#1#2#3#4#5{%
  \reset@font\fontsize{#1}{#2pt}%
  \fontfamily{#3}\fontseries{#4}\fontshape{#5}%
  \selectfont}%
\fi\endgroup%
\begin{picture}(4174,3069)(293,-2208)
\put(4467,378){\makebox(0,0)[lb]{\smash{\SetFigFont{12}{14.4}{\familydefault}{\mddefault}{\updefault}{\color[rgb]{0,0,0}$TN$}%
}}}
\put(4185,-350){\makebox(0,0)[lb]{\smash{\SetFigFont{12}{14.4}{\familydefault}{\mddefault}{\updefault}{\color[rgb]{0,0,0}$(TM,J^c)$}%
}}}
\put(1865,-1355){\makebox(0,0)[lb]{\smash{\SetFigFont{12}{14.4}{\familydefault}{\mddefault}{\updefault}{\color[rgb]{0,0,0}$\varphi$}%
}}}
\put(293,-1209){\makebox(0,0)[lb]{\smash{\SetFigFont{12}{14.4}{\familydefault}{\mddefault}{\updefault}{\color[rgb]{0,0,0}$\Delta^+$}%
}}}
\put(860,-1771){\makebox(0,0)[lb]{\smash{\SetFigFont{12}{14.4}{\familydefault}{\mddefault}{\updefault}{\color[rgb]{0,0,0}$v$}%
}}}
\put(1416,-222){\makebox(0,0)[lb]{\smash{\SetFigFont{12}{14.4}{\familydefault}{\mddefault}{\updefault}{\color[rgb]{0,0,0}$\varphi^c$}%
}}}
\put(3687,123){\makebox(0,0)[lb]{\smash{\SetFigFont{12}{14.4}{\familydefault}{\mddefault}{\updefault}{\color[rgb]{0,0,0}$d\varphi^c(0)(v)$}%
}}}
\put(4213,-2150){\makebox(0,0)[lb]{\smash{\SetFigFont{12}{14.4}{\familydefault}{\mddefault}{\updefault}{\color[rgb]{0,0,0}$(M,J)$}%
}}}
\put(4096,-1543){\makebox(0,0)[lb]{\smash{\SetFigFont{12}{14.4}{\familydefault}{\mddefault}{\updefault}{\color[rgb]{0,0,0}$N$}%
}}}
\end{picture}

%% file: fefferman.pstex_t
\begin{picture}(0,0)%
\includegraphics{fefferman.pstex}%
\end{picture}%
\setlength{\unitlength}{3947sp}%
\begingroup\makeatletter\ifx\SetFigFont\undefined%
\gdef\SetFigFont#1#2#3#4#5{%
  \reset@font\fontsize{#1}{#2pt}%
  \fontfamily{#3}\fontseries{#4}\fontshape{#5}%
  \selectfont}%
\fi\endgroup%
\begin{picture}(5932,3549)(-144,-2269)
\put(526,-1306){\makebox(0,0)[lb]{\smash{\SetFigFont{12}{14.4}{\familydefault}{\mddefault}{\updefault}$D$}}}
\put(2675,460){\makebox(0,0)[lb]{\smash{\SetFigFont{12}{14.4}{\familydefault}{\mddefault}{\updefault}$f^*$}}}
\put(4947,-462){\makebox(0,0)[lb]{\smash{\SetFigFont{12}{14.4}{\familydefault}{\mddefault}{\updefault}$(T^*M',f_*(J)^c)$}}}
\put(2701,-1191){\makebox(0,0)[lb]{\smash{\SetFigFont{12}{14.4}{\familydefault}{\mddefault}{\updefault}$f$}}}
\put(-144,379){\makebox(0,0)[lb]{\smash{\SetFigFont{12}{14.4}{\familydefault}{\mddefault}{\updefault}$W_{U'}$}}}
\put(3971,-1226){\makebox(0,0)[lb]{\smash{\SetFigFont{12}{14.4}{\familydefault}{\mddefault}{\updefault}$D'$}}}
\put(5411,-2201){\makebox(0,0)[lb]{\smash{\SetFigFont{12}{14.4}{\familydefault}{\mddefault}{\updefault}$M'$}}}
\put( -4,-448){\makebox(0,0)[lb]{\smash{\SetFigFont{12}{14.4}{\familydefault}{\mddefault}{\updefault}$(T*M,J^c)$}}}
\put(918,1088){\makebox(0,0)[lb]{\smash{\SetFigFont{12}{14.4}{\familydefault}{\mddefault}{\updefault}$\Sigma(\partial D)$}}}
\put(1896,904){\makebox(0,0)[lb]{\smash{\SetFigFont{12}{14.4}{\familydefault}{\mddefault}{\updefault}$W_U$}}}
\put(2861,984){\makebox(0,0)[lb]{\smash{\SetFigFont{12}{14.4}{\familydefault}{\mddefault}{\updefault}$f^*(W_{U'})$}}}
\put(4491,1124){\makebox(0,0)[lb]{\smash{\SetFigFont{12}{14.4}{\familydefault}{\mddefault}{\updefault}$\Sigma(\partial D')$}}}
\put( 21,-2211){\makebox(0,0)[lb]{\smash{\SetFigFont{12}{14.4}{\familydefault}{\mddefault}{\updefault}$(M,J)$}}}
\end{picture}